\documentclass[a4paper,10pt]{article}
\usepackage[utf8]{inputenc}  

\usepackage{amsthm}
\usepackage{amssymb}
\usepackage{pifont}
\usepackage{amsmath}
\usepackage{mathrsfs}
\usepackage{dsfont}
\usepackage{yhmath}
\usepackage{listings}
\usepackage{enumitem}
\usepackage{float}
\usepackage{graphicx}
\usepackage{changepage}
\usepackage{listings}
\usepackage{adjustbox}

\usepackage{pgf,tikz}
\usetikzlibrary{arrows}
\usetikzlibrary{patterns}

\usepackage[all]{xy}

\usepackage{subcaption}
\usepackage{ifthen}

\usepackage[colorlinks=true,linkcolor=blue,citecolor=green]{hyperref}

\newcommand{\R}{\mathbb{R}}
\newcommand{\N}{\mathbb{N}}
\newcommand{\Z}{\mathbb{Z}}
\newcommand{\Q}{\mathbb{Q}}
\newcommand{\C}{\mathbb{C}}

\newcommand{\F}{\mathbb{F}}

\renewcommand{\O}{{\cal O}}

\newcommand{\GL}{\mathrm{{GL}}}
\newcommand{\SL}{\mathrm{{SL}}}

\newcommand{\Stab}{\mathrm{Stab}}

\def\lint{[\hskip -1.5pt[}
\def\rint{]\hskip -1.5pt]}

\def\bsys{\left\{\begin{array}}
\def\esys{\end{array}\right.}

\newcommand{\eps}{\varepsilon}

\newcommand{\cmark}{\textcolor{green}{\text{\ding{51}}}}
\newcommand{\xmark}{\text{\ding{55}}}

\renewcommand{\1}{\mathbf{1}}

\newcommand{\grp}[1]{\langle #1 \rangle}
\newcommand{\diag}{\mathrm{diag}}
\newcommand{\id}{\mathrm{id}}
\newcommand{\fix}{\mathrm{Fix}}
\newcommand{\pstab}{\mathrm{PStab}}

\newcommand{\step}[2]{\medskip \noindent \underline{{\it Step #1:}} #2}
\newcommand{\steppf}[3]{\medskip \noindent \underline{{\it Step #1:}} #2
\begin{proof} #3
\end{proof}}

\theoremstyle{plain}
\newtheorem{theorem}{Theorem}[section]
\newtheorem{lemma}[theorem]{Lemma}
\newtheorem{minilemma}[theorem]{Mini-lemma}
\newtheorem{proposition}[theorem]{Proposition}
\newtheorem{corollary}[theorem]{Corollary}

\theoremstyle{definition}
\newtheorem{definition}[theorem]{Definition}
\newtheorem{notation}[theorem]{Notation}

\newtheorem{remark}[theorem]{Remark}
\newtheorem{example}[theorem]{Example}
\newtheorem{nexample}[theorem]{Non-exemple}

\newcommand{\lem}[2][None]
{\begin{lemma}\ifthenelse{\equal{#1}{None}}{}{\label{#1}}
#2
\end{lemma}}

\newcommand{\minlem}[2][None]
{\begin{minilemma}\ifthenelse{\equal{#1}{None}}{}{\label{#1}}
#2
\end{minilemma}}

\newcommand{\prop}[2][None]
{\begin{proposition}\ifthenelse{\equal{#1}{None}}{}{\label{#1}}
#2
\end{proposition}}

\newcommand{\cor}[2][None]
{\begin{corollary}\ifthenelse{\equal{#1}{None}}{}{\label{#1}}
#2
\end{corollary}}

\newcommand{\theo}[2][None]
{\begin{theorem}\ifthenelse{\equal{#1}{None}}{}{\label{#1}}
#2
\end{theorem}}

\newcommand{\defi}[2][None]
{\begin{definition}\ifthenelse{\equal{#1}{None}}{}{\label{#1}}
#2
\end{definition}}

\newcommand{\nota}[2][None]
{\begin{notation}\ifthenelse{\equal{#1}{None}}{}{\label{#1}}
#2
\end{notation}}

\newcommand{\rem}[2][None]
{\begin{remark}\ifthenelse{\equal{#1}{None}}{}{\label{#1}}
#2
\end{remark}}

\newcommand{\expl}[2][None]
{\begin{example}\ifthenelse{\equal{#1}{None}}{}{\label{#1}}
#2
\end{example}}

\newcommand{\nexpl}[2][None]
{\begin{nexample}\ifthenelse{\equal{#1}{None}}{}{\label{#1}}
#2
\end{nexample}}

\newcommand{\demo}[1]{\begin{proof} #1
\end{proof}}

\newcommand{\demode}[2]{\begin{proof}[Proof of #1] #2
\end{proof}}

\title{Finite quotients of abelian varieties with a Calabi-Yau resolution}
\author{C\'ecile Gachet\footnote{Universit\'e C\^ote d'Azur, CNRS, LJAD, France}}

\begin{document}

\maketitle

\begin{abstract}
Let $A$ be an abelian variety, and $G\subset\mathrm{Aut}(A)$ a finite group acting freely in codimension two. We discuss whether the singular quotient $A/G$ admits a resolution that is a Calabi-Yau manifold. While Oguiso constructed two examples in dimension 3, we show that there are none in dimension 4. We also classify up to isogeny the possible abelian varieties $A$ in arbitrary dimension.
\end{abstract}

\section{Introduction}

Since singularities are a byproduct of the Minimal Model Program, studying singular varieties with trivial canonical class, or singular $K$-trivial varieties, is an important question in the birational classification of complex algebraic varieties. From this point of view, the recent generalization of the Beauville-Bogomolov decomposition theorem for smooth $K$-trivial varieties (\cite{Beauv84}) to klt $K$-trivial varieties (\cite{GGK,Druel,HoeringPeternell,BGL}) is highly relevant. It indeed establishes that, after a finite quasiétale cover, any klt $K$-trivial variety is a product of a smooth abelian variety, some irreducible holomorphic symplectic varieties with canonical singularities, also called hyperkähler varieties, and some Calabi-Yau varieties with canonical singularities. These three main families of $K$-trivial varieties are the subject of large, mostly disjoint realms of the literature, ranging from the well-known theory of abelian varieties (exposed notably in the reference books \cite{BirkLang,Shimura}), through the thriving study of hyperkähler varieties (see \cite{DebarreHK,AmVerbSurvey,Kamenova} for surveys),
to the unruly ``Zoo of Calabi-Yau varieties'', populated by a huge amount of examples (\cite{KreuzerSkarke1998,KreuzerSkarke2000} for K3 surfaces and Calabi-Yau threefolds embedded as hypersurfaces in toric varieties only), and whose boundedness is yet not established (see \cite{WilsonPic2, WilsonPic3, ChenDiCerboHanSvaldi,DiCerboSvaldi,BirkarDiCerboSvaldi} for recent breakthroughs).

A new feature appearing in the context of singular $K$-trivial varieties is that birational $K$-trivial varieties may have different Beauville-Bogomolov decomposition types. For example, Kummer surfaces are $K3$ surfaces, but also minimal resolutions of finite quasiétale quotients of abelian surfaces. Similar examples of dimension 3 exist, as in \cite{III0}, as well as higher dimensional examples, {\it cf.} \cite{CynkHulek,CynkSchuett,PR,AndreaWis,Burek}.
In arbitrary dimension, it is known that a crepant resolution or terminalization only changes the type of a klt $K$-trivial variety if its decomposition entails an abelian factor (\cite[Prop.4.10]{Druel}).\footnote{It reflects the more general fact that the Beauville-Bogomolov decomposition type of a klt $K$-trivial variety $X$ with non-trivial fundamental group $\pi_1(X_{\mathrm{reg}})$ is not captured by its algebra of global holomorphic differential forms $H^{0}\left(X,{\Omega_X}^{[\cdot]}\right)$. Many examples supporting this fact are exposed in \cite[Sec.14]{GGK}, and most notably, smooth $K$-trivial threefolds with Beauville-Bogomolov decomposition of pure abelian type and algebra of global differential forms generated by the volume form (as for a Calabi-Yau threefold) are classified in \cite{OgSak}.}

This paper aims at describing changes of the type of a $K$-trivial variety through a birational morphism in the extremal case, i.e., when a singular variety with Beauville-Bogomolov decomposition of purely abelian type is resolved by a Calabi-Yau manifold. We work in the following set-up: 
By a Calabi-Yau manifold, we mean a smooth simply-connected complex projective variety of dimension $n$ with trivial canonical bundle, without any global holomorphic differential form of degree $1\le i\le n-1$. Extending the terminology of \cite{OgFibSp}, we define $n$-dimensional Calabi-Yau manifolds of type $n_0$ as follows.

\theo[thm-LuTaji]{\textup{\cite{SBW},\cite[Rem.1.5]{LuTaji}} Let $X$ be a Calabi-Yau manifold of dimension $n$. The following are equivalent:
\begin{enumerate}[label = \textup{(\roman*)}]
\item There is a nef and big divisor $D$ on $X$ such that $c_2(X)\cdot D^{n-2}=0$.
\item There is an abelian variety $A$ and a finite group $G$ acting freely in codimension $2$ on $A$ such that $X$ is a crepant resolution of $A/G$.
\end{enumerate}
If it satisfies these conditions, $X$ is called a Calabi-Yau manifold of type $n_0$.}

Calabi-Yau threefolds of type $\textup{III}_0$ appear naturally when classifying extremal contractions of Calabi-Yau threefolds \cite{OgFibSp}, and fit in a more general circle of ideas on how the cubic intersection form and the second Chern class determine the birational geometry of a Calabi-Yau threefold (see, e.g., the work of Wilson \cite{WilKahCone}, Oguiso and Peternell \cite{OgPet}). Calabi-Yau threefolds of type $\textup{III}_0$ were classified by Oguiso, as we now recall.

\theo[thm-Oguiso]{\textup{\cite{III0}} There are exactly two Calabi-Yau threefolds $X_3$, $X_7$ of type $\textup{III}_0$. They are the unique crepant resolution of ${E_j}^3$ quotiented by the group generated by $j\id_3$, and of ${E_{u_7}}^3$ quotiented by the group generated by: 
$$\left(\begin{array}{ccc}
0 & -8 & 7-10u_7\\
1 & -6-2u_7 & 11-u_7\\
0 & -1-2u_7 & 6+3u_7\\
\end{array}\right).$$
where
$j=e^{2i\pi/3}$,
$\zeta_7=e^{2i\pi/7}$,
$u_7=\zeta_7+{\zeta_7}^2+{\zeta_7}^4=\frac{-1+i\sqrt{7}}{2}$, and for
any complex number $z\in\C\setminus\R$, we denote by $E_z$ the elliptic curve $\C/(\Z\oplus z\Z)$.}

Our first theorem restricts the isogeny type of $A$ in arbitrary dimension.

\theo[theo-babymain]{Let $A$ be an abelian variety of dimension $n$ and $G$ be a finite group acting freely in codimension $2$ on $A$. If $A/G$ has a resolution that is a Calabi-Yau manifold, then $A$ is isogenous to ${E_j}^n$ or to ${E_{u_7}}^n$ and $G$ is generated by its elements that admit fixed points in $A$.}

Moreover, the local geometry of $A/G$ is generally quite similar to the $3$-dimensional model (see Theorem \ref{theo-main} below). Two important consequences of this are the following results.

\theo[theo-freecod3]{Let $A$ be an abelian variety and $G$ be a finite group acting freely in codimension $3$ on $A$. Then the quotient $A/G$ has a resolution $X$ with $c_1(X)=0$ if and only if $G$ acts freely on $A$.}

\theo[theo-dim4]{Let $A$ be an abelian fourfold and $G$ be a finite group acting freely in codimension $2$ on $A$, yet not freely on A. If $A/G$ has a resolution $X$ with $c_1(X)=0$, then there is a finite étale Galois cover $\tilde{X}$ of $X$ that is a product of the form $E\times X_3$ or $E\times X_7$, where $E$ is an elliptic curve and $X_3$, $X_7$ are the two Calabi-Yau threefolds of $\textup{III}_0$ type.}

In particular, there are no new, irreducible examples arising in dimension 4. We know of no irreducible example in dimension $n\ge 5$ either, and the discussion at the end of the introduction leads us to conjecture that there are none. The importance of quotients that have a simply-connected, Calabi-Yau resolution within the wider range of quotients that have a resolution $X$ with $c_1(X)=0$ is highlighted by a decomposition result in spirit of the Beauville-Bogomolov decomposition, see Lemma \ref{lem_bbc2}.

\medskip

On the way, we make a substantial use of a necessary criterion due to Yamagishi \cite[Theorem 1.1]{Yama} for certain quotient singularities to admit crepant resolutions. This criterion can be stated as follows.

\prop[prop-yama]{Let $G\subset\SL_n(\C)$ be a finite group acting freely in codimension 1 on $\C^n$, and let $U\subset \C^n$ be a $G$-stable simply-connected analytic neighborhood of $0\in\C^n$. If the singularity $U/G$ admits a crepant resolution $X_G$, then the group $G$ is generated by junior elements.}

It is worth noting that this criterion is a weaker version of the following conjecture by Ito and Reid \cite[4.5]{ItoReid}.

\medskip

\noindent{\bf Conjecture.} Let $G\subset\SL_n(\C)$ be a finite group acting freely in codimension 1 on $\C^n$, and let $U\subset \C^n$ be a $G$-stable simply-connected analytic neighborhood of $0\in\C^n$. If the singularity $U/G$ admits a crepant resolution $X_G$, then every maximal cyclic subgroup of $G$ is generated by a junior element.

\medskip

On our way to proving Theorem \ref{theo-main} below, we look for more local arguments on whether particular, explicit finite quotient singularities admit a crepant resolution. This leads us to find a counterexample to the conjecture \cite[4.5]{ItoReid} above, which is strong evidence that the necessary criterion \cite[Theorem 1.1]{Yama}, albeit far from sufficient (see, e.g., \cite[Example 3.11]{Sato}), would be hard to improve. Our counterexample features a faithful representation of the group $\SL_2(\F_3)\simeq Q_8\rtimes\Z_3$ into $\SL_6(\C)$, and is presented in more detail in Remark \ref{rem-macaulay}.

\bigskip

The structure of the paper is as follows. Sections \ref{sec-McKay} to \ref{sec-pftheomain} build up to the proof of the main technical result.

\theo[theo-main]{Let $A$ be an abelian variety of dimension $n$ and $G$ be a finite group acting freely in codimension $2$ on $A$. If $A/G$ has a resolution that is a Calabi-Yau manifold, then
\begin{enumerate}[label = \textup{(\arabic*)}]
\item $A$ is isogenous to ${E_j}^n$ or to ${E_{u_7}}^n$, and $G$ is generated by its elements that admit fixed points in $A$.
\item For every translated abelian subvariety $W\subset A$, there is a integer $k\ge 0$ such that the {\it pointwise stabilizer}
$$ \pstab(W):=\{g\in G\mid \forall w\in W,\, g(w)=w\}$$ 
is isomorphic to ${\Z_3}^k$ if $A$ is isogenous to ${E_j}^n$, or to ${\Z_7}^k$ if $A$ is isogenous to ${E_{u_7}}^n$.
\item For every translated abelian subvariety $W\subset A$, if $ \pstab(W)$ is isomorphic to
\begin{itemize}
\item ${\Z_3}^k$ for some $k\ge 1$, then there are $k$ generators of it such that their matrices are similar to $\diag(\1_{n-3},j,j,j)$, and the $j$-eigenspaces of these matrices are in direct sum. 
\item ${\Z_7}^k$ for some $k\ge 1$, then there are $k$ generators of it such that their matrices are similar to $\diag(\1_{n-3},\zeta_7,{\zeta_7}^2,{\zeta_7}^4)$, and all eigenspaces of these matrices with eigenvalues other than 1 are in direct sum.
\end{itemize}
\end{enumerate}
}

Our starting point in Section \ref{sec-McKay} is to recall the necessary condition \cite[Theorem 1.1]{Yama}, phrased above as Proposition \ref{prop-yama} for a local quotient singularity to admit a crepant resolution. This criterion sheds light on the importance of the so-called junior elements of $\SL_n(\C)$ when studying the birational geometry of finite quotient singularities in dimension $n$. A junior element is a matrix $M\in\SL_n(\C)$ of finite order $d$, whose eigenvalues $(e^{2i\pi a_k/d})_{1\le k\le n}$ with $0\le a_k\le d-1$ satisfy $\sum a_k = d$.

Matrices inducing actions on abelian varieties satisfy a rationality requirement \cite[Proposition 1.2.3]{BirkLang}, which translates into arithmetic constraints on their characteristic polynomial. These constraints allow us to classify matrices of junior elements $g$ acting on $n$-dimensional abelian varieties up to similarity: In Section \ref{sec-JuniorClassification}, we prove that if a junior element $g$ acts on an abelian variety in a way that the generated group $\grp{g}$ acts freely in codimension $2$, then the matrix of $g$ is of one of twelve possible types (see Proposition \ref{prop-junclass}). In particular, the order of $g$ and the number of non-trivial eigenvalues of $g$ are bounded independently of the dimension $n$.

\medskip

The next step is to show that ten out of the twelve types of junior elements can not belong to $G$, for a mix of local and global reasons. The proof spreads throughout Sections \ref{sec-CyclicPointwiseStabilizer}, \ref{sec-RuleOutMostJuniors}, \ref{sec-codim5} and \ref{sec-codim6}. Let us sketch the idea of the argument in the simplest case, namely if $g$ is a junior element of composite order other than 6, with at most four non-trivial eigenvalues. If such a junior element $g$ belongs to $G$, then some non-trivial power $g^{\alpha}$ is not junior, and has a larger fixed locus in $A$. Fix an irreducible component $W$ of that larger fixed locus that is not in the fixed locus of $g$: the pointwise stabilizer $\pstab(W)\subset G$ does not contain $g$, but the power $g^{\alpha}$. Now, as $W$ has codimension less than 4, Section \ref{sec-CyclicPointwiseStabilizer} shows that $\pstab(W)$ is cyclic generated by one junior element $h$, and thus, up to possibly replacing $h$ by another junior generator of $\fix(W)$, one has $g^{\alpha}=h^{\alpha}$. For well-chosen $\alpha$, this is enough to yield $g=h$, and a contradiction. 

This idea excludes seven out of the twelve types of junior elements (see Subsection \ref{subsec-RuleOutSevenJuniors}). The three types of junior elements of order 6 are excluded by technical variations in the next sections. Ruling them out works along with classifying pointwise stablizers in higher codimension: In codimension 4, Section \ref{sec-CyclicPointwiseStabilizer} establishes cyclicity of the pointwise stabilizers and Section \ref{sec-RuleOutMostJuniors} deduces that junior elements with four non-trivial eigenvalues do not exist; in codimension 5 (Section \ref{sec-codim5}), we first prove that junior elements with five non-trivial eigenvalues do not exist (Subsection \ref{subsec-ruleout65}), then deduce cyclicity of the pointwise stabilizers (Subsection \ref{subsec-EigenspaceDirectSumCod5}). In codimension 6 (Section \ref{sec-codim6}), we first classify pointwise stabilizers which do not contain junior elements of type $\diag(\1_{n-6},\omega,\omega,\omega,\omega,\omega,\omega)$: they are isomorphic to $\Z_3,\Z_7,\Z_3\times\Z_3,\Z_7\times\Z_7,$ or $\SL_2(\F_3)$ (Subsection \ref{subsec-pstabcod6}). We use this partial classification to rule out junior elements with six non-trivial eigenvalues (Subsection \ref{subsec-ruleout66}), and we then finally refine the study of pointwise stabilizers in codimension 6 by ruling out $\SL_2(\F_3)$ (Subsection \ref{subsec-JuniorsCommute}).

There finally remain two types of possible junior elements, which are those already appearing in dimension 3 in \cite{III0}: $\diag(\1_{n-3},j,j,j)$ and $\diag(\1_{n-3},\zeta_7,{\zeta_7}^2,{\zeta_7}^4)$.

\medskip

This description of pointwise stabilizers in codimension up to 6 implies that any two junior elements admitting a common fixed point commute. Together with a simple argument about the isogeny type of $A$ (see Section \ref{sec-MakeDisjoint}), it concludes the proof of Theorem \ref{theo-main}. In fact, the idea that the existence of certain automorphisms on an abelian variety determines the isomorphism type of some special abelian subvarieties is general (\cite{Shimura}), and it applies crucially throughout this paper, starting in Section \ref{sec-CyclicPointwiseStabilizer}. From there, it is not so surprising that we are able to determine the isogeny type of $A$, interpreting the fact that $A/G$ admits a Calabi-Yau resolution as an irreducibility property of the $G$-equivariant Poincaré decomposition of $A$.

\medskip

Under the additional assumption that the group $G$ is abelian, Theorem \ref{theo-main} and the results of Section \ref{sec-MakeDisjoint} suffice to generalize Theorem \ref{theo-dim4} to higher dimensions, i.e., to the statement that, if $A$ is an abelian variety of dimension $n$ and $G$ is a finite group acting freely in codimension $2$ on $A$ such that $A/G$ admits a Calabi-Yau resolution $X$, then $n=3$ and $X$ is $X_3$ or $X_7$.

Also note that $G$ is abelian if and only if any two junior elements $g,h$ of $G$ commute, which by our results can be checked {\it via} their matrices acting on a vector space $V$ of dimension $3$, $4$, $5$, or $6$. Standard finite group theory allows us to explicitly bound the order of $\grp{g,h}$ depending on this dimension and the isogeny type of $A$. If the dimension is 3 or 4, the bounds are reasonable enough to launch a computer-assisted search through all possible abstract groups $\grp{g,h}$. Among these, the only groups which, in a faithful 3 or 4-dimensional representation, are generated by two junior elements of the same type, are $\Z_3$, $\Z_7$, and the finite simple group $\SL_3(\F_2)$ of order 168. But a geometric argument on fixed loci excludes $\SL_3(\F_2)$, whence the wished contradiction.  This reproves the classification of \cite{III0} in dimension 3, and settles Theorem \ref{theo-dim4}.

When $V$ has dimension 5 or 6, we could also bound the order of $\grp{g,h}$ explicitly. For example, we could consider the image of the faithful representation $M\oplus\overline{M}$ in $\SL_{2\dim(V)}(\Q)$, and use the classification of irreducible maximal finite integral matrix groups in dimension less than 12 by V. Felsch, G. Nebe, W. Plesken, and B. Souvignier to obtain a bound on the order of $\grp{g,h}$. But the bounds obtained in this way are too large for the \texttt{SmallGroup} library. One needs to better understand the arising matrix groups of larger order, and build a reasonably smaller finite list of possibilities for the abstract group $\grp{g,h}$. It will then remain to figure out geometric ways for ruling out those potential groups in the list other than $\Z_3$, $\Z_7$, $\Z_3\times\Z_3$, and $\Z_7\times\Z_7$.

\medskip

Some of our proofs resort to computer-searches among all finite groups of certain fixed orders (relying on the \texttt{SmallGroup} library of \texttt{GAP}). The computer-assisted results used in Subsection \ref{subsec-GroupThPointwiseStabilizer} were actually originally proven by hand using elementary representation theory and Sylow theory. Such arguments being standard in finite group theory, we chose to keep their exposition concise for the sake of readability, and preferred invoking computer-checked facts as black boxes when needed. This approach also has the advantage of better separating abstract group-theoretic arguments on $G$ from properties of the particular representation $G\hookrightarrow\GL(H^0(T_A))$. All programs used are available in the Appendix.

\medskip

\medskip

\noindent \textbf{Acknowledgments.} I am grateful to my advisor Andreas Höring for fruitful discussions. This paper has greatly benefit from an anonymous referee's report, and I want to thank them for their time, care, commitment, and for bringing the highly relevant results of \cite{Yama} to my attention. Finally, I thank Julia Schneider for suggesting to state and prove Theorem \ref{theo-freecod3}, and Stéphane Druel, for proposing the current, better statement of Theorem \ref{theo-freecod3}.


\tableofcontents

\newpage

\section{Some results in McKay correspondence}\label{sec-McKay}

\defi{Let $g$ be a matrix in $\GL_n(\C)$. Assume that it has finite order $d$. Since $g^d = \id$, $g$ is diagonalizable and has eigenvalues of the form $e^{2i\pi a_k/d}$, for integers $0\le a_1\le\ldots\le a_n\le d-1$. (We allow redundancies.)

\noindent We define the {\it ranked vector of eigenvalues} of $g$ as the tuple $(e^{2i\pi a_k/d})_{1\le k\le n}$.

\noindent The {\it age} of $g$ is set to be the number $\frac{a_1+\ldots+a_n}{d}.$ If it equals 1, we say that $g$ is {\it junior}.}

\defi[def-juniorab]{If $A$ is an abelian variety of dimension $n$ and $g\in\mathrm{Aut}(A)$ has finite order, then $g$ can be written as:
$$g:[z]\in A\mapsto [M(g)z+T(g)]\in A,$$
where $M(g)$ is a matrix of finite order in $\GL_n(\C)$, $T(g)$ a vector in $\C^n$. If $g$ fixes any point $a$ of $A$, it can be represented locally in a neighborhood of $a$ by its matrix $M(g)$. Hence, it makes sense to say that the automorphism $g$ is {\it junior} if $g$ fixes at least one point in $A$ and the matrix $M(g)$ is junior.}

\rem{Note that if $g\in\mathrm{Aut}(A)$ admits a fixed point, then $\grp{g}$ contains no translation, so $g$ and its matrix $M(g)$ have the same order.}

Junior elements play a key role in the study of finite quotient singularities, as the following theorem emphasizes.

\theo[thm-IR]{\textup{\cite{ItoReid}} Let $G$ be a finite subgroup of $\SL_n(\C)$. Suppose that the finite Gorenstein quotient $\C^n/G$ has a minimal model $X$. Then there is a natural one-to-one correspondence between conjugacy classes of junior elements in $G$ and prime exceptional divisors in $X$.}

\rem{Note that such a minimal model $X$ always exists as a relative minimal model of a resolution $\tilde{X}\to\C^n/G$, by \cite[1.30.6]{Singbook}.}

Quotient singularities are $\Q$-factorial, so they can not be resolved by small birational morphisms. This yields a simple corollary of the theorem.

\cor[cor-IR]{\textup{\cite{ItoReid}} Let $G$ be a finite subgroup of $\SL_n(\C)$. If the Gorenstein quotient singularity $\C^n/G$ admits a crepant resolution, then there is a junior element $g\in G$.}

In fact, \cite[Par.4.5]{ItoReid} conjectures that under the same hypotheses, if the singularity $\C^n/G$ admits a crepant resolution, then any maximal cyclic subgroup of $G$ contains a junior element. A counterexample to this conjecture is however presented in Remark \ref{rem-macaulay}.

We will also use a slightly stronger version of Corollary \ref{cor-IR} for quotients by finite subgroups of $\GL_n(\C)$. The formulation is inspired by \cite[Theorem 2.3]{MorrisonStevens}.

\prop[prop-Reid]{\textup{\cite{Reid,RSB,Tai}} Let $G$ be a finite subgroup of $\GL_n(\C)$. If the Gorenstein quotient singularity $\C^n/G$ admits a crepant resolution, then there is a junior element $g\in G$.}

\subsection{A local necessary criterion}

The following result, due to \cite{Yama}, provides a valuable necessary criterion for the existence of a crepant resolution to an affine quotient singularity, in the spirit of \cite{ItoReid}.

\prop[prop-genjun]{Let $G\subset\SL_n(\C)$ be a finite group acting freely in codimension 1 on $\C^n$, 
If the affine quotient singularity $\C^n/G$ admits a crepant resolution, then the group $G$ is generated by junior elements.}


The following corollary allows us to apply this proposition to some finite groups which are {\it a priori} contained in $\GL_n(\C)$ and not $\SL_n(\C)$. We recall that a {\it $K$-trivial variety} is a smooth projective variety $X$ with $K_X\sim\O_X$.

\cor[cor-pstabjun]{Let $G$ be a finite group acting freely in codimension 1 on an abelian variety $A$. Suppose that $A/G$ has a $K$-trivial resolution $X$. Then for every point $a\in A$, the stabilizer
$\Stab(a):=\{g\in G\mid g(a)=a\}$
is generated by junior elements.}

\demo{
Since $X$ is $K$-trivial, we have $h^n(X,\O_X)=1$, where $n$ denotes the dimension of $X$. Since $A/G$ has rational singularities, this yields $h^n(A/G,\O_{A/G})=1$, hence there exists a non-zero $G$-invariant element in $H^n(A,\O_A)$. Since $A$ is an abelian variety, the action of $G$ on $H^n(A,\O_A)=\bigwedge^n H^1(A,\O_A)$ is given by the dual of the determinant representation $\det M: G\to \GL_n(\C)\to \C^*$. The fact that there is a non-zero invariant element shows that $\det M$ is trivial, i.e., $M(G)\subset \SL_n(\C)$.

Fix $a\in A$. Note that, near the point $a$, the variety $A/G$ is locally analytically isomorphic to the affine quotient singularity $\C^n/\Stab(a)$ near the origin $0$ -- with the action by $M:\Stab(a)\to\SL_n(\C)$. Since $A/G$ admits a crepant resolution, the affine quotient singularity $\C^n/\Stab(a)$ admits a crepant resolution too.
Applying Proposition \ref{prop-genjun} yields that the subgroup $M(\Stab(a))<\SL_n(\C)$ is generated by junior elements, and this concludes.
}

\subsection{A global result along the same lines}

We also prove a global result along the same lines as Proposition \ref{prop-genjun}.

\lem[lem-makeitglob]{Let $G$ be a finite group acting freely in codimension 1 on an abelian variety $A$. Suppose that $A/G$ has a resolution $X_G$ with $c_1(X_G)=0$ and $\pi_1(X_G)=\{1\}$ Then $G$ is generated by its elements admitting fixed points in $A$.}

\demo{Let $G_0\triangleleft G$ be the normal subgroup of $G$ generated by elements admitting fixed points. We want to prove that $G_0 = G$. We have a commutative diagram given by the fiber product:
\[
\xymatrix{
X_0\ar[d]^{\eps_0}\ar[r]^{\tilde{q}} & X_G\ar[d]^{\eps_G} \\
A/G_0\ar[r]^{q} & A/G
}
\]

By definition of $G_0$, for every $a\in A$, the stabilizers of $a$ in $G$ and $G_0$ coincide. Hence, $q$ is étale, and $\tilde{q}$ is étale too by base change.
But $X_G$ is simply-connected and $X_0$ is connected, so $\mathrm{deg}(\tilde{q})=1$ and $G_0 = G$.}

\rem{If $G$ is a finite group acting freely in codimension 1 on an abelian variety $A$ so that $A/G$ has a resolution $X$ with $c_1(X)=0$ and $\pi_1(X)=\{1\}$, then $G$ may still contain elements that admit no fixed point. Without loss of generality, we can assume that $G$ contains no translation, up to replacing $A$ by an isogenous abelian variety, but that is the best we can do.}

\subsection{A reduction result inspired by the Beauville-Bogomolov decomposition}

We conclude this section with a reduction result in the spirit of the Beauville--Bogomolov decomposition. We reiterate that we include simply connectedness in the definition of a Calabi--Yau variety.

\lem[lem_bbc2]{Let $X$ be a smooth projective variety with $c_1(X)=0$. Assume that $X$ is a resolution of a quotient $A/G$, where $A$ is an abelian variety, and $G$ is a finite group acting freely in codimension $2$ on $A$. Then, there is a finite étale Galois cover $p:\tilde{X}\to X$ such that
$$\tilde{X}= B\times\prod_{i=1}^r Y_i,$$ 
where $B$ is an abelian variety, and each $Y_i$ is a Calabi-Yau variety of dimension at least 3 that resolves a quotient $B_i/H_i$, where $B_i$ is an abelian variety and $H_i$ is a finite group acting freely in codimension $2$ on $B_i$.}

\demo{Let $n$ denote the dimension of $X$. Pulling back an ample divisor on $A/G$ to $X$, we have a nef and big Cartier divisor $D$ on $X$ such that $c_2(X)\cdot D^{n-2}=0$. By the Beauville--Bogomolov decomposition theorem \cite{Beauv84}, we can take a finite étale Galois cover $p:\tilde{X}\to X$ that decomposes as
$$\tilde{X}= B\times\prod_{i=1}^r Y_i\times \prod_{j=1}^s Z_j,$$ where $B$ is an abelian variety, each $Y_i$ is a Calabi-Yau variety of dimension at least 3, and each $Z_j$ is an irreducible holomorphic symplectic variety. Since $h^1(Y_i,\O_{Y_i})=0$ for each $i$ and $h^1(Z_j,\O_{Z_j})=0$ for each $j$, we can use \cite[Exercise III.12.6]{HartBook} to write 
$$p^*D=p_B^*D_B+\sum_{i=1}^r p_i^*L_i+\sum_{j=1}^s q_j^*M_j,$$
where $D_B$, $L_i$, $M_j$ are divisors on $B$, $Y_i$, $Z_j$ respectively, and $p_B,p_i,q_j$ are the projections onto $B$, $Y_i$, $Z_j$.
By the projection formula and since $D^n > 0$, it is easy ot check that each of the divisors $D_B$, $L_i$, $M_j$ is nef and big.

As the tangent bundle of $\tilde{X}$ decomposes into a direct sum of the tangent bundles of $B$, of the $Y_i$, and of the $Z_j$, we also note that
$$c_2(\tilde{X}) = \sum_{i=1}^r p_i^*c_2(Y_i) + \sum_{j=1}^s q_j^*c_2(Z_j).$$
Since the second Chern class of a $K$-trivial variety of dimension $m$ has non-negative intersection number with the product of $(m-2)$ nef divisors \cite{Miyaoka}, and since $c_2(X)\cdot D^{n-2}=0$, we obtain
$c_2(Y_i)\cdot L_i^{\dim Y_i - 2} = 0$ for all $i$ and 
$c_2(Z_j)\cdot M_j^{\dim Z_j - 2} = 0$ for all $j$. 
By \cite[Rem.1.5]{LuTaji}, \cite{SBW}, this shows that each $Y_i$ resolves a quotient $B_i/H_i$, where $B_i$ is an abelian variety and $H_i$ is a finite group acting freely in codimension $2$ on $B_i$, and that each $Z_j$ resolves a quotient $C_j/K_j$, where $C_j$ is an abelian variety and $K_j$ is a finite group acting freely in codimension $2$ on $C_j$.

To conclude, we note that each of the finite quotients $C_j/K_j$ has a symplectic resolution $Z_j$ and is smooth in codimension $2$. Hence \cite[Thm, Cor.1]{NamiSymplectic} applies, and $C_j/K_j$ has terminal singularities. As a terminal $\Q$-factorial variety that is $K$-trivial, it cannot have a $K$-trivial resolution. So there are no irreducible holomorphic symplectic factors $Z_j$ in $\tilde{X}$, and this concludes this proof.}

\section{The twelve types of junior elements on an abelian variety}\label{sec-JuniorClassification}

Section \ref{sec-McKay} just shows that, if we want a finite singular quotient of an abelian variety $A/G$ to have a resolution $X$ with $c_1(X)=0$, the group $G$ must contain some junior elements. The fact that in our set-up, $G$ must also act freely in codimension $2$ on $A$ is restrictive enough that there are only twelve possibilities for the ranked vector of eigenvalues of a junior element $g\in G$.

\prop[prop-junclass]{Let $A$ be an abelian variety of dimension $n$, and $g\in\mathrm{Aut}(A)$ be a junior element such that $\grp{g}$ acts freely in codimension $2$. Then the order $d$ of $g$ and the ranked vector of eigenvalues of $g$ are in one of the twelve columns of Table \ref{tab-propjunclass}.}

\begin{table}[H]
\begin{adjustbox}{width=\columnwidth,center}
$\begin{array}{|c|c|c|c|}
\hline 
d & 3 & 4 & 6 \\
\hline
(e^{2i\pi a_k/d})
&\left(\1_{n-3},j,j,j\right)
&\left(\1_{n-4},i,i,i,i\right)
&\left(\1_{n-4},\omega,\omega,\omega,-1\right)\\
\hline
\hline
d & 6 & 6 & 7 \\
\hline
(e^{2i\pi a_k/d})
&\left(\1_{n-5},\omega,\omega,\omega,\omega,j\right)
&\left(\1_{n-6},\omega,\omega,\omega,\omega,\omega,\omega\right)
& \left(\1_{n-3},\zeta_7,{\zeta_7}^2,{\zeta_7}^4\right)\\
\hline
\hline
d & 8 & 12 & 15 \\
\hline
(e^{2i\pi a_k/d})
& \left(\1_{n-4},\zeta_8,\zeta_8,\zeta_8^3,\zeta_8^3\right)
& \left(\1_{n-4},\zeta_{12},\zeta_{12},\zeta_{12}^5,\zeta_{12}^5\right)
&\left(\1_{n-4},\zeta_{15},\zeta_{15}^2,\zeta_{15}^4,\zeta_{15}^8\right)\\
\hline
\hline 
d & 16 & 20 & 24\\
\hline
(e^{2i\pi a_k/d})
& \left(\1_{n-4},\zeta_{16},\zeta_{16}^3,\zeta_{16}^5,\zeta_{16}^7\right)
& \left(\1_{n-4},\zeta_{20},\zeta_{20}^3,\zeta_{20}^7,\zeta_{20}^9\right)
& \left(\1_{n-4},\zeta_{24},\zeta_{24}^5,\zeta_{24}^7,\zeta_{24}^{11}\right)\\
\hline
\end{array}$
\end{adjustbox}
\caption{Possible ranked vectors of eigenvalues for junior elements in $G$}
\label{tab-propjunclass}
For $d\in\N$, we denote $\zeta_d =e^{2i\pi/d}$, and in particular $j=e^{2i\pi/3}$ and $\omega=e^{2i\pi/6}$. For $k\in\N$, $\1_k$ refers to a sequence of $k$ times the symbol $1$ in a row.
\end{table}

The proof goes by elementary arithmetic and meticulous case disjunctions. The following terminology should simplify the exposition.

\defi{A {\it multiset} $\mathbf{A}$ is the data of a set $A$ and a function $m: A\to \Z_{>0}$, called the {\it multiplicity function}. Intuitively, a multiset is like a set where elements are allowed to appear more than once.

If a multiset $\mathbf{A}=(A,m)$ is {\it finite}, i.e., its underlying set $A=\{a_1,\ldots,a_k\}$ is finite, we may write $A$ in the following form:
$$\{\{\underset{m(a_1)\mbox{\footnotesize{ times}}}{\underbrace{a_1,\ldots,a_1}},\ldots,\underset{m(a_k)\mbox{\footnotesize{ times}}}{\underbrace{a_k,\ldots,a_k}}\}\}.$$
Double-braces are used to avoid confusion between the multiset and the underlying set.

Let ${\bf A}=(A,m)$ be a finite multiset.

\noindent If $\alpha\in\Z_{>0}$, we denote by ${\bf A}^{*\alpha}$ the multiset $(A, \alpha m)$. 

\noindent If $A$ a subset of $\Q$, and $p,q$ are rational numbers, with $q\ne 0$, we denote by $p+q{\bf A}$ the multiset $(p+qA, m')$, where the multiplicity function $m':p+qA\to\Z_{>0}$ is obtained composing the affine base change $x\in p+qA\mapsto \frac{x-p}{q}\in A$ with $m$.

\noindent The {\it cardinal} of ${\bf A}$ is: 
$$|{\bf A}|:=\displaystyle\sum_{a\in{\bf A}}m(a).$$

\noindent More generally, if $f:A\to \Q$ is a function, we define:
$$\displaystyle\sum_{a\in {\bf A}} f(a):=\displaystyle\sum_{a\in{A}} m(a)f(a).$$

\noindent If ${\bf A}=(A,m)$ and ${\bf B}=(B,n)$ are two multisets, we define their {\it union}:
$${\bf A}\cup{\bf B} := (A\cup B, \1_Am+\1_Bn),$$
where $\1_A,\,\1_B$ are the indicator functions of $A$ and $B$.
}

\nota{For $d\in\N$, we denote by $\Phi_d$ the $d$-th cyclotomic polynomial, and by $\phi(d)$ the degree of $\Phi_d$. In other terms, $\phi$ is the Euler indicator function. 
 
\noindent For integers $a,b$, the greatest common divisor of $a$ and $b$ is denoted $a\wedge b$.}

We establish a sequence of three useful lemmas.

\lem[lem-inteq]{Let $u$ be a positive integer strictly greater than 2. Then we have:
$$\left(2\nmid u\mbox{ and }\frac{\phi(u)^2}{u}\le 8\right)
\mbox{ or }
\left(2\mid u\mbox{ and }\frac{\phi(u)^2}{u}\le 4\right)$$
$$\Leftrightarrow\;u\in\lint 3, 10\rint\cup \{12,14,15,16,18,20,21,24,30,36,42\}.$$}

\demo{Write $u=p_1^{\alpha_1}\cdot p_2^{\alpha_2}\cdots p_k^{\alpha_k},$
where $p_1<\ldots <p_k$ are prime numbers, and $\alpha_1,\ldots,\alpha_k$ positive integers, so that:
$$\frac{\phi(u)^2}{u}=\prod_{i=1}^k (p_i-1)^2p_i^{\alpha_i-2}.$$
Each of the $k$ factors of this product is greater or equal to 1, unless $p_1^{\alpha_1}=2$ in which case the first factor is $\frac{1}{2}$.

Hence, if $u$ satisfies: $$\frac{\phi(u)^2}{u}\le 8\mbox{ or }
\left(2\mid u\mbox{ and }\frac{\phi(u)^2}{u}\le 4\right),$$
then each factor satisfies: 
\begin{equation}\label{eq-factineq}
(p_i-1)^2p_i^{\alpha_i-2}\le 8,
\end{equation}
which yields $p_i\in\{2,3,5,7\}$.
Writing $u=2^\alpha 3^{\beta} 5^{\gamma} 7^{\delta},$
where $\alpha,\beta,\gamma,\delta\ge 0$ and using Inequality (\ref{eq-factineq}) again bounds
$\alpha\le 4,\;\beta\le 2,\;\gamma\le 1,\;\delta\le 1.$
Among the finitely many possibilities left, it is easy to check that the solutions exactly are $u\in\lint 3, 10\rint\cup \{12,14,15,16,18,20,21,24,30,36,42\}$.}

\lem[lem-betterinteq]{Let $u\ge 2$ and $d\ge 3$ be integers, such that $u$ divides $d$. Suppose that there are a positive integer $\alpha$ and a multiset ${\bf A}$ such that:
$$ {\bf A}\cup (d-{\bf A}) = \left\{\left\{a\in\lint 1,d-1\rint\mid u=\frac{d}{d\wedge a}\right\}\right\}^{*\alpha},$$
and such that the quantity:
$$S_{{\bf A},d}(u) := \sum_{a\in {\bf A}} \frac{a}{u(a\wedge d)}$$
satisfies $S_{{\bf A},d}(u)\le 1$.
Then $u,\, \frac{1}{d}{\bf A},\, \alpha,\, S_{{\bf A},d}(u)$ are classified in Table \ref{tab-betterinteq}.
}

\newpage

\begin{table}
\centering
\renewcommand{\arraystretch}{1.4}
$\begin{array}{|c|c|c|c|}
\hline 
u & \alpha & \frac{1}{d}A & S_{A,d}(u)\le 1 \\ 
\hline
2 & 1 & \left\{\frac{1}{2}\right\} & \frac{1}{2} \\
\cline{2-4}
& 2 & \left\{\frac{1}{2},\frac{1}{2}\right\} & 1 \\
\hline 
 & 1 & \left\{\frac{1}{3}\right\}, \left\{\frac{2}{3}\right\} & \frac{1}{3},\frac{2}{3} \\
\cline{2-4}
3 & 2 & \left\{\frac{1}{3},\frac{1}{3}\right\}, \left\{\frac{1}{3},\frac{2}{3}\right\} & \frac{2}{3},1\\
 \cline{2-4}
 & 3 & \left\{\frac{1}{3},\frac{1}{3},\frac{1}{3}\right\} & 1 \\
\hline 
 & 1 & \left\{\frac{1}{4}\right\},\left\{\frac{1}{4}\right\} & \frac{1}{4},\frac{3}{4} \\ 
\cline{2-4}
4 & 2 & \left\{\frac{1}{4},\frac{1}{4}\right\},\left\{\frac{1}{4},\frac{3}{4}\right\} & \frac{1}{2},1\\
\cline{2-4}
& 3 & \left\{\frac{1}{4},\frac{1}{4},\frac{1}{4}\right\} & \frac{3}{4}\\
\cline{2-4}
& 4 & \left\{\frac{1}{4},\frac{1}{4},\frac{1}{4},\frac{1}{4}\right\} & 1\\
\hline 
5 & 1 & \left\{\frac{1}{5},\frac{2}{5}\right\},\left\{\frac{1}{5},\frac{3}{5}\right\} & \frac{3}{5},\frac{4}{5} \\ 
\hline 
 & 1 & \left\{\frac{1}{6}\right\}, \left\{\frac{5}{6}\right\} & \frac{1}{6},\frac{5}{6} \\ 
\cline{2-4}
& 2 & \left\{\frac{1}{6},\frac{1}{6}\right\}, \left\{\frac{1}{6},\frac{5}{6}\right\} & \frac{1}{3},1\\
\cline{2-4}
6 & 3 & \left\{\frac{1}{6},\frac{1}{6},\frac{1}{6}\right\} &\frac{1}{2}\\
\cline{2-4}
& 4 & \left\{\frac{1}{6},\frac{1}{6},\frac{1}{6},\frac{1}{6}\right\} &\frac{2}{3}\\
\cline{2-4}
& 5 & \left\{\frac{1}{6},\frac{1}{6},\frac{1}{6},\frac{1}{6},\frac{1}{6}\right\} &\frac{5}{6}\\
\cline{2-4}
& 6 & \left\{\frac{1}{6},\frac{1}{6},\frac{1}{6},\frac{1}{6},\frac{1}{6},\frac{1}{6}\right\} &1\\
\hline 
7 & 1 & \left\{\frac{1}{7},\frac{2}{7},\frac{3}{7}\right\},\left\{\frac{1}{7},\frac{2}{7},\frac{4}{7}\right\} & \frac{6}{7},1 \\ 
\hline
8 & 1 & \left\{\frac{1}{8},\frac{3}{8}\right\},\left\{\frac{1}{8},\frac{5}{8}\right\} & \frac{1}{2},\frac{3}{4} \\ 
\cline{2-4}
& 2 & \left\{\frac{1}{8},\frac{1}{8},\frac{3}{8},\frac{3}{8}\right\} & 1\\
\hline 
9 & 1 & \left\{\frac{1}{9},\frac{2}{9},\frac{4}{9}\right\},\left\{\frac{1}{9},\frac{2}{9},\frac{5}{9}\right\} & \frac{7}{9},\frac{8}{9} \\ 
\hline
10 & 1 & \left\{\frac{1}{10},\frac{3}{10}\right\},\left\{\frac{1}{10},\frac{7}{10}\right\} & \frac{2}{5},\frac{4}{5}\\ 
\cline{2-4}
& 2& \left\{\frac{1}{10},\frac{1}{10},\frac{3}{10},\frac{3}{10}\right\} & \frac{4}{5}\\
\hline 
12 & 1 & \left\{\frac{1}{12},\frac{5}{12}\right\},\left\{\frac{1}{12},\frac{7}{12}\right\} & \frac{1}{2},\frac{2}{3}\\
\cline{2-4}
& 2  & \left\{\frac{1}{12},\frac{1}{12},\frac{5}{12},\frac{5}{12}\right\} & 1 \\
\hline   
14 & 1 & \left\{\frac{1}{14},\frac{3}{14},\frac{5}{14}\right\},\left\{\frac{1}{14},\frac{3}{14},\frac{9}{14}\right\} & \frac{9}{14},\frac{13}{14}\\ 
\hline 
15 & 1 & \left\{\frac{1}{15},\frac{2}{15},\frac{4}{15},\frac{7}{15}\right\},\left\{\frac{1}{15},\frac{2}{15},\frac{4}{15},\frac{8}{15}\right\} & \frac{14}{15}, 1 \\ 
\hline 
16 & 1& \left\{\frac{1}{16},\frac{3}{16},\frac{5}{16},\frac{7}{16}\right\} & 1 \\ 
\hline 
18 & 1& \left\{\frac{1}{18},\frac{5}{18},\frac{7}{18}\right\},\left\{\frac{1}{18},\frac{5}{18},\frac{11}{18}\right\} & \frac{13}{18},\frac{17}{18}\\ 
\hline 
20 & 1 & \left\{\frac{1}{20},\frac{3}{20},\frac{7}{20},\frac{9}{20}\right\}  & 1 \\ 
\hline 
24 & 1 & \left\{\frac{1}{24},\frac{5}{24},\frac{7}{24},\frac{11}{24}\right\} & 1 \\ 
\hline 
\end{array}$
\renewcommand{\arraystretch}{1}
\caption{Possibilities for $u,\,\frac{1}{d}A,\,\alpha,S_{A,d}(u)$ such that $S_{A,d}(u)\le 1$}
\label{tab-betterinteq}

\end{table}

\newpage

\demo{We start by noting that the underlying sets satisfy 
$$A\cup (d-A)=\left\lbrace a\in\lint 1,d-1\rint\mid u=\frac{d}{d\wedge a}\right\rbrace,$$
and are in an order-preserving bijection with $L:=\{\ell\in\lint 1,u-1\rint\mid \ell\wedge u = 1\}$ {\it via}
\[\arraycolsep=1.2pt\def\arraystretch{1.7}
\begin{array}{lllll}
&f: & a\in A\cup (d-A) &\mapsto \frac{a}{a\wedge d}=\frac{ua}{d} &\in L,\\
\mbox{and }&g: & \ell\in L  &\mapsto \frac{d\ell}{u} &\in A\cup (d-A).
\end{array}
\]
So $|A|\ge \frac{\phi(u)}{2}$. Since $f$ is injective, the restriction $f|_A$ takes at least $\frac{\phi(u)}{2}$ distinct values in its image set inside $L$, so that:

\begin{equation}\label{sadeq}
1\ge S_{{\bf A},d}(u) = \frac{1}{u}\sum_{a\in{\bf A}}f(a) \ge \frac{\alpha}{u}\left(\sum_{\substack{1\le \ell\le u/2 \\ \ell\wedge u = 1}} \ell\right).
\end{equation}

\noindent Let us denote by $\Sigma(u)$ the sum $\displaystyle\sum_{\substack{1\le \ell\le u/2 \\ \ell\wedge u = 1}} \ell.$ 
We have the following coarse estimates:
$$u\ge \Sigma(u) \ge \displaystyle\sum_{\ell=1}^{\phi(u)/2} \ell \ge \frac{\phi(u)^2}{8},
\mbox{ and, if $u$ is even, }u\ge \Sigma(u) \ge \displaystyle\sum_{\ell=1}^{\phi(u)/2} (2\ell-1) \ge \frac{\phi(u)^2}{4}.$$

\noindent Applying Lemma \ref{lem-inteq}, these coarse estimates yield finitely many possibilities for $u$. Computing explicitly $\frac{1}{u}\Sigma(u)$ for the possible values and applying Inequality \eqref{sadeq} again, we exclude a few of them, finally obtaining that:
$$u\in\lint 2, 10\rint\cup \{12,14,15,16,18,20,24\}.$$ 

For each $u$, we then list by hand the finitely many possibilities for the multiplicity $\alpha$ and the multiset $\frac{1}{d}{\bf A}$, and this is how we construct Table \ref{tab-betterinteq}.}

\lem[lem-partsad]{Let $k\in\N$. For each $m\in\lint 1,k\rint$, let $u_m\ge 2$ and $d_m\ge 3$ be integers, such that $u_m$ divides $d_m$, and suppose that there are a positive integer $\alpha_m$ and a multiset ${\bf A_m}$ such that:
$$ {\bf A_m}\cup (d_m-{\bf A_m}) = \left\{\left\{a\in\lint 1,d_m-1\rint\mid u_m=\frac{d_m}{d_m\wedge a}\right\}\right\}^{*\alpha_m}.$$
Suppose additionally that:
$$\sum_{m=1}^k S_{{\bf A_m},d_m}(u_m) = 1.$$
Then the data of $k$ and of all $u_m,\,\alpha_m,\,\frac{1}{d_m}{\bf A_m}$ is classified in Table \ref{tab-partitions}.}

\demo{It is easily derived by hand from Table \ref{tab-betterinteq}.}

\newpage

\begin{table}
\begin{adjustbox}{scale=0.98,center}
\renewcommand{\arraystretch}{1.4}
$\begin{array}{|c|c|c|c|}
\hline 
u_1,\ldots,u_k & \alpha_1,\ldots,\alpha_k & \frac{1}{d_1}A_1,\ldots,\frac{1}{d_k}A_k & \mbox{freeness in codimension $2$}\\ 
\hline
2 & 2 & \left\{\frac{1}{2},\frac{1}{2}\right\} & \xmark \\
\hline
2,3,6 & 1,1,1 & \left\{\frac{1}{2}\right\},\left\{\frac{1}{3}\right\},\left\{\frac{1}{6}\right\}& \xmark\\
\hline
2,4 & 1,2 & \left\{\frac{1}{2}\right\},\left\{\frac{1}{4},\frac{1}{4}\right\}& \xmark\\
\hline
2,6& 1,3 & \left\{\frac{1}{2}\right\},\left\{\frac{1}{6},\frac{1}{6},\frac{1}{6}\right\} & \cmark\\
\hline
2,8& 1,1 & \left\{\frac{1}{2}\right\},\left\{\frac{1}{8},\frac{3}{8}\right\}& \xmark\\
\hline 
2,12& 1,1 & \left\{\frac{1}{2}\right\},\left\{\frac{1}{12},\frac{5}{12}\right\}& \xmark\\
\hline
3 & 2 & \left\{\frac{1}{3},\frac{2}{3}\right\} & \xmark\\
\cline{2-4}
 & 3 & \left\{\frac{1}{3},\frac{1}{3},\frac{1}{3}\right\} & \cmark\\
\hline
3,4,6 & 1,2,1 & \left\{\frac{1}{3}\right\},\left\{\frac{1}{4},\frac{1}{4}\right\},\left\{\frac{1}{6}\right\}& \xmark\\
\hline
3,6 & 1,2 & \left\{\frac{2}{3}\right\},\left\{\frac{1}{6},\frac{1}{6}\right\} & \xmark\\
\cline{2-4}
& 1,4 & \left\{\frac{1}{3}\right\},\left\{\frac{1}{6},\frac{1}{6},\frac{1}{6},\frac{1}{6}\right\} & \cmark \\
\cline{2-4}
 & 2, 2 & \left\{\frac{1}{3},\frac{1}{3}\right\},\left\{\frac{1}{6},\frac{1}{6}\right\}& \xmark\\
\hline
3,12 & 1,1 & \left\{\frac{1}{3}\right\},\left\{\frac{1}{12},\frac{7}{12}\right\}& \xmark\\
\hline
 4 & 2 & \left\{\frac{1}{4},\frac{3}{4}\right\}& \xmark\\
 \cline{2-4}
 & 4 & \left\{\frac{1}{4},\frac{1}{4},\frac{1}{4},\frac{1}{4}\right\} & \cmark\\
 \hline
 4,6 & 2,3 & \left\{\frac{1}{4},\frac{1}{4}\right\},\left\{\frac{1}{6},\frac{1}{6},\frac{1}{6}\right\} & \xmark\\
 \hline
 4,8 & 1, 1 & \left\{\frac{1}{4}\right\},\left\{\frac{1}{8},\frac{5}{8}\right\}& \xmark\\
 \cline{2-4}
 & 2,1 & \left\{\frac{1}{4},\frac{1}{4}\right\},\left\{\frac{1}{8},\frac{3}{8}\right\}& \xmark\\
 \hline
 4, 12 & 2,1 & \left\{\frac{1}{4},\frac{1}{4}\right\},\left\{\frac{1}{12},\frac{5}{12}\right\}& \xmark\\
\hline 
 5, 10 & 1, 1 & \left\{\frac{1}{5},\frac{2}{5}\right\},\left\{\frac{1}{10},\frac{3}{10}\right\} & \xmark\\
\hline 
 6 & 2 & \left\{\frac{1}{6},\frac{5}{6}\right\}& \xmark\\
 \cline{2-4}
 & 6 & \left\{\frac{1}{6},\frac{1}{6},\frac{1}{6},\frac{1}{6},\frac{1}{6},\frac{1}{6}\right\} & \cmark \\
\hline
 6,8 & 3,1 & \left\{\frac{1}{6},\frac{1}{6},\frac{1}{6}\right\}, \left\{\frac{1}{8},\frac{3}{8}\right\}& \xmark\\
\hline 
6,12 & 2,1 &\left\{\frac{1}{6},\frac{1}{6}\right\}, \left\{\frac{1}{12},\frac{7}{12}\right\}& \xmark\\
\cline{2-4}
& 3,1 &\left\{\frac{1}{6},\frac{1}{6},\frac{1}{6}\right\}, \left\{\frac{1}{12},\frac{5}{12}\right\}& \xmark\\
\hline
7 & 1 & \left\{\frac{1}{7},\frac{2}{7},\frac{4}{7}\right\} & \cmark\\ 
\hline
8 & 2 & \left\{\frac{1}{8},\frac{1}{8},\frac{3}{8},\frac{3}{8}\right\} & \cmark\\
\hline
8,12 & 1,1 & \left\{\frac{1}{8},\frac{3}{8}\right\},\left\{\frac{1}{12},\frac{5}{12}\right\}& \xmark\\
\hline
12 & 2 & \left\{\frac{1}{12},\frac{1}{12},\frac{5}{12},\frac{5}{12}\right\}& \cmark\\ 
\hline   
15 & 1 & \left\{\frac{1}{15},\frac{2}{15},\frac{4}{15},\frac{8}{15}\right\} & \cmark\\ 
\hline 
16 & 1 & \left\{\frac{1}{16},\frac{3}{16},\frac{5}{16},\frac{7}{16}\right\} & \cmark\\ 
\hline 
20 & 1 & \left\{\frac{1}{20},\frac{3}{20},\frac{7}{20},\frac{9}{20}\right\} & \cmark\\ 
\hline 
24 & 1 & \left\{\frac{1}{24},\frac{5}{24},\frac{7}{24},\frac{11}{24}\right\} & \cmark\\ 
\hline 
\end{array}$
\renewcommand{\arraystretch}{1}
\end{adjustbox}
\caption{Classification of the data described in Lemma \ref{lem-partsad}}
\label{tab-partitions}
\end{table}

\newpage

We also recall a simple fact from the theory of abelian varieties:

\lem[lem-blrat]{Let $A$ be an abelian variety of dimension $n$, and $g\in\mathrm{Aut}(A)$ of finite order. Denote by $P(g)$ the characteristic polynomial of $M(g)$. Then $P(g)\overline{P(g)}$ is a product of cyclotomic polynomials.}

\demo{By \cite[Proposition 1.2.3]{BirkLang}, the matrix $M(g)\oplus\overline{M(g)}$ in $\GL_{2n}(\C)$ is similar to an element of $\GL_{2n}(\Q)$. Hence, $P(g)\overline{P(g)}$ is a polynomial over $\Q$. Since $g$ has finite order, the roots of this polynomial are roots of unity. Remembering that cyclotomic polynomials are the minimal polynomials of roots of unity over $\Q$, an easy induction shows that there is a product $\Pi$ of cyclotomic polynomial that has the exact same roots as $P(g)\overline{P(g)}$. But since both cyclotomic polynomials and characteristic polynomials are unitary, it means that $P(g)\overline{P(g)}=\Pi$.}

We can now prove Proposition \ref{prop-junclass}.

\demode{Proposition \ref{prop-junclass}.}{Denote by $d$ the order of the junior element $g$, denote by $(e^{2i\pi a_j/d})_{1\le j\le n}$ its ranked vector of eigenvalues, and by $P(g)$ the characteristic polynomial of its matrix $M(g)$.
As $g$ itself acts freely in codimension $2$ and $g$ is junior, it must be that $d\ge 3$.

By Lemma \ref{lem-blrat}, there are positive integers $k$, $(u_m)_{1\le m\le k}$ ordered increasingly, and $(\alpha_m)_{1\le m\le k},$ such that:

\begin{equation}\label{eq-polyn}
\prod_{j=1}^n (X-e^{2i\pi a_j/d})(X-\overline{e^{2i\pi a_j/d}}) = P(g)\overline{P(g)} = \prod_{m=1}^k {\Phi_{u_m}}^{\alpha_m}.
\end{equation}

Note that $\Phi_{u_m}(e^{2i\pi a_j/d})=0$, or equivalently $\Phi_{u_m}(\overline{e^{2i\pi a_j/d}})=0$, if and only if $u_m=\frac{d}{d\wedge a_j}$. We define the following partition of $\lint 1,n\rint$
$$\begin{array}{lll}
\mbox{for $m\in\lint 1,k\rint$, }\; & I_m &:= \{j\in\lint 1,n\rint\mid u_m=\frac{d}{d\wedge a_j}\};\\
& {\bf A_m} &:= \{\{a_j\mid j\in I_m\}\},\mbox{ as a multiset.}
\end{array}$$

\noindent By Identity \ref{eq-polyn}, for $m\in\lint 1,k\rint$ we have:
\begin{equation}\label{eq-ajmult}
{\bf A_m}\cup (d-{\bf A_m})
= \{\{r\in \lint 1,d-1\rint \mid u_m=\frac{d}{d\wedge r}=0\}\}^{*\alpha_m}
\end{equation}

\noindent Moreover, since $g$ is junior:
\begin{equation}\label{eq-juniorsaj}
1=\sum_{j=1}^n \frac{a_j}{d}
= \sum_{m=1}^k \sum_{j\in I_m} \frac{a_j}{d}
= \sum_{m=1}^k \sum_{j\in I_m} \frac{a_j}{u_m(d\wedge a_j)}
= \sum_{m=1}^k S_{{\bf A_m},d}(u_m).
\end{equation}

So, possibly leaving out the data of index 1, if $u_1=1$ (which is determined by the multiplicity $\alpha_1\in\N$, since then ${\bf A_1}=\{\{{\bf 0}_{\alpha_1}\}\}$ and $S_{{\bf A_1},d}(u_1)=0$), Lemma \ref{lem-partsad} applies, showing that there are finitely many possibilities for 
$$k,\, (u_m)_{1\le m\le k},\, (\alpha_m)_{1\le m\le k},\,\left(\frac{1}{d}{\bf A_m}\right)_{1\le m\le k}$$
and listing them. 
We exclude by hand a lot of these possibilities using the assumption that $\grp{g}$ acts freely in codimension $2$ on $A$, i.e., that for all $\ell\in \lint 1,d-1\rint$, there must be distinct indices $j_1(\ell),j_2(\ell),j_3(\ell)\in\lint 1,n\rint$, such that none of the $\frac{\ell a_{j_i(\ell)}}{d}$ is an integer. What remains then is precisely the list in Table \ref{tab-propjunclass}.
}

\section{Cyclicity of the pointwise stabilizers of loci of codimension 3 and 4}\label{sec-CyclicPointwiseStabilizer}

We now know that $G$ is generated by junior elements, which we have classified into twelve different types. However, this is by far insufficient to determine the structure of $G$. Even locally, for $W\subset A$ a subvariety, the pointwise stabilizer
$$\pstab(W):=\{g\in G\mid \forall w\in W,\, g(w)=w\}$$
could as well be cyclic and generated by one junior element, as it could be more complicated, e.g., if it contained non-commuting junior elements.

In this section, we show that in fact, if $W$ has codimension 3 or 4 in $A$, $\pstab(W)$ is trivial or cyclic, generated by one junior element. 
Let us outline the proof.
Subsection \ref{sec-ReduceToEquivariantComplement} reduces to proving this in the case when $W$ is a point in an abelian variety $B$ of dimension 3 or 4. Up to conjugating the whole group $G$ by a translation, we therefore just work on the case $W=\{0\}$.
Assuming $\pstab(W)$ is not trivial, we can then find a junior element $g\in\pstab(W)$, that is of one of the twelve types of Section \ref{sec-JuniorClassification}. Subsection \ref{subsec-AbelianClassification} exhibits a correlation between the type of $g$ and the isogeny type (possibly even isomorphism type) of the abelian variety $B$ on which it acts. A corollary is that if $g,\,h\in\pstab(W)$ are two junior elements, then they should either have the same type, or one is of type $\left(\1_{n-4},\omega,\omega,\omega,-1\right)$ and the other $\left(\1_{n-3},j,j,j\right)$, or one is of type $\left(\1_{n-4},i,i,i,i\right)$ and the other $\left(\1_{n-4},\zeta_{12},\zeta_{12},\zeta_{12}^5,\zeta_{12}^5\right)$. In particular, if $\pstab(W)$ is cyclic, it must indeed be generated by one junior element.
The conclusive Subsection \ref{subsec-GroupThPointwiseStabilizer} is the most technical. For any given abelian three- or fourfold $B$ of one of the types just defined, we classify all finite subgroups of 
$$\mathrm{Aut}(B,0):=\{f\in\mathrm{Aut}(B)\mid f(0)=0,\mbox{ i.e., }T(f)=0\}$$
that act freely in codimension $2$ on $B$ and are generated by junior elements.
The main idea is to bound the order of such groups, to scrutinize the finite list arising, and to rule out all but the cyclic case of the list by the assumption on generators.

\subsection{Reduction to a 3 or 4-dimensional question}\label{sec-ReduceToEquivariantComplement}

\defi[defi-complementary]{Let $A$ be an abelian variety. An {\it abelian subvariety} of $A$ is a closed subvariety of $A$ that is also a subgroup of the abelian group $(A,+)$.
A {\it translated abelian subvariety} of $A$ is the image by a translation of an abelian subvariety of $A$.

We say that two translated abelian subvarieties $B$ and $C$ of $A$ are {\it complementary} if one of the following equivalent statements hold: 
\begin{enumerate}[label = \textup{(\roman*)}]
\item $B\cap C$ is non-empty and, for some $p\in B\cap C$, it holds:
$$H^0(\{p\},T_B)\oplus H^0(\{p\},T_C) = H^0(\{p\},T_A).$$
\item The addition map $i:B\times C \to A$ is an isogeny.
\end{enumerate}
}

\demo{$(i)\Rightarrow (ii)$: as the translation by $(p,p)$, respectively by $2p$, is an isomorphism from $B\times C$ to $(B-p)\times(C-p)$, respectively of $A$, it is enough to prove the statement for $p=0$. As $\dim(A)=\dim(B\times C)$ and the varieties are regular, we simply check that $i$ is quasi-finite. Since $B\cap C$ is the intersection of two abelian subvarieties of $A$ satisfying:
$$H^0(\{0\},T_{B})\cap H^0(\{0\},T_{C}) = \{0\},$$
the set $B\cap -C$ is discrete in $A$, hence finite. For $a\in\mathrm{Im}(i)$, say $a=i(a_B,a_C)$, we can express the fiber $i^{-1}(a)=\{(b+a_B,-b+a_C)\mid b\in B\cap -C\}$, so it is finite, and $i$ is indeed quasi-finite.

\medskip

\noindent $(ii)\Rightarrow (i)$: fix $c_0\in C$. The addition $i$ is onto, so let $(p,c)\in B\times C$ be such that $p+c=2c_0$. Clearly, $p=2c_0-c\in B\cap C$, and as $i$ is locally analytically an isomorphism,
 $$H^0(\{p\},T_B)\oplus H^0(\{p\},T_C) = H^0(\{2p\},T_A)= H^0(\{p\},T_A).$$
}

\rem[rem-comptrans]{If $B$ and $C$ are complementary translated abelian subvarieties of an abelian variety $A$, and $t\in A$ is any point, then $B+t$ and $C$ are complementary as well.
Our notion of complementarity is weaker than the notion defined for abelian subvarieties in \cite[p.125]{BirkLang}.}

Let us now state our reduction result. Note that it applies not only in codimension 3 and 4, but in any higher codimension as well. For a finite group $G$ acting on a projective variety $X$, we denote by 
${\rm Fix}(G)=\{x\in X\mid \pstab(x)=G\}.$
 
\prop[prop-point-wise-stabilizer]{Let $A$ be an abelian variety, $G$ be a finite group acting freely in codimension $2$ on $A$. Suppose that the quotient $A/G$ admits a $K$-trivial resolution. Let $H$ be a non-trivial subgroup of $G$ such that ${\rm Fix}(H)$ is non-empty, and let $W$ be an irreducible component of ${\rm Fix}(H)$. Then:
\begin{enumerate}[label = \textup{(\arabic*)}]
\item The variety $W$ is a translated abelian subvariety of $A$, and for any $t\in W$, there is a translated abelian subvariety $B$ of $A$ which is $\pstab(W)$-stable, contains $t$, and is complementary to $W$ in $A$.
\item If $t$ and $B$ are as in (1), then an element $g\in \pstab(W)$ is junior if and only if $g|_B\in\mathrm{Aut}(B, t)$ is a junior element.
\item The group $\pstab(W)\subset\mathrm{Aut}(B,t)$ is generated by junior elements.
\end{enumerate}}

\demode{Proposition \ref{prop-point-wise-stabilizer}.}{Up to conjugating the $G$-action by the translation by $t$, we can assume that $t=0\in W$. Let us establish (1). First, the fact that $W$ is an abelian subvariety of $A$ follows from a strong induction on the order of $H$, using \cite[13.1.2.a)]{BirkLang} for the induction step. Second, as $G$ is finite, we can take a $G$-invariant polarization $L$ on $A$. We can apply \cite[Prop.13.5.1]{BirkLang}: there is a unique complementary abelian subvariety $(B,L|_B)$ to $(W,L|_{W})$ in $(A,L)$, and it is $\pstab(W)$-stable. By Remark \ref{rem-comptrans}, $B$ and $W$ are complementary in our sense as well.

We now prove (2): let $g\in \pstab(W)$. As $g$ fixes all points of $B\cap W$, its restriction $g|_B$ has a fixed point. As $g(B)=B$, we have:
$$M(g)=\left(\begin{array}{cc}
\id_{\dim(W)} & 0\\
0 & M(g|_B)\\
\end{array}\right),$$ 
and thus $g$ is indeed junior if and only if $g|_B$ is.

We move on to (3). Take a general point $w\in W$ such that $ \Stab(w)= \pstab(W)$ and apply Corollary \ref{cor-pstabjun} at the point $w$: It shows that $\pstab(W)$ is generated by junior elements.
}

\subsection{The abelian varieties corresponding to the twelve juniors}\label{subsec-AbelianClassification}

Let $A$ be an abelian variety of dimension $n$, $G$ be a finite group acting freely in codimension $2$ on $A$ such that $A/G$ has a resolution $X$ with $c_1(X)=0$. By Proposition \ref{prop-Reid}, $G\subset\mathrm{Aut}(A)$ must entail a junior element presented in Table \ref{tab-propjunclass} (up to its translation part, and up to similarity for its linear part). The fact that, in some coordinates, a given matrix of Table \ref{tab-propjunclass} acts as an automorphism on the abelian variety $A$ imposes some restrictions. Using the theory of abelian varieties with complex multiplication, these restrictions are investigated by Proposition \ref{prop-singab}.

\nota[def-abelianplayers]{Let us define the following quadratic integers.
$$u_7=\frac{-1+i\sqrt{7}}{2},\,u_8=i\sqrt{2},\,u_{15}=\frac{1+i\sqrt{15}}{2},\,u_{20}=i\sqrt{5},\,u_{24}=i\sqrt{6},\,$$
We also define the following conjugated algebraic integers, which are quadratic over the totally real extension $\Q(2\sqrt{2})$ of $\Q$.
$$u_{16}=i\sqrt{4+2\sqrt{2}},\, v_{16}=i\sqrt{4-2\sqrt{2}}.$$

\noindent For $z\in\C\setminus\R$, we define the elliptic curve $E_z:=\C/{\Z\oplus z\Z}$. If $z$ is a quadratic integer, then we denote by $\Z[z]$ the $\Z$-algebra that it generates. It holds $\Z[z]=\Z\oplus z\Z\subset\C$.

\noindent We also define the simple abelian surface $S_{u_{16},v_{16}}:=\C^2/{\Z[(u_{16},v_{16})]}$.}

\rem{Note that the simplicity of $S_{u_{16},v_{16}}$ follows from \cite[Prop.27]{Shimura}.}

With these notations, we can state the main result of the subsection.

\prop[prop-singab]{Let $A$ be an abelian variety. Suppose that there is a junior element $g\in\mathrm{Aut}(A)$, and that $\grp{g}$ acts freely in codimension $2$ on $A$. Denote by $W$ an irreducible component of $\mathrm{Fix}(g):=\{a\in A\mid g(a)=a\}$. 
Let $B$ be a complementary to $W$ in $A$. Then the isogeny type of $B$ is entirely determined by the type of the junior element $g$ by Table \ref{tab-abvarjun}, unless $g$ is of type $(\1_{n-4},\omega,\omega,\omega,-1)$. 
Moreover, the isomorphism type of a $\grp{g}$-stable complementary $B_{\mathrm{st}}$ to $W$ in $A$ is also entirely determined by the type of $g$, unless $g$ is of type $(\1_{n-4},\omega,\omega,\omega,-1)$ or $(\1_{n-5},\omega,\omega,\omega,\omega,j)$.}

\newpage

\begin{table}[H]
\begin{center}
\renewcommand{\arraystretch}{1.4}
$\begin{array}{|c|c|c|}
\hline 
 \mbox{type of }g & \mbox{isogeny type of }B & \mbox{isomorphism type of }B_{\mathrm{st}} \\
 \hline
\left(\1_{n-3},j,j,j\right) & {E_j}^3 & {E_j}^3\\
\hline
\left(\1_{n-4},i,i,i,i\right) & {E_i}^4 & {E_i}^4\\
\hline
\left(\1_{n-4},\omega,\omega,\omega,-1\right) & E\times {E_j}^3\mbox{ for some }E & \mbox{not determined} \\
\hline
\left(\1_{n-5},\omega,\omega,\omega,\omega,j\right) & {E_j}^5 & \mbox{not determined}\\
\hline
\left(\1_{n-6},\omega,\omega,\omega,\omega,\omega,\omega\right) & {E_j}^6 & {E_j}^6\\
\hline
\left(\1_{n-3},\zeta_7,{\zeta_7}^2,{\zeta_7}^4\right) & {E_{u_7}}^3 & {E_{u_7}}^3\\
\hline
\left(\1_{n-4},\zeta_8,\zeta_8,\zeta_8^3,\zeta_8^3\right) & {E_{u_8}}^4 & {E_{u_8}}^4\\
\hline
\left(\1_{n-4},\zeta_{12},\zeta_{12},\zeta_{12}^5,\zeta_{12}^5\right) & {E_i}^4 & {E_i}^4\\
\hline
\left(\1_{n-4},\zeta_{15},\zeta_{15}^2,\zeta_{15}^4,\zeta_{15}^8\right) & {E_{u_{15}}}^4 & {E_{u_{15}}}^4\\
\hline
\left(\1_{n-4},\zeta_{16},\zeta_{16}^3,\zeta_{16}^5,\zeta_{16}^7\right) & {S_{u_{16},v_{16}}}^2 & {S_{u_{16},v_{16}}}^2\\
\hline
\left(\1_{n-4},\zeta_{20},\zeta_{20}^3,\zeta_{20}^7,\zeta_{20}^9\right) & {E_{u_{20}}}^4 & {E_{u_{20}}}^4\\
\hline
\left(\1_{n-4},\zeta_{24},\zeta_{24}^5,\zeta_{24}^7,\zeta_{24}^{11}\right) & {E_{u_{24}}}^4 & {E_{u_{24}}}^4\\
\hline
\end{array}$
\renewcommand{\arraystretch}{1}
\end{center}
\caption{Correspondence between types of junior elements and types of abelian varieties.}
\label{tab-abvarjun}
\end{table}

\nota{Let $V$ be a $\C$-vector space, $f:V\to V$ be a linear map. We denote by $\mathrm{EVal}(f)$ the set of eigenvalues of $f$ in $\C$, by ${\bf EVal}(f)$ the multiset of eigenvalues of $f$ in $\C$ counted with multiplicities. If $\lambda\in\mathrm{EVal}(f)$, we denote by $E_f(\lambda)$ the eigenspace of $f$ for the eigenvalue $\lambda$.

\noindent We denote by $Z(\Phi_d)$ the set of primitive $d$-th roots of unity in $\C$.}

Let us first carry out an important computation, that makes plain where these special types of abelian varieties come from.
Let $k\ge 3$ be an integer.
There is a natural action of $\zeta_k\otimes 1$ on the algebra $\Z[\zeta_k]\otimes\C$. We compute its eigenvalues. By definition, $\Z[\zeta_{k}]\otimes\C$ is the quotient algebra $\C[X]/(\Phi_k)$, multiplication by $\zeta_{k}\otimes 1$ corresponding to multiplication by the class $X+\Phi_k\C[X]$. So $\xi\in\C$ is an eigenvalue with eigenvector $P+\Phi_k\C[X]$ if and only if $P\not\in\C[X]$ and $XP -\xi P\in\Phi_k\C[X]$, or equivalently, $\xi$ is a root of $\Phi_k$ and $P\in \frac{\Phi_k}{X-\xi}\C[X]$. Hence the linear decomposition

\begin{equation}\label{eq-deceigen}
\Z[\zeta_{k}]\otimes\C=
\bigoplus_{\xi\in Z(\Phi_k)} E_{\zeta_{k}\otimes 1}(\xi)
\end{equation}

Now, consider a subset $S_k$ of $Z(\Phi_k)$ such that $S_k\cup\overline{S_k}=Z(\Phi_k)$. For example, if we let $g$ be a junior element of one of the twelve types in Table \ref{tab-propjunclass}, and we assume that $g$ has an eigenvalue of order $k$, we could set $S_k=S_k(g)=\mathrm{EVal}(g)\cap Z(\Phi_k)$. 
This defines a $\Z$-linear inclusion

\begin{equation}\label{eq-fgk}
f(S_k):\Z[\zeta_{k}]\hookrightarrow
\bigoplus_{\xi\in S_k} E_{\zeta_{k}\otimes 1}(\xi)
\;\simeq\;\C^{\phi(k)/2}
\end{equation}

It is worth noting that the $\Z$-linear inclusion $f(S_k)\oplus\overline{f(S_k)}$ corresponds to the natural inclusion of $\Z[\zeta_k]$ in $\Z[\zeta_k]\otimes\C$ given by Identity \ref{eq-deceigen}.

The following lemma is key.

\lem[lem-specialfactor]{If $S_k=S_k(g)$ for a junior element $g$ of Table \ref{tab-propjunclass}, then the corresponding abelian variety $\C^{\phi(k)/2}/\mathrm{Im}(f(S_k))$ is described in Table \ref{tab-specialfactor}.}

\rem{For $k=3,4,6,8,12$, we consider the sets 
$S_k:=\{j\}$, $\{i\}$, $\{j\}$, $\{\zeta_8,{\zeta_8}^3\}$, and $\{\zeta_{12},{\zeta_{12}}^5\}$ respectively, and Lemma \ref{lem-specialfactor} is \cite[Cor.13.3.4, Cor.13.3.6]{BirkLang}. In the other cases, the computation relies on the same ideas as \cite[Cor.13.3.6]{BirkLang}, as we will soon see.}

\begin{table}[H]
\begin{center}
\renewcommand{\arraystretch}{1.4}
$\begin{array}{|c|c|c|c|}
\hline 
k & S_k & \C^{\phi(k)/2}/\mathrm{Im}(f(S_k)) \\
\hline
3 & \{j\} & E_j \\
\hline
4 & \{i\} & E_i \\
\hline
6 & \{\omega\} & E_j\\
\hline
7 & \{\zeta_7,{\zeta_7}^2,{\zeta_7}^4\} & {E_{u_7}}^3\\
\hline
8 & \{\zeta_8,\zeta_8^3\} & {E_{u_8}}^2 \\
\hline
12 & \{\zeta_{12},\zeta_{12}^5\} & {E_i}^2 \\
\hline
 15 & \{\zeta_{15},\zeta_{15}^2,\zeta_{15}^4,\zeta_{15}^8\} & {E_{u_{15}}}^4 \\
\hline
16 & \{\zeta_{16},\zeta_{16}^3,\zeta_{16}^5,\zeta_{16}^7\} & {S_{u_{16},v_{16}}}^2 \\
\hline
20 & \{\zeta_{20},\zeta_{20}^3,\zeta_{20}^7,\zeta_{20}^9\} & {E_{u_{20}}}^4 \\
\hline
24 & \{\zeta_{24},\zeta_{24}^5,\zeta_{24}^7,\zeta_{24}^{11}\} & {E_{u_{24}}}^4 \\
\hline
\end{array}$
\renewcommand{\arraystretch}{1}
\end{center}
\caption{Computing $\C^{\phi(k)/2}/\mathrm{Im}(f(S_k))$ for given $S_k$ stemming from a junior element.}
\label{tab-specialfactor}
\end{table}

To complete the proof Lemma \ref{lem-specialfactor} for $k=7,15,16,20,24$, we use a part of \cite[Proof of Thm.3, Page 46]{Shimura}, recalled here without proof.

\lem[lem-shimu]{Let $K=\Q(\alpha)$ be a totally imaginary extension of $\Q$ of degree $2m$, that is a quadratic extension of a totally real extension $K_0$ of $\Q$. Let $F$ be a finite Galois extension of $K$, of degree $2r$ over $\Q$. Let $\{\varphi_i\}_{1\le i\le r}$ be morphisms of $\Q$-algebras defined from $F$ to $\C$ such that:
 $$\mathrm{Hom}_{\Q-\mathrm{alg}}(F,\C)=\mathrm{Vect}_{\Q}\left(\varphi_1,\overline{\varphi_1},\ldots,\varphi_r, \overline{\varphi_r}\right).$$
Suppose also that no two of the restrictions $\varphi_i|_K$ are conjugated. 

Then we can restrict $m$ of these morphisms, defining $\psi_{j}=\varphi_{i_j}|_K$ for some distinct $i_j$ with $j\in\lint 1,m\rint$, such that: 
$$\mathrm{Hom}_{\Q-\mathrm{alg}}(K,\C)=\mathrm{Vect}_{\Q}\left(\psi_1,\overline{\psi_1},\ldots,\psi_m,\overline{\psi_m}\right).$$
We obtain a $\Z$-algebra $\Delta:=\Z[(\psi_{1}(\alpha),\ldots,\psi_{m}(\alpha))]$ that is a lattice of rank $2m$ in $\C^m$. The complex torus $A:=(\C^m/\Delta)^{n/m}$ is an abelian variety of CM-type $(F,\{\varphi_i\}_{1\le i\le f})$.}

\demode{Lemma \ref{lem-specialfactor}.}{Let $F=\Q[\zeta_k]$, $r=\frac{\phi(k)}{2}$. Let us define 
$\{\varphi_i\}_{1\le i\le r}$: Composing $f(S_k)$ defined in Identity \ref{eq-fgk} with the projections on the $r$ eigenspaces, we obtain morphisms of $\Z$-algebras $f_i:\Z[\zeta_k]\to \C,$
which we tensor by $\Q$ and normalize to define morphisms of $\Q$-algebras:
$$\varphi_i=\frac{1}{f_i(1)}(f_i\otimes 1) :\Q[\zeta_k]\to\C.$$

By Identities \ref{eq-deceigen} and \ref{eq-fgk}, the morphisms $\{\varphi_i,\overline{\varphi_i}\}_{1\le i\le r}$ are linearly independent over $\Q$, whereas the morphisms $\{\varphi_i\}_{1\le i\le r}$ define an embedding of $F$ into the $\Q$-algebra of linear endomorphisms of the abelian variety $\C^{\phi(k)/2}/\mathrm{Im}(f(S_k))$. In other words, the abelian variety $\C^{\phi(k)/2}/\mathrm{Im}(f(S_k))$ has CM-type $\left(F,\{\varphi_i\}_{1\le i\le r}\right)$.
This is in fact the sole abelian variety with this CM-type, by \cite{MM}, \cite[Prop.17]{Shimura}, and remembering that $k\in\{7,15,16,20,24\}$. 

We get the wished description of the abelian variety $\C^{\phi(k)/2}/\mathrm{Im}(f(S_k))$ by applying Lemma \ref{lem-shimu} for $K=\Q(u_k)$ and $K_0=\Q$ for $k\in\{7,15,20,24\}$, and for $K=\Q(u_{16})$ and $K_0=\Q(2\sqrt{2})$ for $k=16$, and noting that:

\begin{itemize}
\item $u_7 = \zeta_{7}+\zeta_{7}^2+\zeta_{7}^4,$ 
\item $u_{15} = \zeta_{15}+\zeta_{15}^2+\zeta_{15}^4+\zeta_{15}^8,$
\item $u_{16}=\zeta_{16}+\zeta_{16}^3+\zeta_{16}^5+\zeta_{16}^7$ and $v_{16}=\zeta_{16}^3+\zeta_{16}^5+\zeta_{16}^{9}+\zeta_{16}^{15}$,
\item $u_{20} = \zeta_{20}+\zeta_{20}^3+\zeta_{20}^7+\zeta_{20}^9,$
\item $u_{24} = \zeta_{24}+\zeta_{24}^5+\zeta_{24}^7+\zeta_{24}^{11}.$
\end{itemize}}

The next result follows almost effortlessly from the ideas of \cite[Page 333-334]{III0}.

\lem[lem-coollemma]{Let $B$ be an abelian variety. Suppose that there is an automorphism $g$ of $B$ whose set of eigenvalues is one of the $S_k$ in Table \ref{tab-specialfactor}. Then $B$ is isomorphic to a power of the abelian variety $\C^{\phi(k)/2}/\mathrm{Im}(f(S_k))$.}

\demo{Let $\Lambda$ be a lattice in $\C^n$ such that $B=\C^n/\Lambda$. The action of $g$ on $B$ has a linear part $M(g)\in\GL_n(\C)$ that satifies $M(g)(\Lambda)=\Lambda$. Note that $1$ does not appear in the set of eigenvalues $S_k$ of $M(g)$, hence $g$ and $M(g)$ have the same multiplicative order $k$. This provides $\Lambda$ with the structure of a $\Z[g]$-module, i.e., a $\Z[\zeta_k]$-module structure, since the minimal polynomial of $g$ is $\Phi_k$. As such, $\Lambda$ is finitely-generated and torsion-free. But by \cite{MM}, since $k\in\lint 3,20\rint\cup\{24\}$, the ring of cyclotomic integers $\Z[\zeta_k]$ is a principal ideal domain. So, by the structure theorem for finitely-generated modules over principal ideal domains, $\Lambda\simeq \Z[\zeta_k]^{2n/\phi(k)}$, and the action of $g$ on $\Lambda$ identifies with the multiplication by $\zeta_k$ on $\Z[\zeta_k]^{2n/\phi(k)}$.

The embedding $\Lambda\hookrightarrow H^0(B,T_B)\simeq\C^n$ can be recovered from the action of $g$ on $\Lambda$. Indeed, there is an induced action of $g\oplus\overline{g}$ on $\Lambda\otimes \C = H^0(B,T_{B,\R}\otimes\C)\simeq\C^{2n}$. This action splits into two blocks: $g$ is acting on $H^0(B,T_B)$ and $\overline{g}$ is acting on its supplementary conjugate in $H^0(B,T_{B,\R}\otimes\C)$. By the requirement on its set of eigenvalues $S_k$, $g$ has no eigenvalue in common with $\overline{g}$, and therefore:
$$H^0(B,T_B)=\bigoplus_{\xi\in\mathrm{{\bf EVal}}(g)}E_{g\oplus\overline{g}}(\xi).$$

Hence, the corresponding embedding $\Z[\zeta_k]^{2n/\phi(k)}\hookrightarrow\C^n$ must similarly be given by:
$$\C^n=\bigoplus_{\xi\in\mathrm{{\bf EVal}}(g)}E_{\zeta_k\otimes 1}(\xi),$$
where $\zeta_k\otimes 1$ is the action by componentwise multiplication on $\Z[\zeta_k]^{2n/\phi(k)}\otimes\C$. In other words, this embedding is the blockwise embedding $f(S_k)$, repeated on $\frac{2n}{\phi(k)}$ blocks of dimension $\frac{\phi(k)}{2}$ each. So
$B\simeq \left(\C^{\phi(k)/2}/\mathrm{Im}(f(S_k))\right)^{2n/\phi(k)}.$}

The proof of Proposition \ref{prop-singab} is now easy.

\demode{Proposition  \ref{prop-singab}}{By Proposition \ref{prop-point-wise-stabilizer}, let $B_{\mathrm{st}}$ be a $\grp{g}$-stable complement to $W$ in $A$. For any other complement $B$ to $W$, since $B\times W$ and $B_{\mathrm{st}}\times W$ are isogenous, $B$ and $B_{\mathrm{st}}$ are isogenous. Let us determine the isogeny (and if possible isomorphism) type of $B_{\mathrm{st}}$.

On one hand, if $g$ is of type $(\1_{n-4},\omega,\omega,\omega,-1)$ or $(\1_{n-5},\omega,\omega,\omega,\omega,j)$, then $g|_{B_{\mathrm{st}}}$ has eigenvalues of two different orders. By \cite[Thm.13.2.8]{BirkLang}, there are then two $\grp{g}$-stable complementary translated abelian subvarieties $B_1$ and $B_2$ in $B_{\mathrm{st}}$, such that all eigenvalues of $g|_{B_1}$ have order $k_1=6$, and all eigenvalues of $g|_{B_2}$ have the same order $k_2<6$. By definition, $B_{\mathrm{st}}$ is isogenous to $B_1\times B_2$, and thus its isogeny type can be derived from the isomorphism types of $B_1$ and $B_2$, given by Lemma \ref{lem-coollemma} if $k_1,k_2\ge 3$. However, if $g$ is of type $(\1_{n-4},\omega,\omega,\omega,-1)$, then $k_2=2$ and $B_2$ can be any elliptic curve, and that is why the isogeny type of $B_{\mathrm{st}}$ is not entirely determined in this case.

On the other hand, if $g$ is of any other type, then all eigenvalues of $g|_{B_{\mathrm{st}}}$ are of the same order $k\ge 3$, and Lemma \ref{lem-coollemma} determines the isomorphism type of $B_{\mathrm{st}}$.}

\subsection{Group theoretical treatment of a point's stabilizer in dimension 3 or 4}\label{subsec-GroupThPointwiseStabilizer}

We can now establish the following proposition.

\prop[prop-cyclicincod4]{Let $A$ be an abelian variety, $G\subset\mathrm{Aut}(A)$ be a finite group acting freely in codimension $2$. Suppose that the quotient $A/G$ admits a $K$-trivial resolution. Let $W$ be a subvariety of codimension $m\le 4$ in $A$ such that $\pstab(W)\ne\{1\}$. Then $\pstab(W)$ is a cyclic group generated by one junior element.}

By Proposition \ref{prop-point-wise-stabilizer}, it reduces to proving the following result.

\prop[prop-justshowFcyc]{Let $B$ be an abelian variety of dimension $m\le 4$, $F\subset\mathrm{Aut}(B,0)$ be a finite group acting freely in codimension $2$ and fixing $0\in B$. Suppose that $F$ is generated by junior elements. Then $F$ is a cyclic group generated by one junior element.}

We refer the reader to \cite{Robinson},\cite{Isaacs} for standard facts in finite group theory, and in particular Sylow theory and representation theory. Let us just recall a few notations used in the following.

\nota{We denote by $C_F(H)$, respectively $N_F(H)$, the centralizer, respectively normalizer, of a subset $H$ of a group $F$, i.e.,
\begin{align*}
C_F(H)&:=\{f\in F\mid\forall\, h\in H,\, fh = hf\}\\
N_F(H)&:=\{f\in F\mid fH = Hf\}
\end{align*}
By associativity of the group law, we see that $C_F(H)$ and $N_F(H)$ are subgroups of $F$.}

\nota{Let $F$ be a finite group. A {\it class function} $\xi:F\to\C$ is a function that takes constant values within each conjugacy class in $F$. We denote by $\langle\cdot,\cdot\rangle$ the inner product on the space of class functions defined by
$$\langle\xi,\zeta\rangle := \frac{1}{|F|}\sum_{f\in F}\xi(f)\overline{\zeta(f)}.$$

\noindent Let $V$ be a $\C$-vector space of finite dimension. A {\it representation} of $F$ in $V$ is a group morphism $\rho:F\rightarrow\GL(V)$.
The {\it character} $\chi$ of a representation $\rho$ is the class function $\chi:f\in F\to \mathrm{Tr}(\rho(f))\in\C$.

Let us fix a representation $\rho:F\to\GL(V)$. By Schur's lemma, it decomposes as a direct sum of irreducible representations: $$\rho =\rho_1^{\oplus n_1}\oplus\ldots\oplus\rho_k^{\oplus n_k},$$ 
and, accordingly, if $\chi_i$ denotes the character of $\rho_i$, we have $\chi=n_1\chi_1+\ldots+n_k\chi_k$. 
The characters of the irreducible representations $\rho_1,\ldots,\rho_k$ form an orthonormal family in the space of $\C$-valued class functions of $G$, hence
$$\langle\chi, \chi\rangle=(n_1^2+\ldots+n_k^2)|F|.$$
We refer to $u=n_1^2+\ldots+n_k^2$ as the {\it splitting coefficient} of the representation $\rho$.}

We start proving lemmas towards Proposition \ref{prop-justshowFcyc}.
The first lemma classifies all possible finite order elements in $\mathrm{Aut}(B,0)$ of determinant one acting freely in codimension $2$, when $B$ is an abelian fourfold.

\newpage

\lem[lem-allelements]{Let $B$ be an abelian fourfold, and $g\in\mathrm{Aut}(B,0)$ be a finite order element of determinant one such that $\grp{g}$ acts freely in codimension $2$ on $B$. 
Then the order of $g$ and the matrix of a generator of $\grp{g}$ are given by a row of Table \ref{tab-allelements}. That row of Table \ref{tab-allelements} also gives some necessary conditions on $B$.}

\begin{table}[H]
\begin{center}
\renewcommand{\arraystretch}{1.3}
$\begin{array}{|c|c|c|}
\hline 
\mbox{order of }g & \mbox{a generator of $\grp{g}$ up to similarity}& \mbox{restrictions on }B\\
\hline
1 & \id &\\
\cline{1-2}
2 & \mathrm{-\id} &\\
\cline{1-2}
3 & \diag(j,j^2,j,j^2) & \\
\cline{1-2}
4 & \diag(i,-i,i,-i) & \\
\cline{1-2}
5 & \diag(\zeta_5,\zeta_5^2,\zeta_5^3,\zeta_5^4) & \mbox{not studied}\\
\cline{1-2}
6 & \diag(\omega,\omega^5,\omega,\omega^5) & \\
\cline{1-2} 
8 & \diag(\zeta_{8},\zeta_{8}^3,\zeta_{8}^5,\zeta_{8}^7) &\\
\cline{1-2}
10 & \diag(\zeta_{10},\zeta_{10}^3,\zeta_{10}^7,\zeta_{10}^9) &\\
\cline{1-2} 
12 & \diag(\zeta_{12},\zeta_{12}^5,\zeta_{12}^7,\zeta_{12}^{11})&\\
\hline
\textcolor{red}{3} & \textcolor{red}{\diag(1,j,j,j)} & B\sim E\times {E_j}^3\\
\cline{1-2}
\textcolor{red}{6} & \textcolor{red}{\diag(-1,\omega,\omega,\omega)} & \\
\hline
9 & \diag(j^2,\zeta_9,\zeta_9^4,\zeta_9^7) & B\simeq {E_j}^4\\
\cline{1-2}
18 & \diag(\omega^5,\zeta_{18},\zeta_{18}^7,\zeta_{18}^{13}) & \\
\hline
\textcolor{red}{4} & \textcolor{red}{i\id} & \\
\cline{1-2}
\textcolor{red}{12} & \textcolor{red}{\diag(\zeta_{12},\zeta_{12}^5,\zeta_{12},\zeta_{12}^{5})} & B\simeq {E_i}^4\\
\cline{1-2}
20 & \diag(\zeta_{20},\zeta_{20}^9,\zeta_{20}^{13},\zeta_{20}^{17}) & \\
\hline
\textcolor{red}{7} & \textcolor{red}{\diag(1,\zeta_7,{\zeta_7}^2,{\zeta_7}^4)} & B\sim E\times {E_{u_7}}^3\\
\cline{1-2}
14 & \diag(-1,\zeta_{14},\zeta_{14}^9,\zeta_{14}^{11}) & \\
\hline
\textcolor{red}{8} & \textcolor{red}{\diag(\zeta_{8},\zeta_{8}^3,\zeta_{8},\zeta_{8}^3)} & B\simeq {E_{u_8}}^4\\
\cline{1-2}
24 & \diag(\zeta_{24},\zeta_{24}^{11},\zeta_{24}^{17},\zeta_{24}^{19}) & \\
\hline
\textcolor{red}{15} & \textcolor{red}{\diag(\zeta_{15},\zeta_{15}^2,\zeta_{15}^4,\zeta_{15}^8)} & B\simeq {E_{u_{15}}}^4\\
\cline{1-2}
30 & \diag(\zeta_{30},\zeta_{30}^{17},\zeta_{30}^{19},\zeta_{30}^{23}) & \\
\hline
\textcolor{red}{16} & \textcolor{red}{\diag(\zeta_{16},\zeta_{16}^3,\zeta_{16}^5,\zeta_{16}^7)} & B\simeq {S_{u_{16},v_{16}}}^2\\
\cline{2-2}
& \diag(\zeta_{16},\zeta_{16}^7,\zeta_{16}^{11},\zeta_{16}^{13}) & \\
\hline
\textcolor{red}{20} & \textcolor{red}{\diag(\zeta_{20},\zeta_{20}^3,\zeta_{20}^7,\zeta_{20}^9)} & B\simeq {E_{u_{20}}}^4\\
\hline
\textcolor{red}{24} & \textcolor{red}{\diag(\zeta_{24},\zeta_{24}^5,\zeta_{24}^7,\zeta_{24}^{11})} & B\simeq {E_{u_{24}}}^4\\
\hline
\end{array}$
\renewcommand{\arraystretch}{1}
\end{center}
\caption{Classification of possible elements of $g$ in $\mathrm{Aut}(B,0)$, with \textcolor{red}{colored junior elements}.}
\label{tab-allelements}
\end{table}

\newpage

\demo{Let $\zeta$ be an eigenvalue of $g$ of order $u$, such that $(\phi(u),u)$ is maximal in $\N^2$ for the lexicographic order. By Lemma \ref{lem-blrat}, $\Phi_u$ divides the characteristic polynomial $\chi_{g\oplus\overline{g}}$ in $\Q[X]$, so
$\phi(u)\le 2\,\mathrm{dim}\, B=8.$
Let us discuss cases:
\begin{enumerate}[label= (\arabic*)]
\item If $\phi(u) = 1$, then $u=1$ or $2$. As $g$ acts freely in codimension $2$ and has determinant one, $g=\pm\id_B$.
\item Suppose that $\phi(u)=8$. Then $g$ has four distinct eigenvalues of order $u$, and hence has order $u$. Listing integers of Euler number 8, $u\in\{15,16,20,24,30\}$. There is a generator $g'$ of $\grp{g}$ of which $e^{2i\pi/u}$ is an eigenvalue. Denote its other eigenvalues by $e^{2i\pi a/u}$, $e^{2i\pi b/u}$, $e^{2i\pi c/u}$, with 
\begin{itemize}
\item $a,b,c\in\lint 1,u-1\rint$ coprime to $u$
\item $u$ divides $1+a+b+c$
\item and \begin{align*}\Phi_u(X) = &(X-e^{2i\pi /u})(X-e^{2i\pi (u-1)/u})(X-e^{2i\pi a/u})(X-e^{2i\pi (u-a)/u})\\
&(X-e^{2i\pi b/u})(X-e^{2i\pi (u-b)/u})(X-e^{2i\pi c/u})(X-e^{2i\pi (u-c)/u})
\end{align*}
\end{itemize}
We check by hand the solutions to this system and plug them in Table \ref{tab-allelements}. For example, this is how we add $\diag(\zeta_{15},\zeta_{15}^2,\zeta_{15}^4,\zeta_{15}^8)$.
\item Suppose that $\phi(u)=6$. Then $g$ has three distinct eigenvalues of order $u$, and one eigenvalue of order $v$, with $\phi(v)=1$ or $2$. Since $g^u$ has three trivial eigenvalues and $\grp{g}$ acts freely in codimension $2$, $g^u=\id_B$, so $g$ has order $u$ and $v$ divides $u$. Listing the integers of Euler number 6, $u\in\{7,9,14,18\}$.
Using that $\chi_{g\oplus\overline{g}}=\Phi_u\Phi_v$ or $\Phi_u{\Phi_v}^2$, $g$ has determinant 1 and $\grp{g}$ acts freely in codimension $2$, we work out all possibilities by hand and add them to the table. One example falling in this case is $\diag(1,\zeta_7,{\zeta_7}^2,{\zeta_7}^4)$.
\item Suppose that $\phi(u)=4$. Then $g$ has two distinct eigenvalues of order $u$, and two remaining eigenvalues of respective order $v_1\le v_2$. As $\grp{g}$ acts freely in codimension $2$, $g^u$, which has two trivial eigenvalues, must be trivial, so $g$ has order $u$ and $v_1$ and $v_2$ divide $u$. Similarly, $g^{\mathrm{lcm}(v_1,v_2)}=\id_B$, so $u$ divides $\mathrm{lcm}(v_1,v_2)$. Listing integers of Euler number 4, $u\in\{5,8,10,12\}$. 
\begin{enumerate}[label= (\alph*)]
\item If $v_1$ divides $v_2$, then $v_2=u$. We investigate all possibilities of determinant 1 satisfying Lemma \ref{lem-blrat} by hand and add them to the table. One of them is $\diag(\zeta_5,\zeta_5^2,\zeta_5^3,\zeta_5^4)$.
\item If $v_1$ does not divide $v_2$, then by Lemma \ref{lem-blrat} again, $\phi(v_1)+\phi(v_2)\le 4$. Listing possibilities by hand, we see that $(v_1,v_2)\in\{(2,3),(3,4),(4,6)\}$. From the divisibility relations between $v_1$, $v_2$ and $u$, we obtain that $u=12$, and in fact, $(v_1,v_2)=(3,4)$ or $(4,6)$. In particular, $g$ has order $12$, so $g^6=-\id_B$, and so $g^3$ has four eigenvalues of order 4. But since $v_1=3$ or $v_2=6$, this can not be the case. Contradiction!
\end{enumerate}
\item The last case is when $\phi(u)=2$, i.e., $u=3,4,$ or $6$. In that case, each eigenvalue of $g$ has order $1,2,3,4,$ or $6$. As $\grp{g}$ acts freely in codimension $2$, $g$ has at most one eigenvalue of order $1$ or $2$. 
\begin{enumerate}[label= (\alph*)]
\item Suppose that $g$ has an eigenvalue of order 4. As it has determinant 1, it has an even number of eigenvalues of order 4, so at least two of them. Hence, by freeness in codimension $2$, $g^4=\id_B$, and so $g^2=-\id_B$, i.e., all eigenvalues of $g$ have order 4. There is a generator of $\grp{g}$ similar to either $\diag(i,i,i,i)$, or $\diag(i,-i,i,-i)$.
\item Suppose that $u=3$. Then as $(\phi(v),v)\le (\phi(u),u)$ for any order $v$ of another eigenvalue of $g$, the other eigenvalues have order $1,2,$ or $3$. Hence, there are at least three eigenvalues of order 3, and thus by freeness in codimension $2$, $g^3=\id_B$. So $g$ has order $3$ and there is a generator of $\grp{g}$ similar to either $\diag(1,j,j,j)$, or $\diag(j,j^2,j,j^2)$. 
\item Suppose finally that $u=6$ and $g$ has no eigenvalue of order 4: Then $g$ has order $6$, so $g^3$ has order $2$. Since $g^3$ acts freely in codimension $2$ and has determinant one, we obtain $g^3=-\id_B$. All eigenvalues of $g$ thus have order 2 or 6, so $g$ has at least three eigenvalues of order 6. As $g$ has determinant one, we only have two possibilities: There is a generator of $\grp{g}$ similar to $\diag(-1,\omega,\omega,\omega)$, or to $\diag(\omega,\omega^5,\omega,\omega^5)$.
\end{enumerate}
\end{enumerate}

This discussion constructs the first two columns of the table. The restrictions on $B$ given in the third column are given by the same arguments as in the proof of Lemmas \ref{lem-specialfactor}, \ref{lem-coollemma}.}

\cor[cor-allelements]{Let $B$ and $B'$ be isogenous abelian fourfolds, and let $g\in\mathrm{Aut}(B,0)$ and $h\in\mathrm{Aut}(B',0)$ be junior elements such that $\grp{g}$ and $\grp{h}$ act freely in codimension $2$, and $\mathrm{ord}(h)\le\mathrm{ord}(g)$. Then there are three possibilities:
\begin{itemize}
\item $g$ and $h$ are similar, in particular have the same order;
\item $h$ is similar to $\diag(1,j,j,j)$, $g$ is similar to $\diag(-1,\omega,\omega,\omega)$, and $B$ and $B'$ are isogenous to $E\times {E_j}^3$ for some elliptic curve $E$;
\item $h=i\id_B$, $g$ is similar to $\diag(\zeta_{12},\zeta_{12}^5,\zeta_{12},\zeta_{12}^{5})$ and $B$ and $B'$ are both isomorphic to ${E_i}^4$.
\end{itemize}
}

\demo{If $h$ has order $7$, then by Lemma \ref{lem-allelements}, $B$ is isogenous to $E\times {E_{u_7}}^3$ for some elliptic curve $E$. By uniqueness in the Poincaré decomposition of $B$ \cite[Thm.5.3.7]{BirkLang}, $B$ is not isogenous to any of the other special abelian varieties appearing in Lemma \ref{lem-allelements}. So, by Lemma \ref{lem-allelements} again, $g$ being junior must have order $7$. By Proposition \ref{prop-junclass}, any junior element $k$ of order $7$ acting on a fourfold with $\grp{k}$ acting freely in codimension $2$ are similar to $\diag(1,\zeta_7,{\zeta_7}^2,{\zeta_7}^4)$. So $g$ and $h$ are similar.

The same argument works if $h$ has order $8,15,16,20,24$.

If $h$ has order $3$ or $6$, then by Lemma \ref{lem-allelements}, $B$ is isogenous to $E\times {E_{j}}^3$ for some elliptic curve $E$. By uniqueness in the Poincaré decomposition of $B$ \cite[Thm.5.3.7]{BirkLang}, $B$ is not isogenous to any of the other special abelian varieties appearing in Lemma \ref{lem-allelements}. So, by Lemma \ref{lem-allelements} again, $g$ being junior must have order $3$ or $6$. As we assumed $\mathrm{ord}(h)\le\mathrm{ord}(g)$, the only strict inequality is when $h$ has order $3$ and $g$ has order $6$. In this case, by Proposition \ref{prop-junclass}, $h$ is similar to $\diag(1,j,j,j)$ and $g$ to $\diag(-1,\omega,\omega,\omega)$.

The same argument works if $h$ has order $4$ or $12$.}

We can now prove cyclicity of $F$ when it contains a junior element of order $3$.

\prop[prop-3]{Let $B$ be an abelian fourfold, and let $F$ be a finite subgroup of $\mathrm{Aut}(B,0)$ acting freely in codimension $2$, generated by junior elements. Suppose that $F$ contains an element similar to $\diag(1,j,j,j)$. Then $F$ is cyclic and generated by one junior element.}

\demo{By Corollary \ref{cor-allelements}, $B$ is isogenous to $E\times{E_j}^3$ for some elliptic curve $E$, and any junior element in $\mathrm{Aut}(B,0)$ is similar to $\diag(1,j,j,j)$, or $\diag(-1,\omega,\omega,\omega)$. 

Suppose by contradiction that $F$ is not generated by one junior element. Then there are two junior elements $g, h\in F$ such that $\grp{g}\nsubseteq\grp{h}$ and $\grp{h}\nsubseteq\grp{g}$. Up to possibly replacing them by their square, we have $\tilde{g}$ and $\tilde{h}$ both similar to $\diag(1,j,j,j)$. Their eigenspaces satisfy
$2\le \dim E_{\tilde{g}}(j)\cap E_{\tilde{h}}(j)\le \dim E_{\tilde{g}\tilde{h}^{-1}}(1).$
As $\grp{\tilde{g},\tilde{h}}\subset F$ acts freely in codimension $2$, $\tilde{g}=\tilde{h}$. 
Since $\grp{\tilde{g}}\subset\grp{h}$, $\tilde{g}\ne g$, so $\tilde{g}=g^2$. Similarly, $\tilde{h}=h^2$. Since $g^3=h^3=-\id$, it nonetheless yields $g=h$, contradiction.}

Let us now present our general strategy to prove that $F$ is cyclic.
By Lemma \ref{lem-allelements}, the prime divisors of $|F|$ are $2,3,5,$ and $7$.
Hence, $|F|=2^{\alpha}\cdot 3^{\beta}\cdot 5^{\gamma}\cdot 7^{\delta}$. Since $2^{\alpha}$ (respectively $3^{\beta}$, etc.) is the order of a $2$ (respectively $3$, etc.)-Sylow subgroup of $F$, we can rely on Sylow theory to bound $|F|$, as in the following result.

\prop[prop-card]{Let $B$ be an abelian fourfold, and let $F$ be a finite subgroup of $\mathrm{Aut}(B,0)$ acting freely in codimension $2$, generated by junior elements, containing no junior element of order 3. Then
$$|F|\mbox{ divides }2^4\cdot 3\cdot 5\cdot 7 =1680.$$
}

The proof of this proposition relies on the following two lemmas.

\lem[lem-5sylow]{Let $B$ be an abelian fourfold, and let $F$ be a finite subgroup of $\mathrm{Aut}(B,0)$ acting freely in codimension $2$ with determinant one, containing no junior element of order 3.
Let $p=3,5,$ or $7$ divide $|F|$. Then a $p$-Sylow subgroup $S$ of $F$ is cyclic of order $p$.}

\demo{As $S$ is a $p$-group, its center $Z(S)$ is non-trivial. Hence, it contains an element $g$ of order $p$. Let $h\ne\id\in S$. By Lemma \ref{lem-allelements}, $F$ has no element of order $p^2$, so $h$ has order $p$.
Since $g$ and $h$ commute, they are codiagonalizable. Let $v,w$ be two non-colinear common eigenvectors of them associated with eigenvalues other than 1. Let $\tilde{g}\in\grp{g}$ and $\tilde{h}\in\grp{h}$ satisfy $\tilde{g}(v)=\tilde{h}(v) =\zeta_pv$.

If $p=3$ or $5$, Lemma \ref{lem-allelements} shows that $E_g(1)=E_h(1)=\{0\}$, so $\tilde{g}\tilde{h}^{-1}$ can not have $1$ as an eigenvalue and be of order $p$. So it is trivial, i.e., $\tilde{g}=\tilde{h}$, and $h\in\grp{g}$.

Suppose $p=7$. If $\tilde{g}(w)\ne\tilde{h}(w)$, then by Lemma \ref{lem-allelements}, we have an equality $\{\tilde{g}(w),\tilde{h}(w)\}=\{{\zeta_7}^2w,{\zeta_7}^4w\}.$ So $\tilde{g}\tilde{h}^2$ has eigenvalue ${\zeta_7}^3$ on $v$, and ${\zeta_7}$ or ${\zeta_7}^3$ on $w$, which in either case contradicts Lemma \ref{lem-allelements}.  
So $\tilde{g}(w)=\tilde{h}(w)$, i.e., $\tilde{g}\tilde{h}^{-1}$ has eigenvalue 1 with multiplicity two. By freeness in codimension $2$, $\tilde{g}=\tilde{h}$, hence $h\in\grp{g}$.}

\lem[lem-2sylow]{Let $B$ be an abelian fourfold, and let $F$ be a finite subgroup of $\mathrm{Aut}(B,0)$ acting freely in codimension $2$ with determinant one. If not trivial, a $2$-Sylow subgroup $S$ of $F$ is cyclic or a generalized quaternion group, and its order divides 16.}

\demo{By Lemma \ref{lem-allelements}, the element of order 2 in $F$ is unique: it is $-\id_B$. By \cite[5.3.6]{Robinson}, $S$ is hence either cyclic or a generalized quaternion group.
Moreover, by Lemma \ref{lem-allelements}, $S$ has no element of order $32$.
Hence, the only case where the order of $S$ does not divide $16$, is when $S$ is isomorphic to $Q_{32}$. Let us however show that this is impossible.

Indeed, $Q_{32}$ contains an element $h$ of order 16 and an element $s$ of order $4$ such that $shs^{-1}=h^{-1}$ \cite[pp.140-141]{Robinson}. However, if $h\in S$ is an element of order 16, it can not be conjugated in $S$ to $h^{-1}$, because by Lemma \ref{lem-allelements} they have distinct eigenvalues.}

\demode{Proposition \ref{prop-card}.}{It is straightforward from Lemma \ref{lem-5sylow} and Lemma \ref{lem-2sylow}.}


The following Lemma and Proposition show that if $7$ divides $|F|$, i.e., if $F$ contains a junior element of order 7, then $F$ is cyclic generated by one junior element of order 7.

\lem[lem-syl7factors]{Let $B$ be an abelian fourfold, and let $F$ be a finite subgroup of $\mathrm{Aut}(B,0)$ acting freely in codimension $2$ with determinant one, containing no junior element of order $3$. Suppose that $7$ divide $|F|$. Let $S$ be a $7$-Sylow subgroup of $F$. Then there is a normal subgroup $N$ of $F$ such that $F=N\rtimes S$.}

\demo{By Burnside's normal complement theorem \cite[10.1.8]{Robinson}, it is enough to show that $N_F(S)=C_F(S)$.

Let $g$ be a generator of $S$, with $S= \grp{g}\simeq \Z_7$ by Lemma \ref{lem-5sylow}. By Lemma \ref{lem-allelements}, for an element $f\in N_F(S)$, we have $fgf^{-1}\in \{ g, g^2, g^4\}$, because they are the only elements with the same set of eigenvalues as $g$. Let us assume by contradiction that $N_F(S)\neq C_F(S)$, i.e., there exists $f\in N_F(S)$ such that $fgf^{-1}\neq g$. In particular, either $fgf^{-1}=g^2$, or $fgf^{-1}=g^4$, in which case $f^2gf^{-2}=(g^{4})^4=g^2$. Let $\tilde{f}$ be the element of $\{f,f^2\}$ such that $\tilde{f}g\tilde{f}^{-1}=g^2$.
Looking at the action of $\tilde{f}$ on the eigenspaces of $g$ in coordinates diagonalizing $g$, we see that
$$\tilde{f}=\left(\begin{array}{cccc}
t & 0 & 0 & 0\\
0 & 0 & z & 0\\
0 & 0 & 0 & y\\
0 & x & 0 & 0\\
\end{array}\right),$$
for some complex numbers $x,y,z,t$ with $xyzt=1$, and so $\chi_{\tilde{f}}=(X-t)(X^3-t^{-1})$. Note that $\tilde{f}$ has order $3k$ for some integer $k\ge 1$, and that $g\tilde{f}^3$ has order the least common multiple of $k$ and $7$, since $g$ and $\tilde{f}^3$ commute.
But by Lemma \ref{lem-allelements} for the element $g\tilde{f}^3$, this yields $k\in\{1,2,7,14\}$. By Lemma \ref{lem-allelements} for $\tilde{f}$, this yields $3k\in\{3,6\}$. Since by assumption $F$ contains no junior element of order 3, we obtain by Lemma \ref{lem-allelements} that $\tilde{f}$ has characteristic polynomial $(X^2+X+1)^2$ or $(X^2-X+1)^2$, contradiction.}

\prop[prop-7]{Let $B$ be an abelian fourfold, and let $F$ be a finite subgroup of $\mathrm{Aut}(B,0)$ acting freely in codimension $2$, generated by junior elements, containing no junior element of order 3. Suppose that $7$ divides $|F|$. Then $F$ is cyclic and generated by one junior element.}

\demo{Let $S$ be a $7$-Sylow subgroup of $F$. By Lemma \ref{lem-syl7factors}, $F=N\rtimes S$, where $N$ is a normal subgroup of $F$, and by Proposition \ref{prop-card}, $|N|$ divides 240. A simple \texttt{GAP} program in the appendix checks that a group of order dividing 240 cannot have an automorphism of order 7. So $S$ acts trivially on $N$, i.e., $F=N\times S$. But $F$ is generated by its junior elements, which all have order 7 by Corollary \ref{cor-allelements}. So $N$ is trivial, and $F=S$ is cyclic of order 7.}

Now we can focus on the case when $F$ contains no junior element of order 3 or 7. 
We start by showing that, provided $F$ is cyclic, it is generated by one junior element.

\lem[lem-cycliconegen]{Let $F$ be a cyclic group. If $E$ is a subset of $F$ that generates $F$, and all elements of $E$ have the same order, then any one element of $E$ actually generates $F$.}

\demo{Suppose that $F=\Z_d$ and that every element of $E$ has order $k$ dividing $d$. Then $E$ is actually a subset of the subgroup $\Z_{k} < \Z_d$. Since $E$ generates $\Z_d$, it must be $k=d$. So any element $e\in E$ satisfies $\grp{e}=\Z_d=F$.}

\cor[cor-cycliconegen]{Let $B$ be an abelian fourfold, and let $F$ be a finite subgroup of $\mathrm{Aut}(B,0)$ acting freely in codimension $2$, generated by junior elements, containing no junior element of order $3$ or $7$. If $F$ is cyclic, then $F$ is generated by one junior element.}

\demo{Assume that $F$ is cyclic. If $F$ contains one junior element of order $8,15,16,20$, or $24$, then by Corollary \ref{cor-allelements}, all junior elements have the same order and we use Lemma \ref{lem-cycliconegen} to conclude.

Otherwise, the junior elements of $F$ each have order $4$ or $12$. If there are no junior elements of order 12, Lemma \ref{lem-cycliconegen} concludes again. If there is a junior element $g$ of order 12, then a quick computation from Lemma \ref{lem-allelements} shows that $g^3$ is the only junior element of order 4 in $F$, and thus the set of the junior elements of order 12 actually generates $F$ too, so we conclude by Lemma \ref{lem-cycliconegen} again.}

These versions of Lemma \ref{lem-syl7factors} for 3- and 5-Sylow subgroups will be useful too.

\lem[lem-nc5]{Let $B$ be an abelian fourfold, and let $F$ be a finite subgroup of $\mathrm{Aut}(B,0)$ acting freely in codimension $2$, generated by junior elements. Suppose that $p\in\{3,5\}$ divides $|F|$. Let $S$ be a $p$-Sylow subgroup of $F$. Then $N_F(S)/C_F(S)$ is isomorphic to a subgroup of $(\Z_p)^{\times}$.}

\demo{The quotient $N_F(S)/C_F(S)$ acts faithfully by conjugation on $S$, and therefore embeds in $\mathrm{Aut}(S)$, which by Lemma \ref{lem-5sylow} is isomorphic to $(\Z_p)^{\times}$.}

\lem[lem-nc5bis]{Let $B$ be an abelian fourfold, and let $F$ be a finite subgroup of $\mathrm{Aut}(B,0)$ acting freely in codimension $2$, generated by junior elements. Suppose that $5$ divides $|F|$. Let $S$ be a $5$-Sylow subgroup of $F$. Then, if $f\in N_F(S)$ is a junior element of order $8$, $[f]\in N_F(S)/C_F(S)$ cannot have order $4$.}

\demo{Let $f\in N_F(S)$ be a junior element of order $8$ such that $[f]\in N_F(S)/C_F(S)$ has order $4$, and let $g$ be a generator of $S$.
Looking at the action of $f$ on the eigenspaces of $g$ in coordinates diagonalizing $g$,
$$f=\left(\begin{array}{cccc}
0 & 0 & 0 & t\\
x & 0 & 0 & 0\\
0 & y & 0 & 0\\
0 & 0 & z & 0\\
\end{array}\right),$$
with $xyzt=-1$, and so $\chi_f=X^4+1$. By Lemma \ref{lem-allelements}, no junior element of order $8$ has this characteristic polynomial, contradiction.}

We finally prove the following two key propositions, which imply Proposition \ref{prop-justshowFcyc}.

\prop[prop-2sylowiscyc]{Let $B$ be an abelian fourfold, and let $F$ be a finite subgroup of $\mathrm{Aut}(B,0)$ acting freely in codimension $2$, generated by junior elements, containing no junior element of order $3$ or $7$. Then a $2$-Sylow subgroup of $F$ is either trivial, or cyclic.}

\prop[prop-if2sylowcyc]{Let $B$ be an abelian fourfold, and let $F$ be a finite subgroup of $\mathrm{Aut}(B,0)$ acting freely in codimension $2$, generated by junior elements, containing no junior element of order $3$ or $7$. Suppose that a $2$-Sylow subgroup of $F$ is trivial or cyclic. Then $F$ is cyclic.}

\demode{Proposition \ref{prop-if2sylowcyc}}{Let us write $|F|=2^{\alpha}\cdot 3^{\beta}\cdot 5^{\gamma}$ with $\alpha\in\lint 0, 4\rint$, $\beta,\gamma\in \lint 0,1\rint$. By Lemma \ref{lem-5sylow} and by assumption, the Sylow subgroups of $F$ are cyclic, so \cite[pp.290-291]{Robinson} applies and $F$ is a semidirect product:
$F\simeq (\Z_{5^{\gamma}}\rtimes\Z_{3^{\beta}})\rtimes \Z_{2^{\alpha}}.$
Since $3^{\beta}$ is coprime to $\phi(5^{\gamma})$, the group $\Z_{5^{\gamma}}$ has no automorphism of order 3, and thus the first semidirect product is direct:
$$F\simeq (\Z_{5^{\gamma}}\times\Z_{3^{\beta}})\rtimes \Z_{2^{\alpha}}.$$

If $\beta=\gamma = 1$, the group $F$ contains an element of order $15$, so by Lemma \ref{lem-allelements}, $B$ is isomorphic to ${E_{u_{15}}}^4$ and all junior elements of $F$ have order $15$.
However, since $F\simeq \Z_{15}\rtimes \Z_{2^{\alpha}}$, and since $F$ is generated by its junior elements, we must have $\alpha = 0$, and so $F\simeq \Z_{15}$ is cyclic and generated by one junior element.

If $\beta =\gamma = 0$, then $F\simeq\Z_{2^{\alpha}}$ is cyclic.

Else, write $p=3^{\beta}5^{\gamma}$ and $F\simeq \Z_{p}\rtimes\Z_{2^{\alpha}}$. Note that $\Z_p\rtimes\Z_{2^{\alpha-1}}$ is a proper subgroup of $F$ containing all elements whose order divides $2^{\alpha -1}p$. As $F$ is generated by its junior elements, their orders cannot all divide $2^{\alpha -1}p$: There is a junior element $g\in F$ of order $2^{\alpha}$ or $2^{\alpha}p$. If $g$ has order $2^{\alpha}p$, $\grp{g}=F$ and so $F$ is cyclic.
If $g$ has order $2^{\alpha}$, we can write $F\simeq \grp{u}\rtimes\grp{g}$, where $u$ is an element of $F$ of order $p$. The discussion now depends on $\alpha$ and $p$.
\begin{enumerate}[label= (\arabic*)]
\item By Lemma \ref{lem-allelements}, if $g$ has order $4$, then $g=i\id$ commutes with every element of $F$, so the semidirect product is direct and $F$ is cyclic.
\item If $p=5$ and $g$ has order $8$, by Lemma \ref{lem-nc5bis}, $g^2$ and $u$ commute, so $g^2u$ has order 20. Since $g$ is junior of order 8, by Lemma \ref{lem-allelements}, $B$ is isomorphic to ${E_{u_8}}^4$. So by Lemma \ref{lem-allelements} again, $B$ has no automorphism of order 20, contradiction.
\item If $p=5$ and $g$ has order 16, by Lemma \ref{lem-nc5}, $g^4$ and $u$ commute, so $g^4u$ has order 20. But since $g$ is junior of order 16, by Lemma \ref{lem-allelements}, $B$ has no automorphism of order 20, contradiction. 
\item If $p=3$ and $g$ has order 16, by Lemma \ref{lem-nc5}, $g^2u$ has order 24. But since $g$ is junior of order 16, by Lemma \ref{lem-allelements}, $B$ has no automorphism of order 24, contradiction.
\item If $p=3$ and $g$ has order 8, then $F\simeq\Z_3\rtimes\Z_8$. With \texttt{GAP}, we check in the Appendix that:
\begin{itemize}
\item The irreducible representations of $F$ have rank 1 or 2.
\item No irreducible character of $F$ takes value $j$ or $j^2$, so $F\subset\mathrm{Aut}(B,0)$ has no irreducible subrepresentation of rank 1.
\item The only two irreducible representations of $F$ of rank 2 sending $-\id\in F$ to $-\id$ indeed are complex conjugates, so all elements of $F\subset\mathrm{Aut}(B,0)$ have characteristic polynomials in $\Q[X]$.
\end{itemize}
However, $g\in F$ is a junior element of order 8, which by Lemma \ref{lem-allelements} has a non-rational characteristic polynomial, contradiction.
\end{enumerate}}

We prove Proposition \ref{prop-2sylowiscyc} by contradiction.

\demode{Proposition \ref{prop-2sylowiscyc}}{Suppose that $2$ divides $|F|$ and that a $2$-Sylow subgroup of $F$ is not cyclic. We first show that any junior element in $F$ has order $15$, $20$ or $24$.

By contradiction and by Proposition \ref{prop-junclass}, let $g\in F$ be a junior element of order $4,8,12,$ or $16$. If $g$ has order $12$, then $g^3\in F$ is a junior element of order $4$, and $F$ thus contains a junior element $\tilde{g}$ of order $4,8,$ or $16$. Let $S$ be a $2$-Sylow subgroup containing that junior element. By assumption, $S$ is not cyclic, so by Lemma \ref{lem-2sylow}, $S$ is isomorphic to $Q_8$ or to $Q_{16}$.
Clearly, $Q_8$ and $Q_{16}$ have no element of order 16, and no element of order 4 in their centers, so $\tilde{g}$ has order $8$. As $Q_8$ has no element of order $8$, $S$ is isomorphic to $Q_{16}$. But we easily check with \texttt{GAP} that:
\begin{itemize}
\item The irreducible representations of $Q_{16}$ have rank 1 or 2.
\item The only irreducible representations of $Q_{16}$ of rank $r$ sending the unique element of order $2$ to $-\id_r$ are two complex conjugates representations with $r=2$, so all elements of $S\subset\mathrm{Aut}(B,0)$ have characteristic polynomials in $\Q[X]$.
\end{itemize}
However, $\tilde{g}\in S$ is a junior element of order 8, which by Lemma \ref{lem-allelements} has a non-rational characteristic polynomial, contradiction. 

So any junior element in $F$ has order $15$, $20$ or $24$.
We also know that:
\begin{itemize}
\item $F$ has exactly one element of order 2, by Lemma \ref{lem-allelements}.
\item A $2$-Sylow subgroup of $F$ is isomorphic to $Q_8$ or $Q_{16}$, by Lemma \ref{lem-2sylow}.
\item $|F|$ divides $240$, by Proposition \ref{prop-card}.
\item $F$ has no element of order 60 or 40, by Lemma \ref{lem-allelements}.
\item If $F$ has elements of orders $o,o'\in\{15,20,24\}$, then $o=o'$, by Lemma \ref{lem-allelements}. 
\end{itemize}
We check with \texttt{GAP} that there are only five groups satisfying all these properties, namely the groups indexed \texttt{(40,4)},\texttt{(40,11)},\texttt{(80,18)},\texttt{(48,8)}, and \texttt{(48,27)} in the \texttt{SmallGroup} library. The function \texttt{StructureDescription} then shows that they are respectively of the form 
$\Z_5\rtimes Q_8$, $\Z_5\times Q_8$, $\Z_5\rtimes Q_{16}$, $\Z_3\rtimes Q_{16}$, and $\Z_3\times Q_{16}.$ Note that only $\Z_5\times Q_8$, $\Z_5\rtimes Q_{16}$ are generated indeed by their elements of orders ($15,24$, or) $20$. Checking the irreducible character tables of these two cadidates with \texttt{GAP} shows that they have no appropriate four-dimensional representation (see Appendix for programs supporting this discussion.)

This concludes the proof of Proposition \ref{prop-2sylowiscyc}.}

\demode{Proposition \ref{prop-justshowFcyc}}{If $F$ contains a junior element of order $3$, then Proposition \ref{prop-3} applies and shows that $F$ is cyclic generated by one junior element. If $F$ contains no junior element of order $3$, but one of order $7$, then Proposition \ref{prop-7} applies and shows that $F$ is cyclic generated by one junior element. Finally, if $F$ contains no junior element of order $3$ or $7$, Proposition \ref{prop-2sylowiscyc} shows that its $2$-Sylow subgroups are cyclic or trivial, Proposition \ref{prop-if2sylowcyc} deduces that $F$ is cyclic and Corollary \ref{cor-cycliconegen} proves that $F$ is generated by one junior element.}

\section{Ruling out junior elements in codimension 4}\label{sec-RuleOutMostJuniors}

The aim of this section is to rule out eight out of the twelve types of junior elements presented in Proposition \ref{prop-junclass}, namely those which fix  pointwise at least one subvariety of codimension 4, but no subvariety of codimension 3.

\prop[prop-nojun37]{Let $A$ be an abelian variety of dimension $n$, $G$ a group acting freely in codimension $2$ on $A$ such that $A/G$ has a $K$-trivial resolution. Then, if $g\in G$ is a junior element, the matrix $M(g)$ cannot have eigenvalue $1$ with multiplicity exactly $n-4$.}

\rem{Whether the local affine quotients corresponding to these eight types of junior elements admit a crepant resolution is actually settled by toric geometry in \cite{Sato}. In fact, by \cite[Thm.3.1]{Sato},
$$\C^4/\grp{i\id},\quad
\C^4/\grp{\diag{(\omega,\omega,\omega,-1)}},\quad
\C^4/\grp{\diag(\zeta_8,\zeta_8,\zeta_8^3,\zeta_8^3)},$$
$$\C^4/\grp{\diag(\zeta_{12},\zeta_{12},\zeta_{12}^5,\zeta_{12}^5)},\quad
\C^4/\grp{\diag(\zeta_{15},\zeta_{15}^2,\zeta_{15}^4,\zeta_{15}^8)}$$
have a crepant Fujiki-Oka resolution, and by \cite[Prop.3.9]{Sato},
$$\C^4/\grp{\diag(\zeta_{16},\zeta_{16}^3,\zeta_{16}^5,\zeta_{16}^9)},\;
\C^4/\grp{\diag(\zeta_{20},\zeta_{20}^3,\zeta_{20}^7,\zeta_{20}^9)},\;
\C^4/\grp{\diag(\zeta_{24},\zeta_{24}^5,\zeta_{24}^7,\zeta_{15}^{11})}$$
admit no toric crepant resolution. They could nevertheless have a non-toric crepant resolution.}

In light of this remark, the proof of Proposition \ref{prop-nojun37} must crucially involve global arguments.

\subsection{Ruling out junior elements of order 4, 8, 12, 15, 16, 20, 24}\label{subsec-RuleOutSevenJuniors}

In this subsection, we rule out the seven types of junior elements or order other than $3,$ $6,$ $7$.

\prop[prop-globfailjun]{Let $A$ be an abelian variety, $G$ a group acting freely in codimension $2$ on $A$ such that $A/G$ has a $K$-trivial resolution $X$. Then any junior element of $G$ has order $3,$ $6,$ or $7$.}

\rem[rem-notrans]{Let $A$ be an abelian variety, $G$ be a group acting freely in codimension $2$ on $A$. As translations in $G$ form a normal subgroup $G_0$, we can write:
$$(A/G_0)/(G/G_0)\simeq A/G.$$ 
Clearly, $A/G_0$ is isogenous to $A$ and $G/G_0$ still acts freely in codimension $2$ on it, except that it contains no translation.
Hence, we can assume without loss of generality that $G$ contains no translation (other than $\id$). In particular, any element of $G$ has the same finite order as its matrix.}

\demode{Proposition \ref{prop-globfailjun}.}{By contradiction, suppose that $g\in G$ is a junior element of order $d\in\{4,8,12,15,16,20,24\}$, of maximal order among the junior elements of $G$ of such orders. Up to conjugating the whole group $G$ by an appropriate translation, we may assume that $g$ fixes $0\in A$. In particular, by Proposition \ref{prop-junclass}, $g$ fixes pointwise an abelian subvariety $W$ of $A$ of codimension 4, so Proposition \ref{prop-cyclicincod4} shows that the pointwise stabilizer $\pstab(W)$, which contains $\grp{g}$, is cyclic and generated by a junior element $h\in G$. By our maximality assumption on the order of $g$, we have $\grp{g} =\grp{h}$, i.e., $\pstab(W)=\grp{g}$. By Proposition \ref{prop-point-wise-stabilizer}, we can define a $\grp{g}$-stable complementary abelian subvariety $B$ to $W$ in $A$.
The key to the proof is that a well-chosen power $g^{\alpha}$ of $g$ has strictly more fixed points in $B$ than $g$, as many distinct eigenvalues as $g$, but is not a junior element. We set the integer $\alpha$ depending on the order $d$ of $g$ as follows, check with Proposition \ref{prop-junclass} that $g^{\alpha}$ is not junior and has as many distinct eigenvalues as $g$. Applying \cite[Prop.13.2.5(c)]{BirkLang} shows that $(g^{\alpha})|_B$ has strictly more fixed points than $g|_B$ in $B$.

\begin{table}[H]
$$\begin{array}{|c|c|c|c|c|c|c|c|}
\hline
d & 4 & 8 & 12 & 15 & 16 & 20 & 24\\
\hline
\alpha & 2 & 2 & 4 & 3 & 2 & 4 & 3\\
\hline
\end{array}$$
\caption{Definition of a certain $\alpha\in\lint 0,d-1\rint$ depending on $d$}
\end{table}

Let $\tau\in B$ be a fixed point of $g^{\alpha}$ that is not fixed by $g$. Note that $W+\tau$ is pointwise fixed by $g^{\alpha}$. By Proposition \ref{prop-cyclicincod4}, $\pstab(W+\tau)=\grp{h}$ for some junior element $h$.

By Proposition \ref{prop-point-wise-stabilizer}, there is an $\grp{h}$-stable translated abelian subvariety $B'$ of $A$ containing $\tau$ such that $B'$ and $W+\tau$ are complementary. By uniqueness in Poincaré's complete reducibility theorem \cite[Thm.5.3.7]{BirkLang}, the abelian varieties $B$ and $B'$ are isogenous. Hence, the relationship between the orders and similarity classes of $g$ and $h$ is described by Corollary \ref{cor-allelements}.

Let us discuss the special case when $g$ and $h$ do not have the same order. By Corollary \ref{cor-allelements}, and since $g$ does not have order $3$ or $6$ by assumption, $g$ and $h$ must have orders in $\{4,12\}$. By the maximality assumption on the order of $g$, we have $g$ of order $12$, hence $g^{\alpha}$ of order $3$, and $h$ of order $4$. That contradicts the fact that $g^{\alpha}\in\pstab(W+\tau)=\grp{h}$.

We can now assume that $g$ and $h$ have the same order $d\in\{4,8,12,15,16,20,24\}$, and similar matrices.
Recall that $g^{\alpha}\in\grp{h}$. Since $g$ and $h$ have the same order, it implies $\grp{g^{\alpha}}=\grp{h^{\alpha}}$, i.e., $g^{\alpha}=h^{u\alpha}$ for some $u$ coprime to $\frac{d}{\alpha}$. Since $g$ and $g^{\alpha}$, and $h$ and $h^{u\alpha}$ have the same number of distinct eigenvalues, it follows from $g^{\alpha}=h^{u\alpha}$ that the eigenspaces of $g$ and $h$ are the same, i.e., $g$ and $h$ commute. Since we assumed that $g$ does not fix $\tau$, we still have that $g\not\in\grp{h}$. We discuss three cases separately.

\begin{enumerate}[label={(\arabic*)}]
\item If $d=4$, then in appropriate coordinates, we have:
\begin{align*}
M(g)&=\diag(\1_{n-4},i,i,i,i)\\
M(h)&=\diag(\1_{n-4},i,i,i,i)
\end{align*}
so $g=h$, contradiction!
\item If $d=8$ or $12$, then in appropriate coordinates, we have:
\begin{align*}
M(g)&=\diag(\1_{n-4},\zeta_d,\zeta_d,{\zeta_d}^m,{\zeta_d}^m)\\
M(h)&=\diag(\1_{n-4},{\zeta_d}^{m},{\zeta_d}^{m},\zeta_d,\zeta_d)
\end{align*}
for some integer $m\in\lint 2,d-1\rint$ such that $2+2m = d$. In particular, $m^2\equiv 1\;\mathrm{mod}\, d$, so $g=h^m\in \grp{h}$, contradiction!
\item Else, $d=15,$ $16,$ $20,$ or $24$. There is an integer $u'$ coprime to $d$ such that, in appropriate coordinates,
\begin{align*}
M(g)&=\diag(\1_{n-4},\zeta_d,{\zeta_d}^a,{\zeta_d}^b,{\zeta_d}^c)\\
M(h^{u'})&=\diag(\1_{n-4},\zeta_d,{\zeta_d}^{\sigma(a)},{\zeta_d}^{\sigma(b)},{\zeta_d}^{\sigma(c)})
\end{align*}
for some distinct integers $a,b,c\in\lint 2,d-1\rint$ coprime to $d$, and permutation $\sigma$ of $\{a,b,c\}$. If $\sigma=\id$, then $g=h^{u'}\in\grp{h}$, contradiction!  Nevertheless, let us prove that $\sigma=\id$.
Note that 
\begin{align*}
(h^{u'-u})^{\alpha}
&=(h^{u'}g^{-1})^{\alpha}(h^{-u}g)^{\alpha}\\
&=\diag(\1_{n-3},{\zeta_d}^{(\sigma(a)-a)\alpha},{\zeta_d}^{(\sigma(b)-b)\alpha},{\zeta_d}^{(\sigma(c)-c)\alpha}),
\end{align*}
and thus $(h^{u'-u})^{\alpha}$ fixes a translated abelian variety $W'\supset W+\tau$ of codimension at most 3. By Proposition \ref{prop-cyclicincod4}, $\pstab(W')$ is trivial, or cyclic and generated by one junior element $k$ of order $3$ or $7$. In the second case, as $k\in\pstab(W+\tau)$, $k$ restricts to an automorphism of the fourfold $B'$, which also has $h$ junior of order $d\ne 3,6,7$ acting on it. This contradicts Corollary \ref{cor-allelements}. Hence, $(h^{u'-u})^{\alpha}\in \pstab(W')=\{\id\}$, so for any $\ell\in\{a,b,c\}$, $(\sigma(\ell)-\ell)\alpha$ is a multiple of $d$.
However, $\alpha$ was chosen so that $g^{\alpha}$ and $g$ have the same number of distinct eigenvalues, i.e., $a\alpha,b\alpha,c\alpha$ are distinct modulo $d$. In particular, $\sigma(\ell)\alpha=\ell\alpha$ modulo $d$ if and only if $\sigma(\ell)=\ell$. So $\sigma=\id$, contradiction!
\end{enumerate}}

\subsection{Ruling out junior elements of order 6 with four non-trivial eigenvalues}\label{subsec-ruleout64}

In this subsection, we conclude the proof of Proposition \ref{prop-nojun37} by ruling out the one remaining type of junior element fixing at least one subvariety of codimension 4, but no subvariety of codimension 3. It is the type of junior element of order 6, and matrix similar to $\diag(\1_{n-4}, \omega,\omega,\omega,-1).$

\prop[prop-ruleout6cod4]{Let $A$ be an abelian variety, $G$ a group acting freely in codimension $2$ on $A$ such that $A/G$ has a $K$-trivial resolution. Then there is no junior element of $G$ with matrix similar to 
$\diag(\1_{n-4},\omega,\omega,\omega, -1).$}

The proof involves general arguments which we will use later, hence we factor it into a general lemma.

\lem[lem-generalGW]{Let $A$ be an abelian variety of dimension $n$, $G$ a group acting freely in codimension $2$ on $A$, without translations. Assume moreover that $G$ is generated by junior elements, and that for any point $a\in A$, the stabilizer subgroup $\pstab(a)$ of $G$ is trivial, or generated by junior elements. Suppose that $g\in G$ fixes $0\in A$ and has order $d$. Let $W$ be the abelian subvariety in $A$ containing $0$ that $g$ fixes pointwise, and denote by $G_W$ the subgroup of $G$ generated by 
$$G_{\mathrm{gen}} = {G_{\mathrm{gen}}}^{-1} = \left\{h\in G\mid \exists\,\tau\in A\mbox{ such that } h\in \pstab(W+\tau)\right\}.$$ 
Then

\begin{enumerate}[label =\textup{(\arabic*)}]
\item There is an $M(G_W)$-stable complementary abelian subvariety $B$ to $W$, which induces a representation $\rho:G_W\rightarrow \mathrm{Aut}(B,0)$ by $\rho(h):=M(h)|_B$.
\item If we denote by $pr_W$, $pr_B$ the projections induced by the splitting of the tangent space, then, for any $h\in G_W$, 
its matrix and translation parts satisfy
\begin{itemize}
\item $M(h)= pr_W + \rho(h) pr_B$
\item $pr_W(T(h))=0$, i.e., $T(h)\in B$
\end{itemize}
\item The representation $\rho$ is faithful and takes values in $\SL(H^0(T_B))$. 
\item The abelian subvariety $B$ is in fact $G_W$-stable. 
\item For any element $h\in G_W$, 
\begin{itemize}
\item if $h$ fixes a point $\tau\in A$, then it fixes the point $pr_B(\tau)\in B$;
\item if $h$ has no fixed point in $A$, then $1$ is an eigenvalue of $\rho(h)$.
\end{itemize}
\end{enumerate}
Moreover, if we also assume that the codimension $k$ of $W$ in $A$ is at most $7$, and that there is an integer $1\le \alpha\le d-1$ such that $M(g^{\alpha})$ is similar to 
$\diag(\1_{n-k},-\1_k)$, then
\begin{enumerate}[label =\textup{(\arabic*)}]
\setcounter{enumi}{5}
\item For any $h\in G_W$, the elements $h$ and $g^\alpha$ commute and
\begin{itemize}
\item either there is a point $\tau\in A$ such that $h\in\pstab(W+\tau)\cup g^{\alpha}\pstab(W+\tau)$;
\item or there is no such point, and $1$ and $-1$ are eigenvalues of $\rho(h)$.
\end{itemize}
\item For any $h\in G_W$, the translation part $T(h)$ of $h$ is a $2$-torsion point of $B$. 
\item If $h\in G_W$ has even order and fixes a point in $A$, all fixed points of $h$ in $B$ are of $2$-torsion. 
\item If $h\in G_W$ is a junior element of order 3, then $h$ fixes a non-zero $2$-torsion point in $B$.
\end{enumerate}
}

\demo{Item (1) follows immediately from \cite[Prop.13.5.1]{BirkLang}, since $M(G_W)$ is a finite group of group automorphisms of the abelian variety $A$, and $W$ is $M(G_W)$-stable.

\medskip

Item (2) is proven by induction on the number of generators used to write an element $h\in G_W$. First, if $h\in G_W$ is in $G_{\mathrm{gen}}$, there is a point $\tau\in A$ such that $h\in\pstab(W+\tau)$. In particular, for $w\in W$ and $b\in B$;
$$M(h)(w+b)=h(w+\tau)-h(\tau)+M(h)(b)=w+\rho(h)(b),$$ 
as wished.
Moreover, $T(h)=(\id-M(h))(\tau)$, so $pr_W(T(h))=0$.

Second, if $h_1,h_2\in G_W$ satisfy (2), then 
$$M(h_1h_2)=M(h_1)M(h_2) = pr_W + \rho(h_1h_2) pr_B,$$
since $\rho$ is a group morphism and $pr_W pr_B=pr_B pr_W=0$.
Moreover, $T(h_1h_2)=T(h_1)+M(h_1)T(h_2)$, and the fact that $pr_W(T(h_1h_2))=0$ easily follows from the induction assumption, notably using $pr_W(\id-M(h_1))=0$.

\medskip

For (3), let $h\in G_W$ and note that $\rho(h)=\id_B$ if and only if $M(h)=pr_W+pr_B=\id_A$, so $\rho$ is faithful if and only if $M$ is, which it is since we assumed that $G$ acts on $A$ without translations. Note that the set $M(G_{\rm gen})$ is contained in $\SL(H^0(T_A))$. Hence, the generated group $M(G_W)$ is contained in $\SL(H^0(T_A))$ too. By the matrix description in Item (2), $\rho(G_W)$ is contained in $\SL(H^0(T_B))$.

\medskip

Regarding (4) we note that, for $h\in G_W$, $h(B)=M(h)B+T(h)= B+T(h)=B$ by the translation description in Item (2).

\medskip

For (5), the fact that the identity $h(\tau)=\tau$ implies $h(pr_B(\tau))=pr_B(\tau)$ is clear from the description of the matrix part in (2). That concludes if $h$ has a fixed point. Assume now that $h$ has no fixed point in $A$. In other words, the element $T(h)$ does not belong to $\mathrm{Im}(\id-M(h))$. By the description of the translation part in Item (2), $T(h)$ belongs to $B$. Hence, the image $\mathrm{Im}(\id_B-M(h)|_B)$ is strictly contained in $B$. So $1$ appears among the eigenvalues of $\rho(h)=M(h)|_B$, as wished.

\medskip

We now prove (6). Note that $\rho(g^{\alpha})=-\id_B$ commutes with any element of $\rho(G_W)$, and thus, as $\rho$ is faithful by Item (3), $g^\alpha$ is in the center of $G_W$. 
Let $h\in G_W$ such that neither $h$, nor $g^{\alpha}h$ has any fixed point in $A$. By Item (5), this implies that $1$ appears among the eigenvalues of both $\rho(h)$ and $\rho(g^{\alpha}h)$. So $1$ and $-1$ appear among the eigenvalues of $\rho(h)$, as wished.
Now if $h$ or $g^{\alpha}h$ has a fixed point $\tau$ in $A$, then it must fix $W+\tau$ pointwise, by the matrix description in Item (2), as wished.

\medskip

For (7), we use that $h$ commutes with $g^{\alpha}$ by Item (6), that $g(0)=0$, that $h(0)=T(h)\in B$ by the translation description in Item (2), and that $g^{\alpha}|_B = -\id_B$. It yields 
$$0 = g^{\alpha}(h(0)) - h(g^{\alpha}(0))=g^{\alpha}(T(h))-T(h) = -2T(h),$$ 
so $T(h)$ is of $2$-torsion. 

\medskip

To prove (8), we have an element $h\in G_W$ of even order, that fixes a point $\tau$ in $A$. For some positive integer $\beta$, the power $h^{\beta}$ has order $2$. Since the whole group $G$ acts freely in codimension 2 on $A$, each of the matrices $\rho(h^{\beta})$ and $\rho(g^{\alpha}h^{\beta})$ is either trivial, or has $-1$ as an eigenvalue of multiplicity at least $3$, and in fact at least $4$ since $\rho(G_W)$ is contained in $\SL(H^0(T_B))$ by Item (3). Since $\rho(g^{\alpha})=-\id_B$, this means that either $\rho(g^{\alpha}h^{\beta})=\id_B$, or that $k\ge 8$, a contradiction with our assumption. Since $\rho$ is faithful by Item (3), we obtain $h^{\beta}=g^{\alpha}$. Hence, every point of $B$ fixed by $h$ is also fixed by $g^{\alpha}$, hence is a $2$-torsion point of $B$.

\medskip

To prove (9), we have a junior element $h\in G_W$ of order $3$. By Item (5), it has a fixed point $\tau\in B$. Since $M(h)$ is similar to ${\rm diag}(\1_{n-3},j,j,j)$, we note that $h$ pointwise fixes a translated abelian subvariety of codimension $3$ of the form $V+\tau$ in $A$. Let $C$ be a $\grp{M(h)}$-stable complementary to the abelian subvariety $V\cap B$ in the abelian subvariety $B$. We write $\tau = v + c$, with $v\in V\cap B$ and $c\in C$. It gives
$h(c)=h(\tau-v)=\tau-v=c$, i.e., $h$ fixes the point $c\in C$. Moreover, since $M(h|_{C})=j\id_{C}$, where $j$ is the third root of unity, it holds
$$0=h(c) - c = (j-1)c + T(h).$$
It shows that $T(h)\in C\simeq {E_j}^3$. Multiplying this equality by $j^2-1$, and using that $T(h)$ is of $2$-torsion by Item (7), we see that $3c$ is of $2$-torsion too. Since $h(c)=c$ and $3T(h)=T(h)$, this point $3c$ in $C\subset B$ is fixed by $h$. If $3c$ is non-zero, this concludes the proof of (9). Assume now that $3c=0$, i.e., $T(h)=j^2T(h)$ in $C$. Since $T(h)$ is also $2$-torsion in $C$, this implies $T(h)=0$, and so $h$ fixes $0$ in $A$, and pointwise fixes the abelian subvariety $V\cap B$ of dimension $3$. It clearly contains many non-zero $2$-torsion points, and any of them is fixed by $h$, as wished.}

We can now come back to our proposition.

\demode{Proposition \ref{prop-ruleout6cod4}.}{By Remark \ref{rem-notrans}, we can assume that $G$ contains no translation other than $\id_A$. By contradiction, suppose that there is an element $g\in G$ such that $g(0)=0$ and, in some coordinates, 
$$M(g) = \diag(\1_{n-4}, \omega,\omega,\omega,-1).$$

We import the notations of Lemma \ref{lem-generalGW}, whose hypotheses are satisfied by $g$ for $k=4$, $d=6$, $\alpha = 3$. 
The proof of the proposition now goes in three steps. First, we show that every element of $\rho(G_W)$ is similar to an element of $\grp{\rho(g)}\simeq\grp{\diag(\omega,\omega,\omega,-1)}$. Second, we deduce that $G_W=\grp{g}$. Third, we use global considerations on fixed loci to derive a contradiction from this description of $G_W$.

\step{1}{We show that every element of $\rho(G_W)$ is similar to an element of $\grp{\rho(g)}$.

\medskip

By Lemma \ref{lem-generalGW} (1) and (4), there is a $G_W$-stable complementary $B$ to $W$. As $\rho(g)$ acts on it, $B$ is isogenous to $E\times {E_j}^3$ for some elliptic curve $E$. Fix an element $\tau\in A$ such that $\pstab(W+\tau)$ is non-trivial. By Proposition \ref{prop-cyclicincod4}, the group $\pstab(W+\tau)$ is cyclic generated by one junior element $u$, and by Corollary \ref{cor-allelements}, $\rho(u)$ is similar to $\rho(g)$ (if of order 6) or to $\rho(g^2)$ (if of order 3) in $\GL(H^0(T_B))$. Hence, any element of $\rho(\pstab(W+\tau))$ is similar to a power of $\rho(g)$. Varying $\tau\in A$, this shows that every element in $\rho(G_{\rm gen})$ is similar to a power of $\rho(g)$. 

Now, assume that we have an element $h\in G_W$ such that $\rho(h)$ is not similar to a power of $\rho(g)$. Then Lemma \ref{lem-generalGW} (6) shows that $1$ and $-1$ are eigenvalues of $\rho(h)$. Applying Lemma \ref{lem-generalGW} (6) again to $h^2$, either $\rho(h)^2$ is similar to a power of $\rho(g)$, or $1$ and $-1$ are eigenvalues of $\rho(h)^2$.
If $1$ and $-1$ are eigenvalues of $\rho(h)$ and of $\rho(h)^2$, which both have determinant 1, then $\rho(h)$ has to be similar to $\diag(1,-1,i,i)$, or to $\diag(1,-1,-i,-i)$. Since $\rho(h)$ defines an automorphism of $B$, by \cite[Thm.13.2.8, Thm.13.3.2]{BirkLang}, $B$ must thus be isogenous to $S\times {E_i}^2$ for some abelian surface $S$. We already know that $B$ is isogenous to $E\times {E_j}^3$, and this contradicts the uniqueness of the Poincaré decomposition of $B$ up to isogeny \cite[Thm.5.3.7]{BirkLang}.
Hence, $\rho(h)^2$ is similar to a power of $\rho(g)$, and $1$ is an eigenvalue of multiplicity at least 2 for it. Hence, $\rho(h)^2=\id_B$, and thus $\rho(h)$, which has determinant one, is similar to $\diag(1,1,-1,-1)$.

We just proved that for any element $h\in G_W$ such that $\rho(h)$ is not similar to a power of $\rho(g)$, it holds that $\rho(h)$ is similar to $\diag(1,1,-1,-1)$. Now, take such an element $h\in G_W$, and note that for $hg\in G_W$, the matrix $\rho(hg)$ is still not similar to a power of $g$ by definition of $h$. Since $\rho(h)$ is similar to $\diag(1,1,-1,-1)$ and $\omega$ has multiplicity three as an eigenvalue of $\rho(g)$, both $\omega$ and $-\omega$ appear as eigenvalues of $\rho(hg)$, thus not similar to $\diag(1,1,-1,-1)$, a contradiction.
This concludes Step 1.}

\step{2}{We show that $G_W=\grp{g}.$

\medskip

By Step 1 and since $\rho$ is faithful, every element of $G_W$ has order $1,2,3,$ or $6$. By Step 1 and since $\rho$ is faithful again, any element of $G_W$ of order $2$ coincides with $g^3$. So $G_W$ has exactly one element of order $2$, and thus $|G_W|=2\cdot 3^{\beta}$ for some $\beta\ge 1$. We now show that $\beta = 1$. Let $S$ be a 3-Sylow subgroup of $G_W$, and $s\in Z(S)$ of order $3$. Arguing by contradiction, let $t\in S\setminus \grp{s}$. By Step 1, every non-trivial element of $\rho(S)$ is similar to $\diag(1,j,j,j)$ or $\diag(1,j^2,j^2,j^2)$. Such are $\rho(t)$ and $\rho(s)$, which also commute since $s$ is central. Looking at their images by $\rho$, either $ts$ or $t^2s$ is trivial, a contradiction. So $\beta = 1$, and $G_W$ has order 6. Hence, $G_W=\grp{g}$. This concludes Step 2.}

\step{3}{We conclude this proof.

\medskip

By \cite[Cor.13.2.4, Prop.13.2.5(c)]{BirkLang}, the number of fixed points of $g$ and $g^3$ on  the abelian fourfold $B$ are respectively 4 and 256. Let $\tau$ be a point of $B$ fixed by $g^3$ but not by $g$.
By Proposition \ref{prop-cyclicincod4}, there is a junior element $h$ generating  the cyclic group $\pstab(W+\tau)$. By Step 2 and by definition of $\tau$, we have $\grp{g^3}\subset \grp{h}\subset G_W = \grp{g}$. Hence, the junior element $h$ has order 6, and $\grp{h}=\grp{g}$,. As both $g$ and $h$ are the only junior elements of order 6 in their generated cyclic groups, $g=h$. But $h$ fixes $\tau$ and $g$ does not, a contradiction.
}}

\demode{Proposition \ref{prop-nojun37}}{It is straightforward from Propositions \ref{prop-globfailjun} and \ref{prop-ruleout6cod4}.}

\section{The isogeny type of \texorpdfstring{$A$}{A}}\label{sec-MakeDisjoint}

This section proves the first part of Theorem \ref{theo-main}, namely the following proposition, inspired by \cite[Proof of Lem.3.4]{OgFibSp}.

\prop[prop-tbmade]{Let $A$ be an abelian variety of dimension $n$, $G$ be a finite group acting freely in codimension $2$ on $A$. Suppose that $A/G$ has a resolution $X$ which is a Calabi-Yau manifold. Then either $A$ is isogenous to ${E_j}^n$ and $G$ is generated by junior elements of order 3 and 6, or $A$ is isogenous to ${E_{u_7}}^n$ and $G$ is generated by junior elements of order 7.
}

\demo{By the $M(G)$-equivariant Poincaré's complete reducibility theorem \cite[Theorem 13.5.2, Proposition 13.5.4, and the paragraph before]{BirkLang}, there are $M(G)$-stable abelian subvarieties $Y_1,\ldots, Y_s$ of $A$ such that:
\begin{enumerate}[label=({\arabic*})]
\item For any $i\in\lint 1,s\rint$, $Y_i$ is isogenous to a power of a $M(G)$-stable $M(G)$-simple abelian subvariety of $A$. In particular, by \cite[Prop.13.5.5]{BirkLang}, there is a simple abelian subvariety $Z_i$ of $Y_i$ such that $Y_i$ is isogenous to a power of $Z_i$.
\item For each $i\ne j$, the set of $M(G)$-equivariant homomorphisms satisfies 
$$\mathrm{Hom}_{M(G)}(Y_i,Y_j)=\{0\}.$$
\item The addition map
$Y_1\times\ldots\times Y_s\to A$
is an $M(G)$-equivariant isogeny.
\end{enumerate}

We define
\begin{align*}
Y_I=\prod_{i\in I} Y_i, &\mbox{ where }I=\{i\in\lint 1,s\rint\mid Z_i\sim E_j\}\\
Y_J=\prod_{j\in J} Y_j, &\mbox{ where }J=\{j\in\lint 1,s\rint\mid Z_j\sim E_{u_7}\}\\
Y_K=\prod_{k\in K} Y_k, &\mbox{ where }K=\lint 1,s\rint\setminus (I\cup J).
\end{align*}

The action of $M(G)$ on $Y_I\times Y_J\times Y_K$ is diagonal by (2), and there is a proper surjective finite morphism $A/M(G)\to Y_I/M(G)\times Y_J/M(G)\times Y_K/M(G)$ induced by the $G$-equivariant addition by (3). Composing with projections, we get proper surjective morphisms $f_I,f_J,f_K$ from $A/M(G)$ to $Y_I/M(G)$, to $Y_J/M(G)$, and to $Y_K/M(G)$.

Let $g\in G$ be a junior element. By Propositions \ref{prop-globfailjun} and \ref{prop-ruleout6cod4}, $g$ has order $3$, or $7$, or $6$ and then five or six non-trivial eigenvalues. By Proposition \ref{prop-singab}, $A$ thus contains an abelian subvariety isogenous to ${E_j}^3$, or to ${E_{u_7}}^3$. Hence, $\dim\, Y_I + \dim\, Y_J\ge 3$, so one of the two quotients $Y_I/M(G)$, $Y_J/M(G)$ has positive dimension. 
Moreover, by Proposition \ref{prop-singab} again, if $g$ has order $3$ or $6$, $M(g)$ acts trivially on $Y_J$ and $Y_K$, and if $g$ has order $7$, it acts trivially on $Y_I$ and $Y_K$. Hence, $M(g)$ acts with determinant 1 on each of the three factors.

But by Lemma \ref{lem-makeitglob} and Corollary \ref{cor-pstabjun}, the group $G$ is generated by its junior elements . Thus, by \cite{Khinich,Watanabe}, the quotients $Y_I/M(G)$, $Y_J/M(G)$ and $Y_K/M(G)$ are thus normal Gorenstein varieties.

We now pullback a certain holomorphic form, namely either the volume form of $Y_I/M(G)$ if $\dim Y_I > 0$, or the volume form of $Y_J/M(G)$ if $\dim Y_I = 0$ (in which case $\dim Y_J > 0$). This yields an $M(G)$-invariant non-zero global holomorphic form of positive degree $y$ on $A$. Note that the sections of $\Omega_{A}^{\cdot}$ are invariant by translations of $A$, so that we in fact have a $G$-invariant non-zero global holomorphic $y$-form on $A$. It pulls back to $X$, which is a Calabi-Yau variety. Hence $y=n$, and either $A\sim {E_j}^{n}$ or $A\sim {E_{u_7}}^{n}$. The order of junior elements generating $G$ is given accordingly by Propositions \ref{prop-singab}, \ref{prop-globfailjun}.}

\section{Junior elements and pointwise stabilizers in codimension 5}\label{sec-codim5}

In this section, we extend the results of Sections \ref{sec-CyclicPointwiseStabilizer} and \ref{sec-RuleOutMostJuniors} to codimension $k=5$. In the first subsection, we exclude the one type of junior element with exactly five non-trivial eigenvalues. In the second subsection, we prove the following result.

\prop[prop-fixor5]{Let $A$ be an abelian variety on which a finite group $G$ acts freely in codimension $2$. Suppose that $A/G$ has a $K$-trivial resolution. Let $W$ be a translated abelian subvariety of codimension $k\le 5$ in $A$ such that $\pstab(W)\ne\{1\}$. Then $ \pstab(W)$ is a cyclic group, generated by one junior element $g$ of order $3$ or $7$.}

\subsection{Ruling out junior elements of order 6 with five non-trivial eigenvalues}\label{subsec-ruleout65}

\prop[prop-jun6cod5]{Let $A$ be an abelian variety, $G$ a group acting freely in codimension $2$ on $A$ such that $A/G$ has $K$-trivial resolution. Then there is no junior element of $G$ whose matrix is similar to 
$\diag(\1_{n-5}, \omega,\omega,\omega,\omega,j).$}

\demo{Suppose by contradiction that there is an element $g\in G$ such that $g(0)=0$ and, in some coordinates, 
$$M(g) = \diag(\1_{n-5}, \omega,\omega,\omega,\omega,j).$$ 
Then there is an abelian subvariety $W$ of codimension 4 in $A$ which is pointwise fixed by $g^3$. By Proposition \ref{prop-cyclicincod4}, $\pstab(W)$ is cyclic, generated by one junior element $h$. As $g^3\in\grp{h}$, $h$ has even order. However, by Propositions \ref{prop-globfailjun} and \ref{prop-ruleout6cod4}, it must have order 3 or 7, contradiction!}

\subsection{The pointwise stabilizer for loci of codimension 5}\label{subsec-EigenspaceDirectSumCod5}

For proving Proposition \ref{prop-fixor5}, it is enough to establish the following result.

\prop[prop-fixor5local]{Let $B$ be an abelian fivefold isogenous to either ${E_j}^5$ or ${E_{u_7}}^5$, and let $p=3$ in the first case, $p=7$ in the second case. Let $F$ be a finite subgroup of $\mathrm{Aut}(B,0)$ generated by junior elements of order $p$, satisfying the following condition

\begin{center}
\begin{minipage}{0.9\textwidth}
\centering
For any non-trivial subgroup $F_0<F$ such that
$\dim\,{\rm Fix}(F_0) \ge 1$,

$F_0$ is generated by one junior element of order $p$.
\end{minipage}
\end{center}

\noindent Then $F$ is cyclic.}

\demode{Proposition \ref{prop-fixor5} admitting Proposition \ref{prop-fixor5local}.}{Let $W$ be a translated abelian subvariety of codimension $k\le 5$ in $A$ such that $\{1\}\ne  \pstab(W)< G$. Propositions \ref{prop-cyclicincod4} and \ref{prop-nojun37} show that if $k\le 4$, then $k=3$ and $\pstab(W)$ is cyclic, generated by one junior element. By Proposition \ref{prop-junclass}, the junior generator has order $3$ or $7$. 

From now on, we assume $k=5$. By Remark \ref{rem-notrans}, we can assume that $G$ contains no non-trivial translation, and up to conjugating the whole group $G$ by a translation, we can assume that $0\in W$. By Proposition \ref{prop-point-wise-stabilizer}, we have a $\pstab(W)$-stable complementary abelian fivefold $B$ to $W$. Consider the group $F=\pstab(W)\subset\mathrm{Aut}(B,0)$. It is generated by junior elements by Proposition \ref{prop-point-wise-stabilizer} (3), and each of them has order $3$ or $7$, by Propositions \ref{prop-globfailjun}, \ref{prop-ruleout6cod4}, and \ref{prop-jun6cod5}. By uniqueness of the Poincaré decomposition of $B$ \cite[Thm.5.3.7]{BirkLang}, the group $\mathrm{Aut}(B,0)$ cannot contain both a junior element of order $3$ and a junior element of order $7$. Let $p$ denote the order of all junior elements in $\mathrm{Aut}(B,0)$.

We claim that the following condition is satisfied
\begin{center}
\begin{minipage}{0.9\textwidth}
\centering
For any non-trivial subgroup $F_0<F$ such that
$\dim\,{\rm Fix}(F_0) \ge 1$,

$F_0$ is generated by one junior element of order $p$.
\end{minipage}
\end{center}
Indeed, let $F_0$ be such a subgroup, and let $B_0$ be an abelian subvariety of $B$ of maximal dimension in ${\rm Fix}(F_0)$. By definition, we have $F_0 < \pstab(B_0)$, and $B_0$ has dimension at least $1$, hence codimension at most 4 in $B$. By Propositions \ref{prop-cyclicincod4}, \ref{prop-globfailjun}, \ref{prop-ruleout6cod4}, the group $\pstab(B_0)$ is cyclic, generated by one junior element of order $3$ or $7$. Since $F_0$ is non-trivial and $3$ and $7$ are prime numbers, we obtain that $F_0=\pstab(B_0)$ is cyclic, generated by one junior element of order $3$ or $7$. In fact, that element has order $p$ as defined at the end of the previous paragraph. This shows our claim.
 
Let us prove that $F$ is cyclic. By contradiction, we assume that there are two junior elements $g,h\in F$ such that $\grp{g}\neq\grp{h}$. From our claim on subgroups of $F$, we have that $\dim\,{\rm Fix}(\grp{g,h})\le 0$. Recall that the origin $0\in B$ is fixed by any element of $F$, including $g$ and $h$. Hence, $M(g)$ and $M(h)$ have no common eigenvector of eigenvalue $1$.
By the $\grp{g,h}$-equivariant Poincaré's complete reducibility theorem \cite[Thm.13.5.2, Prop.13.5.4, and the paragraph before]{BirkLang}, the abelian fivefold $B$ is thus isogenous to ${E_j}^5$ (if $p=3$) or ${E_{u_7}}^5$ (if $p=7$). By Proposition \ref{prop-fixor5local}, we obtain that $F$ is cyclic, a contradiction.}

To establish Proposition \ref{prop-fixor5local}, we start by stating and proving a lemma.

\lem[lem-fdgroupdim5]{Let $B$ be an abelian fivefold isogenous to either ${E_j}^5$ or ${E_{u_7}}^5$, and let $p=3$ in the first case, $p=7$ in the second case. Let $F$ be a finite subgroup of $\mathrm{Aut}(B,0)$ generated by junior elements of order $p$, satisfying the following condition

\begin{center}
\begin{minipage}{0.9\textwidth}
\centering
For any non-trivial subgroup $F_0<F$ such that
$\dim\,{\rm Fix}(F_0) \ge 1$,

$F_0$ is generated by one junior element of order $p$.
\end{minipage}
\end{center}
Then $F$ is a $p$-group.}

\demo{By Cauchy's theorem, it suffices to prove that any element $g$ of $F$ of prime order is in fact of order $p$. Let $g$ be an element of $F$ of prime order $q$.

If $1$ is an eigenvalue of $g$, then $\dim\,{\rm Fix}(\grp{g}) \ge 1$, so by our subgroup condition on $F$, the group $\grp{g}$ is generated by one junior element of order $p$, thus $q=p$, as wished. 

From now on, we assume that $1$ is not an eigenvalue of $g$. As $g$ has prime order, by Lemma \ref{lem-blrat}, the characteristic polynomial $\chi_{g\oplus \overline{g}}$ is a power of the cyclotomic polynomial $\Phi_q$. Hence, $\deg(\Phi_q)=q-1$ divides $10$, so $q\in\{2,3,11\}$. Let us analyze these cases.
\begin{itemize}
\item If $q=2$, then the only possible eigenvalue for $g$ is $-1$. But $g\in F$ has determinant $1$ and $\dim\, B = 5$, a contradiction.
\item If $q=11$, we have $\chi_{g}\overline{\chi_{g}}=\Phi_{11}$, which by \cite[Prop.2.4]{Wash} is irreducible both over $\Q[j]$ and over $\Q[\zeta_7]$. This contradicts the fact that $\chi_g$ is a polynomial in $\Q[j][X]$ if $p=3$, and in $\Q[\zeta_7][X]$ if $p=7$.
\item If $q=3$, then $\chi_{g}\overline{\chi_{g}} = {\Phi_3}^5$, yet by \cite[Prop.2.4]{Wash}, $\Phi_3$ is irreducible over $\Q[\zeta_7]$. Hence, we have $p\neq 7$. So $p=3$, and thus $q=p$ in this case.
\end{itemize}
We have shown that $q=p$ in every case, and this concludes the proof.}

\demode{Proposition \ref{prop-fixor5local}}{By Lemma \ref{lem-fdgroupdim5}, the group $F$ is a $p$-group. Hence, there is a central element $g\in Z(F)$ of order $p$. Assume by contradiction that $F$ is not cyclic. Hence, we can find a junior element $h\in F\setminus \grp{g}$ of order $p$. Since $h\not\in\grp{g}$, the group $\grp{g,h}$ is not cyclic.
By our assumption on subgroups of $F$, and since $g$ and $h$ both fix $0\in B$, we deduce that $E_g(1)\cap E_h(1)=\{0\}$. As $g$ is central, $g$ and $h$ are codiagonalizable.

If $p=7$, this yields that $gh$ has at least four eigenvalues (counted with multiplicity) of order 7. By Lemma \ref{lem-blrat}, we have a splitting of the characteristic polynomial 
$$\chi_{gh\oplus\overline{gh}}={\Phi_1}^{\alpha}\cdot{\Phi_7}^{\beta}$$
Counting eigenvalues of order $7$ yields $3\beta\ge 4$, whereas computing degrees shows that $10=\alpha+6\beta$, a contradiction.

\medskip

If $p=3$, the condition for subgroups of $F$ with positive common fixed locus, together with having determinant one, forces every element of order $3$ in $F$ to be similar to one of the matrices
$$\diag(1,1,j,j,j),\;\diag(1,1,j^2,j^2,j^2),\;\diag(j,j,j,j,j^2),\;\diag(j,j^2,j^2,j^2,j^2).$$
Typically, an element similar to $\diag(1,j,j,j^2,j^2)$ cannot occur, as it would have a positive-dimensional fixed locus, yet not be the power of a junior element. Let $\chi$ be the character of the representation $\grp{g,h}\subset\mathrm{Aut}(B,0)$ of rank $4$. Let $a$ be the number of junior elements in $\grp{g,h}$, and $b$ the number of elements in $\grp{g,h}$ similar to $\diag(j,j,j,j,j^2)$. 
Since $E_g(1)\cap E_h(1)=\{0\}$, we have
$$0=\langle \chi,\1\rangle = 5 + a(2+3j)+a(2+3j^2)
+ b(4j+j^2)+ b(4j^2+j).$$
Since $g$ and $h$ commute and have order $3$, the group $\grp{g,h}$ is isomorphic to $\Z_3\times\Z_3$, thus $1+2a+2b=9$. This yields $0=-15+6a$, which has no integral solution, a contradiction.}

\section{Junior elements and pointwise stabilizers in codimension 6}\label{sec-codim6}

The goal of this section is to extend the results of Sections \ref{sec-CyclicPointwiseStabilizer}, \ref{sec-RuleOutMostJuniors}, \ref{sec-codim5} to codimension $k=6$. For the first time in our study of pointwise stabilizers, and for the second time in this paper after Section \ref{sec-MakeDisjoint}, we need to assume the existence of a Calabi-Yau resolution, and not just a $K$-trivial resolution of the singular quotient $A/G$.
Indeed, in dimension 6, products of the two examples of \cite{III0} can arise. Locally, they create non-cyclic pointwise stabilizers. Globally, they yield $K$-trivial resolutions of singular quotients $A/G$ that are not Calabi--Yau manifolds.

We start by proving the following partial classification of pointwise stabilizers in codimension 6 in Subsection \ref{subsec-pstabcod6}.

\prop[prop-fixor6int]{Let $A$ be an abelian variety on which a finite group $G$ acts freely in codimension $2$. Suppose that $A/G$ has a resolution $X$ which is a Calabi-Yau manifold. Let $W$ be a translated abelian subvariety of codimension $k\le 6$ in $A$ such that $\pstab(W)\ne\{1\}$ contains no junior element of type $\diag(\1_{n-6},\omega,\omega,\omega,\omega,\omega,\omega)$.
Then $\pstab(W)$ is one of the following.
\begin{itemize}
\item A cyclic group generated by one junior element of order $3$ or $7$.
\item An abelian group generated by two junior elements $g$ and $h$ of order both $3$ or both $7$, satisfying
$E_g(1)\cap E_h(1)=H^0(W,T_W).$
\item $\SL_2(\F_3)$, and the representation $M:\pstab(W)\hookrightarrow\mathrm{Aut}(A,0)$ decomposes as $1^{\oplus n-6}\oplus\sigma^{\oplus 3}$, where $\sigma$ is the unique irreducible $2$-dimensional faithful representation of $\SL_2(\F_3)$ over the splitting field $\Q[j]$.
\end{itemize}}

We then use this result to rule out the existence of junior elements with six non-trivial eigenvalues in Subsection \ref{subsec-ruleout66} by a mix of local and global arguments, and finally refine Proposition \ref{prop-fixor6int} in Subsection \ref{subsec-JuniorsCommute} to the following result.

\prop[prop-fixor6]{Let $A$ be an abelian variety on which a finite group $G$ acts freely in codimension $2$. Suppose that $A/G$ has a resolution $X$ which is a Calabi-Yau manifold. Let $W$ be a translated abelian subvariety of codimension $k\le 6$ in $A$ such that $\pstab(W)\ne\{1\}$.
Then $\pstab(W)$ is one of the following.
\begin{itemize}
\item A cyclic group generated by one junior element of order $3$ or $7$.
\item An abelian group generated by two junior elements $g$ and $h$ of order both $3$ or both $7$, satisfying
$E_g(1)\cap E_h(1)=H^0(W,T_W).$
\end{itemize}}

\subsection{The pointwise stabilizers for loci of codimension 6}\label{subsec-pstabcod6}

For proving Proposition \ref{prop-fixor6int}, it is enough to establish the following result.

\prop[prop-fixor6intlocal]{Let $B$ be an abelian sixfold isogenous to either ${E_j}^6$ or ${E_{u_7}}^6$, and let $p=3$ in the first case, $p=7$ in the second case. Let $F$ be a finite subgroup of $\mathrm{Aut}(B,0)$ generated by junior elements of order $p$, satisfying the following condition

\begin{center}
\begin{minipage}{0.9\textwidth}
\centering
For any non-trivial subgroup $F_0<F$ such that
$\dim\,{\rm Fix}(F_0) \ge 1$,

$F_0$ is generated by one junior element of order $p$.
\end{minipage}
\end{center}

Suppose also that $\omega\id_B\not\in F$. Then $F$ is one of the following.
\begin{itemize}
\item A cyclic group generated by one junior element of order $p$.
\item An abelian group generated by two junior elements $g$ and $h$ of order $p$ satisfying
$E_1(g)\cap E_1(h)=\{0\}.$
\item $\SL_2(\F_3)$, and the representation $M:F\hookrightarrow\mathrm{Aut}(B,0)$ decomposes as $\sigma^{\oplus 3}$, where $\sigma$ is the unique irreducible $2$-dimensional faithful representation of $\SL_2(\F_3)$ over the splitting field $\Q[j]$. In this case, $p=3$.
\end{itemize}}

\demode{Proposition \ref{prop-fixor6int} admitting Proposition \ref{prop-fixor6intlocal}.}{Let $W$ be a translated abelian subvariety of codimension $k\le 6$ in $A$ such that $\pstab(W)\ne\{1\}$ contains no junior element of type $\diag(\1_{n-6},\omega,\omega,\omega,\omega,\omega,\omega)$. Proposition \ref{prop-fixor5} settles the cases when $k\le 5$, so we can assume $k=6$. Up to conjugating the whole group $G$ by a translation, we can assume that $0\in W$, and apply Proposition \ref{prop-point-wise-stabilizer} to obtain a $\pstab(W)$-stable complementary abelian sixfold $B$ to $W$. By Proposition \ref{prop-tbmade} and as an abelian subvariety of $A$, $B$ is isogenous to either ${E_j}^6$ or ${E_{u_7}}^6$.

Let $F=\pstab(W)\subset\mathrm{Aut}(B,0)$. It is generated by junior elements by Proposition \ref{prop-point-wise-stabilizer} (3), which have order $3$ or $7$ by Propositions \ref{prop-nojun37}, \ref{prop-ruleout6cod4}, \ref{prop-jun6cod5} and since, by assumption, $\omega\id_B\not\in F$. Let $p$ be the order of the junior elements generating $F$.

Let $F_0$ be a non-trivial subgroup of $F$ such that $\dim\,{\rm Fix}(F_0)\ge 1$. Let $B_0$ be an abelian subvariety of positive dimension in $B$ contained in ${\rm Fix}(F_0)$. Note that $F_0<\pstab(B_0)$, and that $B_0$ has codimension at most $5$ in $B$. By Proposition \ref{prop-fixor5}, $\pstab(B_0)$ is cyclic, generated by one junior element of order $p$, so since $p$ is prime, $F_0=\pstab(B_0)$, and it is again cyclic, generated by one junior element of order $p$.

So Proposition \ref{prop-fixor6intlocal} applies, and proves Proposition \ref{prop-fixor6int}.}

To establish Proposition \ref{prop-fixor6intlocal}, we need numerous lemmas.

\lem[lem-fdgroup]{Let $B$ be an abelian sixfold isogenous to either ${E_j}^6$ or ${E_{u_7}}^6$, and let $p=3$ in the first case, $p=7$ in the second case. Let $F$ be a finite subgroup of $\mathrm{Aut}(B,0)$ generated by junior elements of order $p$, satisfying the following condition

\begin{center}
\begin{minipage}{0.9\textwidth}
\centering
For any non-trivial subgroup $F_0<F$ such that
$\dim\,{\rm Fix}(F_0) \ge 1$,

$F_0$ is generated by one junior element of order $p$.
\end{minipage}
\end{center}
Let $g\in F$ be an element of prime order $q$. Then, we have $q\in\{2,3,7\}$.}

\demo{If 1 is an eigenvalue of $g$, then the subgroup condition applies to $\grp{g}$ and yields $q=p$, as wished. 
From now on, we assume that $1$ is not an eigenvalue of $g$. By Lemma \ref{lem-blrat}, the characteristic polynomial $\chi_{g\oplus \overline{g}}$ is thus a power of $\Phi_q$, so $q-1$ divides $12$, so $q\in\{2,3,5,7,13\}$.
\begin{itemize}
\item If $q=13$, then $\Phi_{13}=\chi_g\overline{\chi_{g}}$. But by \cite[Prop.2.4]{Wash}, $\Phi_{13}$ is irreducible over $\Q[j]$ and $\Q[\zeta_7]$, a contradiction.
\item If $q=5$, then $\Phi_{5}^3=\chi_g\overline{\chi_{g}}$. But by \cite[Prop.2.4]{Wash}, the cyclotomic polynomial $\Phi_{5}$ is irreducible over $\Q[j]$ and $\Q[\zeta_7]\supset\Q[u_7]$, a contradiction.
\end{itemize}
}

Let us now describe the $2$-, $3$-, and $7$-Sylow subgroups of $F$.

\lem[lem-the2sylcod6]{Let $B$ be an abelian sixfold isogenous to either ${E_j}^6$ or ${E_{u_7}}^6$, and let $p=3$ in the first case, $p=7$ in the second case. Let $F$ be a finite subgroup of $\mathrm{Aut}(B,0)$ generated by junior elements of order $p$, satisfying the following condition
\begin{center}
\begin{minipage}{0.9\textwidth}
\centering
For any non-trivial subgroup $F_0<F$ such that
$\dim\,{\rm Fix}(F_0) \ge 1$,

$F_0$ is generated by one junior element of order $p$.
\end{minipage}
\end{center}
If $2$ divides $|F|$, then a $2$-Sylow subgroup $S$ of $F$ is isomorphic to $Q_8$.}

\demo{Since $-\id_B$ is the unique element of order 2 that can belong to $F$, by \cite[5.3.6]{Robinson}, $S$ is cyclic or a generalized quaternion group. Let us show that $S$ has no element of order $8$. By contradiction, let $s\in S$ be of order $8$. Since $s^4=-\id_B$, all eigenvalues of $s$ have order $8$, so the characteristic polynomial $\chi_{s\oplus\overline{s}}$ is a power of $\Phi_8$. Comparing degrees yields ${\Phi_8}^3=\chi_{s}\overline{\chi_s}$. But by \cite[Prop.2.4]{Wash}, $\Phi_8$ is irreducible over $\Q[j]$ and $\Q[\zeta_7]$, a contradiction! So $S$ is isomorphic to $\Z_2,\Z_4$, or $Q_8$.

If $S$ is cyclic, then by \cite[10.1.9]{Robinson}, there is a normal subgroup $N$ of $F$ such that $F=N\rtimes S$. But all junior elements of $F$ have odd order, so they belong to $N$ and cannot generate $F$, contradiction!  So $S$ is isomorphic to $Q_8$. 
}

\lem[lem-2749]{Let $B$ be an abelian sixfold. Let $g\in\mathrm{Aut}(B,0)$ be an element of finite order. Then $g$ cannot have order $27,49,$ or $63$.}

\demo{It is an immediate consequence of Lemma \ref{lem-blrat}.}

\lem[lem-theothersylcyclic]{Let $B$ be an abelian sixfold isogenous to either ${E_j}^6$ or ${E_{u_7}}^6$, and let $p=3$ in the first case, $p=7$ in the second case. Let $F$ be a finite subgroup of $\mathrm{Aut}(B,0)$ generated by junior elements of order $p$, satisfying the following condition
\begin{center}
\begin{minipage}{0.9\textwidth}
\centering
For any non-trivial subgroup $F_0<F$ such that
$\dim\,{\rm Fix}(F_0) \ge 1$,

$F_0$ is generated by one junior element of order $p$.
\end{minipage}
\end{center}
Let $q=7$ if $p=3$, $q=3$ if $p=7$.
If $q$ divides $|F|$, a $q$-Sylow subgroup $S$ of $F$ is cyclic and has order $3,7$, or $9$.}

\demo{As $S$ is a $q$-group, there is an element $g\in Z(S)$ of order $q$. Let $h\in S\setminus\grp{g}$ be another element of order $q$.
Because $q\not\in\{2,p\}$, the elements $g,h\in F$ cannot be powers of junior elements, so by the subgroup condition in $F$, neither $g$ nor $h$ admits 1 as an eigenvalue. By Lemma \ref{lem-blrat}, since they have determinant one, and since $B$ is isogenous to ${E_{u_7}}^6$ if $q=3$ and ${E_j}^6$ is $q=7$, the matrices $g$ and $h$ are similar to
\begin{align*}
\diag(j,j,j,j^2,j^2,j^2) &\mbox{ if $q=3$}\\
\diag(\zeta_7,{\zeta_7}^2,\zeta_7^3,{\zeta_7}^4,\zeta_7^5,\zeta_7^6) &\mbox{ if $q=7$}
\end{align*}
One can then find a non-trivial element of $\grp{g,h}$ with 1 as an eigenvalue. But as $g$ and $h$ commute, it has order $q\not\in\{2,p\}$, contradiction.
So $\grp{g}$ is the unique subgroup of order $q$ in $S$. By \cite[5.3.6]{Robinson}, $S$ is thus cyclic, and its order is restricted by Lemma \ref{lem-2749}.
}

\lem[lem-the3and7syl]{Let $B$ be an abelian sixfold isogenous to either ${E_j}^6$ or ${E_{u_7}}^6$, and let $p=3$ in the first case, $p=7$ in the second case. Let $F$ be a finite subgroup of $\mathrm{Aut}(B,0)$ generated by junior elements of order $p$, satisfying the following condition
\begin{center}
\begin{minipage}{0.9\textwidth}
\centering
For any non-trivial subgroup $F_0<F$ such that
$\dim\,{\rm Fix}(F_0) \ge 1$,

$F_0$ is generated by one junior element of order $p$.
\end{minipage}
\end{center}
Then a $p$-Sylow subgroup $S$ of $F$ is either cyclic, or the direct product of two cyclic groups. In fact, $S$ is isomorphic to one of the following groups
\begin{align*}
\Z_3,\;\Z_9,\;\Z_3\times\Z_3,\mbox{ or }\Z_3\times\Z_9 &\mbox{ if $p=3$}\\
\Z_7,\mbox{ or }\Z_7\times\Z_7&\mbox{ if $p=7$}
\end{align*}
}

\demo{Let $g\in Z(S)$ be an element of order $p$. If $\grp{g}$ is the only subgroup of order $p$ in $S$, then by \cite[5.3.6]{Robinson}, $S$ is cyclic. Control on its order follows from Lemma \ref{lem-2749} and concludes this case.

From now on, we assume that $\grp{g}$ is not the only subgroup of order $p$ in $S$. Let $h$ be an element of $S$ whose class $[h]$ has order $p$ and is central in the $p$-group $S/\grp{g}$. If we can show that $\grp{[h]}$ is the only subgroup of order $p$ in the quotient group $S/\grp{g}$, note that it concludes this proof. Indeed, if it is the case, then by \cite[5.3.6]{Robinson}, $S/\grp{g}$ is cyclic, and thus $S/Z(S)$ is cyclic too. So $S$ is abelian, and $S\simeq\grp{g}\times C$ for a cyclic group $C$ containing $\grp{h}$. Control on the order of $C$ follows from Lemma \ref{lem-2749}, and then concludes the proof. We now fix an element $k\in S$ such that $[k]$ has order $p$ in $S/\grp{g}$, and are left to show that $[h]$ and $[k]$ span the same subgroup of $S/\grp{g}$.

If $p=7$, then $g$ has an eigenvalue $\zeta$ of order 7 with corresponding eigenspace $E_{g}(\zeta) $ of dimension $1$. By Lemma \ref{lem-2749}, $h$ and $k$ have order $7$ in $S$. As $g$ commutes with $h$ and $k$, we can thus choose $h'\in[h],k'\in[k]$ which both have 1 as an eigenvalue on $E_{g}(\zeta)$. Hence, the group $\grp{h',k'}$ does not act freely in codimension 5 on $B$, so it is cyclic generated by one junior element, and $\grp{h'}=\grp{k'}$ as wished.

If $p=3$, let us show that $j\id_B\in S$. By contradiction, suppose that elements of order $3$ in $S$ are all similar to one of the following matrices
$$\diag (1,1,1,j,j,j),\; \diag (1,1,1,j^2,j^2,j^2),\; \diag (j,j,j,j^2,j^2,j^2).$$
Take $s\in S\setminus\grp{g}$. As $g$ and $s$ commute, a simple computation shows that one of the products $gs$, $g^2s$, $gs^2$, $g^2s^2$ will not fall under these three similarity classes, contradiction.

Hence, we can refine the choice of our central element $g$ from here on by setting $g=j\id_B$. A fortunate consequence of that choice, of Lemma \ref{lem-blrat}, and of the fact that matrices in $S$ all have determinant $1$ is that $g$ has no cubic root in $S$, i.e., every element of order $9$ in $S$ has a class of order 9 in $S/\grp{g}$. Hence, $h$ and $k$ above have order 3. Moreover, recall that $hkh^{-1}k^{-1}\in\grp{g}=\grp{j\id_B}$. If $k$ is conjugated to $jk$ or $j^2k$, then $1,j,$ and $j^2$ each are eigenvalues of $k$, which contradicts the subgroup condition on $F$ again.
Hence, $hkh^{-1}=k$, i.e., $h$ and $k$ commute. They commute with $g$ as well, and thus we can find some non-trivial elements $h'\in [h]$ and $k'\in [k]$ with a common eigenvector of eigenvalue 1. So $\grp{h'}=\grp{k'}$, which concludes this proof.}

\demode{Proposition \ref{prop-fixor6intlocal}}{
We now run (see Appendix) a \texttt{GAP} search through all groups with such $2,3,$ and $7$-Sylow subgroups, which have at most an element of order 2, and no element of order $63$. Among the ninety-four of them, only $\Z_7$ and $\Z_7\times\Z_7$ can be generated by their elements of order $7$, whereas $\Z_3,\Z_3\times\Z_3,\SL_2(\F_3),Q_8\rtimes(\Z_7\rtimes\Z_3),$ and $\Z_3\times (Q_8\rtimes(\Z_7\rtimes\Z_3))$ can be generated by their elements of order $3$.
However, it is easy to check that $Q_8\rtimes(\Z_7\rtimes\Z_3),$ and $\Z_3\times (Q_8\rtimes(\Z_7\rtimes\Z_3))$ have elements of order 28, which by Lemma \ref{lem-blrat} and \cite[Prop.2.4]{Wash} cannot occur in $\mathrm{Aut}_{\Q}({E_j}^6,0)$.

The representation theoretic description is easily obtained from \texttt{GAP} for $\SL_2(\F_3)$, and follows from the condition about freeness in codimension 5 for $\Z_3\times\Z_3$ and $\Z_7\times
\Z_7$.}

\subsection{Ruling out junior elements of order 6 with six non-trivial eigenvalues}\label{subsec-ruleout66}

\prop[prop-jun6cod6]{Let $A$ be an abelian variety, $G$ a group acting freely in codimension $2$ on $A$ such that $A/G$ has a $K$-trivial resolution. Then there is no junior element of $G$ with matrix similar to 
$\diag(\1_{n-6}, \omega,\omega,\omega,\omega,\omega,\omega).$}

In order to prove this, we first reduce to a $6$-dimensional situation, where a lot of local information is given by Proposition \ref{prop-fixor6intlocal}.

\lem[lem-fixtauinfix]{Let $A$ be an abelian variety, $G$ a group acting freely in codimension $2$ on $A$ without translations such that $A/G$ has a $K$-trivial resolution. Suppose that there is an element $g\in G$ such that $g(0)=0$, and with matrix similar to 
$\diag(\1_{n-6}, \omega,\omega,\omega,\omega,\omega,\omega).$
Then there are complementary $\grp{g}$-stable abelian subvarieties $B$ and $W$ in $A$ such that $g|_B = \omega\id_B$ and $g|_W=\id_W$.
Moreover, for any $\tau\in B$, it holds $\pstab(W+\tau)\subset \pstab(W)$, and if $\tau$ is a non-zero 2-torsion point of $B$, we have $\pstab(W+\tau)\simeq\SL_2(\F_3)$.}

\demo{
As in Lemma \ref{lem-generalGW}, we introduce $G_W$, the subgroup of $G$ generated by
$$G_{\mathrm{gen}} = {G_{\mathrm{gen}}}^{-1} = \left\{h\in G\mid \exists\,\tau\in A\mbox{ such that } h\in \pstab(W+\tau)\right\}.$$
By Lemma \ref{lem-generalGW} (1) (4), there is a $G_W$-stable complementary $B$ to $W$. The fact that $\omega\id_B\in\mathrm{Aut}(B,0)$ implies that $B$ is isogenous to ${E_j}^6$, by Proposition \ref{prop-singab}.

We prove that $G_W=\pstab(W)$. From that, it clearly follows that $\pstab(W+\tau)\subset\pstab(W)$ for any $\tau\in A$. Let $h\in G_W$. By Lemma \ref{lem-generalGW} (6) for $k=6$, $d=6$, $\alpha=3$, we have that $g^3$ and $h$ commute, i.e., they are codiagonalizable. As the matrices $M(g^3)$ and $M(g)$ have exactly the same eigenspaces, the matrices $M(g)$ and $M(h)$ commute too. Since $G$ contains no translation, $M$ is faithful, and thus $g$ and $h$ commute.
By Lemma \ref{lem-generalGW} (2), we have $h(0)=T(h)\in B$, hence $g(h(0))=\omega h(0)$. But since $g$ and $h$ commute, we also have $g(h(0))=h(0)$. As by \cite[Cor.13.2.4]{BirkLang}, $\omega\id_B$ has exactly one fixed point on $B$, namely $0$, we obtain $h(0)=0$, whence $h\in \pstab(W)$.

Consider now a non-zero $2$-torsion point $\tau\in B$. Then $g^3$ fixes $\tau$, but $g$ does not, the group $\pstab(W+\tau)$ contains $g^3$ but not $g$. If $\pstab(W+\tau)$ contains a junior element $h$ of type $\diag(\1_{n-6}, \omega,\omega,\omega,\omega,\omega,\omega)$, then $g$ and $h$ coincide, a contradiction. Hence, Proposition \ref{prop-fixor6intlocal} can be applied to $\pstab(W+\tau)$. Since its element $g^3$ has order $2$, it can only be isomorphic to $\SL_2(\F_3)$, within the list given by Proposition \ref{prop-fixor6intlocal}.}

\rem{Note that we showed the following intermediate result: If $G$ contains a junior element $g$ of type $\diag(\1_{n-6},\omega,\omega,\omega,\omega,\omega,\omega)$ such that $g(0)=0$, and $W$ is the maximal abelian subvariety of $A$ fixed by $g$, then the group $G_W$ defined in Lemma \ref{lem-generalGW} equals $\pstab(W)$.}

This description of the pointwise stabilizers of the translations of $W$ by $2$-torsion points yields the following description of the much larger group $\pstab(W)$.

\lem[lem-2sylfix]{Let $A$ be an abelian variety, $G$ a group acting freely in codimension $2$ on $A$ without translations such that $A/G$ has a $K$-trivial resolution. Suppose that there is an element $g\in G$ such that $g(0)=0$, and with matrix similar to 
$\diag(\1_{n-6}, \omega,\omega,\omega,\omega,\omega,\omega).$
Let $B,W$ be as in Lemma \ref{lem-fixtauinfix}.
Then, every element of prime order in $\pstab(W)$ has order $2$ or $3$.
Moreover, a $2$-Sylow subgroup $S_2$ of $\pstab(W)$ is isomorphic to $Q_8$, and a $3$-Sylow subgroup $S_3$ contains an even number of junior elements.
The group $\pstab(W)$ contains exactly $260$ junior elements of order $3$.}

\demo{The group $\pstab(W)$ contains a unique element $g^3$ of order $2$, so by \cite[5.3.6]{Robinson}, its $2$-Sylow subgroup $S_2$ is cyclic or a generalized quaternion group. Moreover, $\pstab(W)$ acts on a complementary abelian variety to $W$, which is isomorphic to ${E_j}^6$ by Proposition \ref{prop-singab}, and the elements of $\pstab(W)$ with $1$ as an eigenvalue pointwise fix an abelian subvariety of $A$ of codimension at least $5$, and therefore by Proposition \ref{prop-fixor5}, they are powers of junior elements. Hence, $\pstab(W)\subset\SL_6(\Q[j])$ has no element of order $8$, so $S_2$ is isomorphic to $\Z_2,\Z_4,$ or $Q_8$. But by Lemma \ref{lem-fixtauinfix}, a copy of $Q_8\subset\SL_2(\F_3)$ embeds in $\pstab(W)$, and therefore $S_2\simeq Q_8$.

We now count the number of junior elements of order $3$ in $\pstab(W)$ easily: Each of them fixes exactly $2^6-1$ non-zero $2$-torsion points of $B$, and every non-zero $2$-torsion point of $B$ is fixed by exactly $4$ junior elements by Lemma \ref{lem-fdgroup}. Since $B$ has $2^{12}-1$ non-zero $2$-torsion points, the number of junior elements in $\pstab(W)$ is $\frac{(2^{12}-1)\cdot 4}{2^6-1}=260$.

Now, consider an element $h\in\pstab(W)$ of prime order $q$. Suppose by contradiction that $q\ne 2,3$. By Lemma \ref{lem-fdgroup}, we have $q=7$. Since $\SL_6(\Q[j])$ has no junior element of order $7$ and by Proposition \ref{prop-fixor5}, $h$ does not admit $1$ as an eigenvalue. Hence, all six eigenvalues of $h$ have order $7$. Note that $h$ acts by conjugation on the set of junior elements of $\pstab(W)$, whose cardinal, which we just computed, is $260$, which is not divisible by $7$. Hence, $h$ commutes with a junior element $k\in\pstab(W)$, so $hk\in\pstab(W)$ has order $21$, and admits three eigenvalues of order $7$, and three eigenvalues of order $21$. By Lemma \ref{lem-blrat}, $\Phi_7\Phi_{21}$ thus divides the characteristic polynomial of $hk\oplus\overline{hk}$, but they have respective degrees $\phi(7)+\phi(21)=18$ and $12$, a contradiction.

The group $\pstab(W)$ contains the element $g^2$, which has order $3$, yet is not junior. Note that $g^2$ commutes with every element of $\pstab{W}$, and thus belongs to any $3$-Sylow subgroup of it. Let $S_3$ be a $3$-Sylow subgroup of $\pstab(W)$.
Now, the group homomorphism $h\in S_3\mapsto g^2h^2\in S_3$ sends any junior element of order 3 to a junior element of order 3, and it is a fixed-point-free involution. Hence, $S_3$ contains an even number of junior elements of order 3 (and any junior element contained in a $S_3$ has order a power of $3$, hence exactly $3$ anyways).}

This result has the following consequence.

\cor[cor-65ins]{Let $A$ be an abelian variety, $G$ a group acting freely in codimension $2$ on $A$ without translations such that $A/G$ has a crepant resolution $X$. Suppose that there is an element $g\in G$ such that $g(0)=0$, and with matrix similar to 
$\diag(\1_{n-6}, \omega,\omega,\omega,\omega,\omega,\omega).$
Let $B,W$ be as in Lemma \ref{lem-fixtauinfix}.
Then the group $ \pstab(W)$ has exactly four 3-Sylow subgroups $S$, $T$, $U$ and $V$. Each $3$-Sylow subgroup contains exactly $65$ junior elements.
}

\demo{By Lemma \ref{lem-2sylfix}, there is a positive integer $\beta$ such that
$$| \pstab(W)|=8\cdot 3^{\beta}.$$
The number $n_3$ of 3-Sylow subgroups in $ \pstab(W)$ is thus either $1$, or $4$.

Let us fix $\tau_0$ as a non-zero 2-torsion point in $B$. By Lemma \ref{lem-fixtauinfix}, there are exactly four junior elements $s,t,u,v$ of order 3 in $\SL_2(\F_3)\simeq\pstab(W+\tau_0)\subset\pstab(W)$. We can check in the multiplication table of $\SL_2(\F_3)$ that the product of any two distinct elements of $\{s,t,u,v\}$ has order 6. Hence, each 3-Sylow subgroup of $\pstab(W)$ contains at most one element of $\{s,t,u,v\}$. So $n_3\ge 4$, hence $n_3=4$.
Denote by $S$, $T$, $U$, and $V$ the four 3-Sylow subgroups of $\pstab(W)$. 

We claim that the junior elements contained in $S$, $T$, $U$, and $V$ form a partition of the set of junior elements of $\pstab(W)$. If this claim is true, then by the second Sylow theorem, these four partitioning pieces are in bijection, so each $3$-Sylow subgroup has $\frac{260}{4}=65$ junior elements.

Let us prove our claim. Consider a junior element $h\in\pstab(W)$. Necessarily, it is of order 3. By Lemma \ref{lem-generalGW} (9), $h$ fixes a non-zero $2$-torsion point $\tau\in B$, hence $h$ belongs to $\pstab(W+\tau)$, and in particular to one of the $3$-Sylow subgroups of $\pstab(W+\tau)$. Since the four $3$-Sylow subgroups of $\SL_2(\F_3)\simeq\pstab(W+\tau)$ are disjoint and coincide with $S\cap\pstab(W+\tau)$, $T\cap\pstab(W+\tau)$, $U\cap\pstab(W+\tau)$, $V\cap\pstab(W+\tau)$, the junior element $h$ is contained in exactly one of them, and thus contained in exactly one of $S$, $T$, $U$, $V$. This proves our claim, and concludes.}

\demode{Proposition \ref{prop-jun6cod6}}{By contradiction, suppose that $G$ contains a junior element $g$ of type $\diag(\1_{n-6},\omega,\omega,\omega,\omega,\omega,\omega)$. By Remark \ref{rem-notrans}, we can assume that $G$ contains no translation other than $\id_A$, and up to conjugating the whole group by a translation, we can assume that $g(0)=0$. Now, Lemma \ref{lem-2sylfix} and Corollary \ref{cor-65ins} apply, but since $65$ is odd, they contradict one another.}

\subsection{Ruling out the pointwise stabilizer \texorpdfstring{$\SL_2(\F_3)$}{SL2(F3)}}\label{subsec-JuniorsCommute}

In this subsection, we prove Proposition \ref{prop-fixor6}. Since Proposition \ref{prop-fixor6int} has been established in the previous section, it is now enough to show the following:

\lem[lem-sl23]{Let $A$ be an abelian variety on which a finite group $G$ acts freely in codimension $2$ without translations. Suppose that $A/G$ has a resolution $X$ that is a Calabi-Yau manifold.
Then, there is no translated abelian subvariety $W$ of codimension $6$ in $A$ such that $\pstab(W)\simeq\SL_2(\F_3) < G$, with representation $M=1^{\oplus n-6}\oplus \sigma^{\oplus 3}$ as in Proposition \ref{prop-fixor6int}.}

This result resembles \cite[Sec.6.1]{AndreaWis}, although working under a different set of assumptions and in dimension 6.

\demo{We prove it by contradiction, using global arguments. Consider such an translated abelian subvariety $W$. Up to conjugating the whole groupe $G$ by a translation, we can assume that $0\in W$. Following Lemma \ref{lem-generalGW}, we define the group $G_W$ and a $G_W$-stable complementary $B$ to $W$. The particular features of the representation $\sigma^{\oplus 3}:\SL_2(\F_3)\to\mathrm{Aut}(B,0)$ yield that $B$ is isogenous to ${E_j}^6$. Let $g\in \pstab(W)\simeq \SL_2(\F_3)$ be the unique element of order $2$. Recall that the matrix of $g$ is given by the representation $\sigma$, so that $g|_B=-\id_B$.

\steppf{1}{If an element $h\in G_W$ fixes no point, then $h$ has even order.}{Indeed, by Lemma \ref{lem-generalGW} (6), either $hg$ fixes a point $\tau$, or $1$ and $-1$ appear as eigenvalues of $h$. Clearly, $h$ has even order in the second case. 
In the first case, $hg$ actually is in $\pstab(W+\tau)$, and Propositions \ref{prop-fixor6int}, \ref{prop-jun6cod6} yield that $\pstab(W+\tau)$ is isomorphic to $\Z_3$, $\Z_3\times\Z_3$, or $\SL_2(\F_3)$. So either $hg$ has order $3$, in which case $h$ has even order $6$, or $hg\in\pstab(W+\tau)\simeq \SL_2(\F_3)$ has order $2,4,$ or $6$. But then, $g\in\pstab(W+\tau)$ since $G_W$ contains no translation. So $h\in\pstab(W+\tau)$ fixes some points, a contradiction!}

\steppf{2}{If an element $h\in G_W$ has prime order $p$, then $p\in\{2,3\}$. Moreover, if $p=3$, $h$ is a junior element or has junior square.}{If $h$ fixes no point, by Step 1, we have $p=2$. If $h$ fixes a point, then by Proposition \ref{prop-fixor6int}, and using the fact that $B\sim {E_j}^6$, we have $p\in\{2,3\}$.

When $p=3$, we thus also have $h\in\pstab(W+\tau)$ for some $\tau\in A$. Apply Proposition \ref{prop-fixor6int} to $\pstab(W+\tau)$. Note that by Proposition \ref{prop-jun6cod6}, $\omega\id_B$ does not appear in $\rho(G_W)$, and as $g|_B=-\id_B$ does, $j\id_B$ does not. In particular, $\pstab(W+\tau)$ can not the group $\Z_3\times\Z_3$ with the representation given by Proposition \ref{prop-fixor6int}. In the remaining two possible cases given by that proposition, every order 3 element of $\pstab(W+\tau)$ is junior or has junior square, and so is $h$.}

\steppf{3}{Any $3$-Sylow subgroup $S$ of $G_W$ is isomorphic to $\Z_3$, generated by one junior element.}{Let $h\in S$ be a non-trivial element. It has odd order, hence it fixes a point by Step 1, and thus it has order 3 by Proposition \ref{prop-fixor6int}. By Step 2, it is thus junior or a square of a junior element.

Let $s\in Z(S)$ be non-trivial, hence again (the square of) a junior element. Let us show that $h\in\grp{s}$. As $h$ and $s$ commute, either they have the same eigenspace for the eigenvalue 1, in which case $h\in\grp{s}$ as wished, or $E_{s|_B}(1)$ and $E_{h|_B}(1)$ are in direct sum, in which case $j\id_B\in\grp{s|_B,h|_B}$, and so $\omega\id_B\in \rho(G_W)$, which contradicts Proposition \ref{prop-jun6cod6}. Hence, $h\in\grp{s}$ and thus $S=\grp{s}\simeq\Z_3$.}

\steppf{4}{If $S_2$, $S_3$ are $2$ and $3$-Sylow subgroups of $G_W$, then $G_W= S_2\rtimes S_3$.}{By Step 3, no two elements of $S_3$ are conjugated in $G_W$, so $N_{G_W}(S_3)=C_{G_W}(S_3)$, and by Burnside's normal complement theorem \cite[10.1.8]{Robinson}, there is a normal subgroup $N\triangleleft G_W$ such that $G_W=N\rtimes S_3$. By Step 2, $N$ is a $2$-group, and it is clearly maximal. As it is normal, it is the unique $2$-Sylow subgroup of $G$, so $N=S_2$.}

\steppf{5}{$S_2$ has order $2^9$.}{We first count the number of junior elements in $G_W$. By Lemma \ref{lem-generalGW} (9), every junior element in $G_W$ fixes at least one $2$-torsion point in $B$. Since it acts trivially on a $3$-dimensional translated abelian subvariety of $B$, it fixes precisely $2^6$ of the $2$-torsion points in $B$. Each $2$-torsion point $\tau$ in $B$ is besides fixed by the four junior elements of $\pstab(W+\tau)\simeq\SL_2(\F_3)$ (by Proposition \ref{prop-fixor6int} and since $g$ of order 2 belongs to $\pstab(W+\tau)$). Hence, there are $\frac{2^{12}\times 4}{2^{6}}=2^8$ junior elements in $G_W$.

Now, note that by Step 3, the number $n_3$ of $3$-Sylow subgroups of $G_W$ equals the number of junior elements in $G_W$. Hence, denoting by $S_3$ a $3$-Sylow subgroup of $G_W$,
$$3|S_2|=|G_W|=n_3|N_{G_W}(S_3)|=n_3|C_{G_W}(S_3)|=2^9\cdot 3,$$
since it is easily checked that $C_{G_W}(S_3)=\grp{g,S_3}\simeq\Z_6<\SL_2(\F_3)$.}

\steppf{6}{Denote by $m_2$, $m_4$ the number of elements of order 2 and 4 in $S_2$. Then $m_2=6\cdot 61+1$ and $m_4=144$.}
{We first describe the order and trace of elements $h\in S_2$ different from $\id_A$ and $g$. By Lemma \ref{lem-blrat}, since $B\sim {E_j}^6$, and by \cite[Prop.2.4]{Wash}, the characteristic polynomial of $\rho(h)=M(h|_B)$ satisfies
$$\chi_{\rho(h)}=(X-1)^{\alpha}(X+1)^{\beta}\Phi_4(X)^{\gamma}\Phi_8(X)^{\delta},$$
with $\alpha,\beta,\gamma,\delta\ge 0$, $\beta$ being even because of the determinant and $\alpha+\beta+2\gamma+4\delta=6$ because of the dimension.
Hence, $\alpha$ is even too. 
If $\alpha\beta = 0$, then by Lemma \ref{lem-generalGW}, there is $\tau\in A$ such that $h\in \pstab(W+\tau)\cup g\pstab(W+\tau)$, so by Proposition \ref{prop-fixor6int}, the only possibility for $h$ other than $\id$ and $g$ satisfies $\chi_{\rho(h)} = {\Phi_4}^3$, hence $\alpha=\beta=0$.
Else, $\alpha$ and $\beta$ are positive.
So, $(\alpha,\beta,\gamma,\delta)$ can be $(0,0,3,0)$,$(2,2,1,0)$,$(2,4,0,0)$, or $(4,2,0,0)$. 
In particular, $h$ has order $2$ or $4$, with order $4$ if and only if $\mathrm{Tr}(h|_B)=0$, and order 2 if and only if $\mathrm{Tr}(h|_B)\in\{-2,2\}$.

Decomposing the representation $\rho|_{S_2}$ into irreducible subrepresentations yields a splitting coefficient $u\in\N$ such that
$u|S_2|=72+4(m_2-1),$
where $m_2$ is the number of elements of order $2$ in $S_2$. 
Denoting by $m_4$ the number of elements of order $4$ in $S_2$ and using Step 5, we rewrite
$(u-4)\cdot 2^9+4m_4=64$. So $u\le 4$.

Note that $h\in G_W$ junior of order $3$ acts by conjugation on the set of elements of order $2$ of the normal subgroup $S_2$, and the only fixed point is the element $g\in C_{G_W}(\grp{h})$. Hence, $m_2-1$ is divisible by $3$. So $u$ is divisible by $3$, and thus $u=3$, and $m_2 = 6\cdot 61 +1$, and $m_4=144$.}

\steppf{7}{But $m_4\ge 6\cdot 2^6$, contradiction!}
{Let us show that the number of elements of $G_W$ of order $4$ fixing a point is exactly $6\cdot 2^6$. By Lemma \ref{lem-generalGW} (8), if $h\in G_W$ has order $4$ and fixes a point, then all its $2^6$ fixed points in $B$ are $2$-torsion points of $B$. Moreover, by Proposition \ref{prop-fixor6int}, for any $\tau\in B$ of $2$-torsion, $\pstab(W+\tau)\simeq\SL_2(\F_3)$ contains exactly six elements of order $4$. Hence the count of $\frac{2^{12}\cdot 6}{2^6}=6\cdot 2^6$ elements of order $4$ fixing a point in $G_W$.}

\noindent And with this contradiction ends the proof of Lemma \ref{lem-sl23}.}

\rem[rem-macaulay]{Local information would not have been enough to rule out $\SL_2(\F_3)$. Indeed, considering a simply-connected neighborhood $U\subset \C^6$ of $0$, which is stable by the action of $\rho^{\oplus 3}:SL_2(\F_3)\hookrightarrow SL_6(\Q[j])$, the quotient $U/SL_2(\F_3)$ admits a crepant resolution. Let us construct it. 

Under the action of $SL_2(\F_3)$ on $\C^6$, exactly four 3-dimensional linear subspaces $Z_1,Z_2,Z_3,Z_4$ have non-trivial point-wise stabilizers $\grp{g_1},\grp{g_2},\grp{g_3},\grp{g_4}\simeq\Z_3$, where $g_1,g_2,g_3,g_4$ are the four junior elements of $\SL_2(\F_3)$. Using \texttt{Macaulay2}, a quick computation shows that the blow-up:
$$\eps : B:=Bl_{{\cal I}_{Z_1}\cap {\cal I}_{Z_2}\cap {\cal I}_{Z_3}\cap {\cal I}_{Z_4}}(\C^6)\to \C^6$$ is a smooth quasiprojective variety with a four-dimensional central fiber $\eps^{-1}(0)$. In particular, $B$ contains exactly four prime exceptional divisors $E_i$ ($1\le i\le 4$), one above each $Z_i$.

By the universal property of the blow-up, the action of $\SL_2(\F_3)$ on $\C^6$ lifts to an action on $B$. The lifted automorphism $\tilde{g_i}$ fixes the exceptional divisor $E_i$ pointwise: hence, locally, for any $x\in B$, $\pstab(x)$ is generated by pseudoreflections. Hence by Chevalley-Shepherd-Todd theorem, the quotient $X:= B/\SL_2(\F_3)$ is smooth.

We are going to prove that the resolution $X\to\C^6/\SL_2(\F_3)$ is crepant. As $\SL_2(\F_3)\subset\GL_6(\C)$ has one conjugacy class of junior elements, by Theorem \ref{thm-IR}, there is exactly one crepant divisor above $\C^6/\SL(2,3)$: A smooth resolution must contain this crepant divisor, and is thus crepant if and only if it contains exactly one exceptional divisor. This is clearly the case for $X$, since the action of $Q_8\subset\SL_2(\F_3)$ on $B$ is transitive on the set of the four prime exceptional divisors $E_1,E_2,E_3,,E_4$.}

\section{Proof of Theorem \ref{theo-freecod3}}

This is now straightforward.

\demode{Theorem \ref{theo-freecod3}}{Let $X$ be a variety with $c_1(X)=0$, that is a resolution of a quotient $A/G$ where $A$ is an abelian variety, and $G$ is a finite group acting freely in codimension $3$ on $A$. By Lemma \ref{lem_bbc2}, there is a finite étale Galois cover $\tilde{X}$ of $X$ which writes as a product
$$\tilde{X} = B\times\prod_{i=1}^r Y_i,$$
where each $Y_i$ is a Calabi-Yau manifold that resolves a quotient $B_i/H_i$, where $B_i$ is an abelian variety and $H_i$ is a finite group acting freely  in codimension $2$ on $B_i$. The fact that $G$ acts freely in codimension 3 on $A$ yields that, for each index $i$, the group $H_i$ acts freely in codimension 3 on $B_i$.

Suppose that $r\ge 1$. Then, since each $Y_i$ satisfies $c_1(Y_i)=0$ and $\pi_1(Y_i)=\{1\}$, by Lemma \ref{lem-makeitglob}, there is a point $bi\in B_i$ such that $\pstab(b_i)$ is a non-trivial subgroup of $H_i$. By Proposition \ref{prop-point-wise-stabilizer}, this subgroup  $\pstab(b_i)$ contains a junior element. By Propositions \ref{prop-globfailjun}, \ref{prop-ruleout6cod4}, \ref{prop-jun6cod5}, and \ref{prop-jun6cod6}, this junior element has eigenvalue $1$ with multiplicity $\dim(B_i)-3$, i.e., it stabilizes a translated abelian subvariety of $B_i$ of codimension $3$, a contradiction to our freeness-in-codimension-3 assumption.

So $\tilde{X}=B$ is an abelian variety, so $X$ itself is a quotient of the form $B/{\rm Gal}$ by a finite group ${\rm Gal}$ acting freely on the abelian variety $B$. As quotient singularities are rational, the resolution map $\eps:X\to A/G$ is in fact an isomorphism, and thus $G$ acts freely on $A$, as wished.}

\section{Concluding the proof of Theorem \ref{theo-main}}\label{sec-pftheomain}

Let us assemble the parts of the previous sections to prove Theorem \ref{theo-main}.

\demode{Theorem \ref{theo-main}.}{Let $A$ be an abelian variety of dimension $n$, and let $G$ be a finite group acting freely in codimension $2$ on $A$, such that $A/G$ has a resolution $X$ that is a Calabi-Yau manifold. By Proposition \ref{prop-tbmade}, either $A$ is isogenous to ${E_j}^n$ and $G$ is generated by junior elements of order 3 and 6, or $A$ is isogenous to ${E_{u_7}}^n$ and $G$ is generated by junior elements of order 7. In particular, $G$ is generated by its elements admitting fixed points. Also note that $G$ contains no junior element of order $6$ by Propositions \ref{prop-ruleout6cod4}, \ref{prop-jun6cod5}, and \ref{prop-jun6cod6}.

Let us show that for any translated abelian subvariety $W\subset A$, the pointwise stabilizer $\pstab(W)$ is abelian. First, it is generated by junior elements by Proposition \ref{prop-point-wise-stabilizer}. Let $g,h$ be two junior elements in $\pstab(W)$. As $g$ and $h$ both fix abelian varieties (containing the origin $0$) of codimension 3, their intersection $W'$ has codimension $3,4,5,$ or $6$ in $A$. Now, by Proposition \ref{prop-fixor6}, $\pstab(W')$ is thus abelian, and therefore $g$ and $h$ commute. 

Moreover, any two junior elements $g$ and $h$ in $\pstab(W)$ have the same order ($3$ if $A\sim {E_j}^n$, $7$ if $A\sim {E_{u_7}}^n$). Hence, using the structure theorem for finite abelian groups,
$\pstab(W)$ is isomorphic to ${\Z_3}^k$ for some $k$ if $A\sim {E_j}^n$, to ${\Z_7}^k$ for some $k$ if $A\sim {E_{u_7}}^n$.
Finally, if $g,h\in\pstab(W)$ are junior elements, then their eigenspaces with eigenvalues other than 1 are in direct sum by Proposition \ref{prop-fixor6}. An induction using that all junior elements of $\pstab(W)$ are codiagonalizable then concludes the proof of Theorem \ref{theo-main}.
}

\section{Proof of Theorem \ref{theo-dim4}}\label{sec-dim4}

In this section, we proceed to the proof of Theorem \ref{theo-dim4}, which in fact splits into two pieces. The first piece describes a slight generalization of the situation in dimension $3$ \cite{III0}. It notably gives an alternative proof of \cite[Key Claim 2]{III0}, replacing the discussion on invariant cohomology and topological Euler characteristics inherent to \cite[§3]{III0} with group theory and a geometric fixed loci argument ruling out the special linear group $\SL_3(\F_2)$.

\theo[theo-dim3gen]{Let $A$ be an abelian variety on which a finite group $G$ acts freely in codimension $2$ without translations. Suppose that $A/G$ has a resolution $X$ which is a Calabi-Yau manifold. Then, for any two junior elements $g,h\in G$ such that $\grp{g}\ne\grp{h}$, the intersection of eigenspaces $E_{M(g)}(1)\cap E_{M(h)}(1)$ does not have codimension $3$ in $H^0(A,T_A)$.}


\theo[theo-dim4gen]{Let $A$ be an abelian variety on which a finite group $G$ acts freely in codimension $2$ without translations. Suppose that $A/G$ has a resolution $X$ which is a Calabi-Yau manifold. Then, for any two junior elements $g,h\in G$ such that $\grp{g}\ne\grp{h}$, the intersection of eigenspaces $E_{M(g)}(1)\cap E_{M(h)}(1)$ does not have codimension $4$ in $H^0(A,T_A)$.}

Let us first show how these two results imply Theorem \ref{theo-dim4}.

\demode{Theorem \ref{theo-dim4} when assuming Theorems \ref{theo-dim3gen}, \ref{theo-dim4gen}}{We argue by contradiction. Suppose that $A$ is an abelian fourfold, that $G$ is a finite group acting freely in codimension $2$, yet not freely, on $A$, and that $A/G$ admits a resolution $X$ with $c_1(X)=0$. By Lemma \ref{lem_bbc2}, $X$ has a finite étale cover $\tilde{X}$ that is one of the following:
\begin{itemize}
\item an abelian fourfold $B$: Then $X$ itself is a quotient of the form $B/{\rm Gal}$ by a finite group ${\rm Gal}$ acting freely on $B$. As quotient singularities are rational, the resolution map $\eps:X\to A/G$ is in fact an isomorphism, and thus $G$ acts freely on $A$, a contradiction to our assumptions;
\item a product of the form $E\times X_3$ or $E\times X_7$, as wished;
\item a Calabi--Yau fourfold resolving a quotient $B/H$ with $B$ an abelian fourfold, and $H$ a finite group acting freely in codimension $2$ on $B$.
\end{itemize}

From here on, we argue by contradiction, assuming that $\tilde{X}$ is a Calabi--Yau fourfold. Up to replacing $B$ by an isogenous variety, we can assume that $H$ contains no translation. 

If $H$ entails two junior elements $g,h$ such that $\grp{g}\ne\grp{h}$, then the intersection of their eigenspaces $E_{M(g)}(1)\cap E_{M(h)}(1)$ is an intersection of two lines in a four-dimensional vector space: It has dimension 0 or 1, hence codimension 3 or 4. This contradicts the conclusions of Theorems \ref{theo-dim3gen} and \ref{theo-dim4gen}.

So $H$ has all of its junior elements contained in $\grp{g}$, and thus by the first point of Theorem \ref{theo-main}, $H=\grp{g}$ and $g$ has order $3$ or $7$, and admits $1$ as an eigenvalue of multiplicity one. Up to conjugating the whole group $H$ by a translation, we can assume $g(0)=0$. Let $E\subset A$ be the elliptic curve containing $0$ and fixed pointwise by $g$, and $T$ be its $\grp{g}$-stable supplementary. The group $H$ acts diagonally on the product $E\times T$ by $\{\id_E\}\times\grp{g|_T}$, and the addition map $E\times T\to B$ is an $H$-equivariant isogeny by \cite[Theorem 13.2.8]{BirkLang}. The volume form on $E$ thus pulls back to a $G$-invariant $1$-form on $B$, and thus to a non-zero global holomorphic $1$-form on the Calabi-Yau resolution $\tilde{X}$ of $B/H$, a contradiction.}

\subsection{Proof of Theorem \ref{theo-dim3gen}}

By Theorem \ref{theo-main}, the proof reduces to the following two cases. The first one is simple.

\prop[prop-dim3gen3]{Let $A$ be an abelian variety isogenous to ${E_j}^n$. Let $g,h\in\mathrm{Aut}(A)$ be two junior elements of order $3$ such that $\grp{g,h}$ contains no translation and no non-junior element fixing points, and $E_{M(g)}(1)= E_{M(h)}(1)$. Then $g=h$.}

\demo{Recall that $M:\mathrm{Aut}(A)\to\mathrm{Aut}(A,0)$ which, to any automorphism of $A$,  associates its matrix, induces a representation of $\grp{g,h}$. As $\grp{g,h}$ contains no translation, $M$ is faithful. Applying Maschke's theorem to the invariant subspace $E_{M(g)}(1)= E_{M(h)}(1)$ in $H^0(T_A)$ yields an $\grp{M(g),M(h)}$-stable supplementary $S$ to it. Let $\rho$ be the faithful representation of $\grp{g,h}$ obtained by restricting $M$ to $S$. By the classification of junior elements in Proposition \ref{prop-junclass}, $\rho(g)=\rho(h)=j\id_S$. But $\rho$ is faithful, and thus $g=h$.}

The second case is the following result.

\prop[prop-dim3gen7]{Let $A$ be an abelian variety isogenous to ${E_{u_7}}^n$. Let $g,h\in\mathrm{Aut}(A)$ be two junior elements of order $7$ such that $\grp{g,h}$ contains no translation and no non-junior element fixing points, and $E_{M(g)}(1)= E_{M(h)}(1)$. Then $\grp{g}=\grp{h}$.}

Its proof relies on two lemmas.

\lem[lem-sl32appears]{Let $A$ be an abelian variety isogenous to ${E_{u_7}}^n$. Let $g,h\in\mathrm{Aut}(A)$ be two junior elements of order $7$ such that $\grp{g,h}$ contains no translation and no non-junior element fixing points, and $E_{M(g)}(1)= E_{M(h)}(1)$. Then $\grp{g,h}$ is isomorphic to $\Z_7$ or $\SL_3(\F_2)$.}

\demo{Up to conjugating both $g$ and $h$ by the same translation, we can assume $g(0)=0$. Let $W$ be the abelian variety fixed by $g$, and let $G:= \grp{g,h}$. Let $G_W$ be the subgroup of $G$ defined by Lemma \ref{lem-generalGW}. Note that the condition $E_{M(g)}(1)= E_{M(h)}(1)$ translates to the fact that $G_W=\grp{g,h}$.

Let $\rho$ be the induced faithful representation of $\grp{g,h}$ defined in Lemma \ref{lem-generalGW} (3).
Let $k\in \grp{g,h}$. If $k$ has a fixed point in $A$, then by assumption, $k$ is junior of order $7$. Else, $k$ has no fixed point in $A$, and thus by Lemma \ref{lem-generalGW} (5), $1$ is an eigenvalue of $\rho(k)$. In that case, since $\rho(k)$ also has determinant 1, and by Lemma \ref{lem-blrat} and \cite[Prop.2.4]{Wash}, 
the characteristic polynomial (defined over $\Q[u_7]$) of the three-by-three matrix $\rho(k)$ is one of the following:
$${\Phi_1}^3,\; \Phi_1{\Phi_2}^2,\;\Phi_1\Phi_3,\;\Phi_1\Phi_4,\;\Phi_1\Phi_6.$$
So, possible prime divisors of $|\grp{g,h}|$ belong to $\{2,3,7\}$.

Let us denote by $\chi$ the character of the representation $\rho$.

Let $S_2$ be a $2$-Sylow subgroup of $\grp{g,h}$. It inherits the restricted representation $\rho|_{S_2}$, with character $\chi|_{S_2}$, and splitting coefficient $v_2$. Since a non-trivial element $g$ of $S_2$ has characteristic polynomial $\Phi_1{\Phi_2}^2$ or $\Phi_1\Phi_4$, it holds $|\chi(g)|^2 = 1$, and thus
$$9 + |S_2|-1 = \langle\chi|_{S_2},\chi|_{S_2}\rangle = v_2|S_2|,$$
with $v_2\in \N$. In other words, $(v_2-1)|S_2|=8$, and thus $|S_2|$ divides $8$.
Let $S_3$, $S_7$ be $3$ and $7$-Sylow subgroups of $\grp{g,h}$: Similarly, we obtain $|S_3|=3$ and $|S_7|=7$. Hence, the order $|\grp{g,h}|$ is a divisor of $8\cdot 3\cdot 7 = 168$. A \texttt{GAP} search (see Appendix) through all groups of such order which have no element of order 12, 14, or 21, and which have an either trivial or non-cyclic $2$-Sylow subgroup \cite[10.1.9]{Robinson} yields three candidates: $\Z_7$, $\Z_7\rtimes\Z_3$, and $\SL_3(\F_2)$. We exclude the second candidate as it is not generated by its elements of order $7$.}

We exclude $\SL_3(\F_2)$ by a geometric argument.

\lem[lem-sl32ruleout]{Let $A$ be an abelian variety isogenous to ${E_{u_7}}^n$. Let $g,h\in\mathrm{Aut}(A)$ be two junior elements of order $7$ such that $\grp{g,h}$ contains no translation and no non-junior element fixing points, and $E_{M(g)}(1)= E_{M(h)}(1)$. Then $\grp{g,h}$ cannot be isomorphic to $\SL_3(\F_2)$.}

\demo{We argue by contradiction. The multiplication table of $\grp{g,h}\simeq \SL_3(\F_2)$ shows that any element $g$ of order $7$ satisfies
\begin{center}
$C_{\grp{g,h}}(\grp{g})=\grp{g}$ and $N_{\grp{g,h}}(\grp{g})/C_{\grp{g,h}}(\grp{g})\simeq\Z_3$.
\end{center} 
Take an element $k\in N_{\grp{g,h}}(\grp{g})$ of order 3. Denote by $W_1,\ldots,W_7$ the seven disjoint translated abelian subvarieties of codimension 3 in $A$ that $g$ fixes pointwise. Then
$$k\left(\bigsqcup_{i=1}^7 W_i\right)=\bigsqcup_{i=1}^7 W_i,$$
and since $3$ and $7$ are coprime, there is some $1\le i\le 7$ such that $k(W_i)=W_i$.
Up to conjugating the whole group $\grp{g,h}$, we can assume that $0\in W_i$. We apply Lemma \ref{lem-generalGW} (2) to $g$, noting that $W=W_i$ and $k\in\grp{g,h}<G_W$. It shows that for any $w\in W_i$, one has $k(w)=w+T(k)$, and $pr_{W_i}(T(k))=0$. As $k(W_i)=W_i$, we obtain $T(k)=0$, so $k$ has fixed points and order $3$. In particular, it is not a power of a junior element, contradiction.}

\demode{Proposition \ref{prop-dim3gen7}}{By Lemmas \ref{lem-sl32appears} and \ref{lem-sl32ruleout}, we have $\grp{g,h}\simeq\Z_7$. But $\Z_7$ has no proper subgroup, so $\grp{g}=\grp{h}$.}

\subsection{Proof of Theorem \ref{theo-dim4gen}}

By Theorem \ref{theo-main}, the proof reduces to the following two cases.

\prop[prop-dim4gen7]{Let $A$ be an abelian variety isogenous to ${E_{u_7}}^n$. Let $g,h\in\mathrm{Aut}(A)$ be two junior elements of order $7$ such that $\grp{g,h}$ contains no translation and no non-junior element fixing points. Then $E_{M(g)}(1)\cap E_{M(h)}(1)$ cannot have codimension $4$ in $H^0(T_A)$.}

\prop[prop-dim4gen3]{Let $A$ be an abelian variety isogenous to ${E_{j}}^n$. Let $g,h\in\mathrm{Aut}(A)$ be two junior elements of order $3$ such that $\grp{g,h}$ contains no translation and no non-junior element fixing points. Then $E_{M(g)}(1)\cap E_{M(h)}(1)$ cannot have codimension $4$ in $H^0(T_A)$.}

Both propositions are proved by classifying matrices of elements in $\grp{g,h}$, and using representation theory to infer contradictory properties of $\grp{g,h}$.
We start with one lemma used in the proof of Proposition \ref{prop-dim4gen7}.

\lem[lem-class7dim4]{Let $A$ be an abelian variety isogenous to ${E_{u_7}}^n$. Let $g,h\in\mathrm{Aut}(A)$ be two junior elements of order $7$ such that $\grp{g,h}$ contains no translation and no non-junior element fixing points, and $E_{M(g)}(1)\cap E_{M(h)}(1)$ has codimension $4$ in $H^0(T_A)$. Then for every $k\in\grp{g,h}$, the trace of $M(k)\oplus\overline{M(k)}$ is at least $2n-8$, and equals $2n-7$ if $k$ is junior of order $7$.}

\demo{Up to conjugating both $g$ and $h$ by the same translation, we can assume $g(0)=0$. Note that the linear automorphism $M(h)$ also fixes $0$, and let $W$ be the abelian variety fixed by $\grp{g,M(h)}$. It is a subvariety of codimension $4$ of $A$. Let $G:= \grp{g,h}$. Let $G_W$ be the subgroup of $G$ defined by Lemma \ref{lem-generalGW}. Note that the fact that $E_{M(g)}(1) \cap_{M(h)}(1) = H^0(T_W)\subset H^0(T_A)$ translates to the fact that $G_W=\grp{g,h}$.
Consider the faithful representation $\rho$ of $\grp{g,h}$ given by Lemma \ref{lem-generalGW}. It has rank $4$.

Let $k\in\grp{g,h}$. If $k$ has a fixed point in $A$, then $k$ is junior of order $7$, and it is clear from Proposition \ref{prop-junclass} that the trace of $M(k)\oplus\overline{M(k)}$ equals $2n-7$. 
Else, by Lemma \ref{lem-generalGW} (5), $1$ is an eigenvalue of the four-by-four matrix $\rho(k)$. Since $\rho(k)$ also has determinant 1, and by Lemma \ref{lem-blrat} and \cite[Prop.2.4]{Wash}, its characteristic polynomial (defined over $\Q[u_7]$) is one of the following:
$${\Phi_1}^4,\;
{\Phi_1}^2{\Phi_2}^2,\;
{\Phi_1}^2\Phi_3,\;
{\Phi_1}^2\Phi_4,\;
{\Phi_1}^2\Phi_6,\;
$$
(where the matrix $\rho(k)$ has trace $4, 0, 1, 2, 3$ respectively), or one of the following 
$$\Phi_1(X)(X^3-\overline{u_7}X^2+u_7X-1),\; \Phi_1(X)(X^3-u_7X^2+\overline{u_7}X-1),$$
(where the matrix $\rho(k)\oplus\overline{\rho(k)}$ has trace $2+u_7+\overline{u_7}=1$). The consequence is that $\rho(k)\oplus\overline{\rho(k)}$ always has non-negative trace, which concludes.}

From this lemma follows a reduction to codimension 3 that concludes the proof of Proposition \ref{prop-dim4gen7}.

\demode{Proposition \ref{prop-dim4gen7}}{Assume by contradiction that $E_{M(g)}(1)\cap E_{M(h)}(1)$ has codimension $4$ in $H^0(T_A)$. Denote by $M$ the usual faithful matrix representation of $\grp{g,h}$ into $\GL(H^0(T_A))$, by $\chi_{M,\grp{g,h}}$ its character, and by $\1$ both the trivial representation of $\grp{g,h}$ and its character. We have
$$\langle \chi_{M,\grp{g,h}}, {\bf 1}\rangle = \sum_{k\in \grp{g,h}}\mathrm{Tr}\, M(k) =\frac{1}{2}\sum_{k\in \grp{g,h}}\mathrm{Tr}\, M(k)+\mathrm{Tr}\, \overline{M(k)} > (n-4)|\grp{g,h}|,$$
by Lemma \ref{lem-class7dim4}, the inequality being strict since $\grp{g,h}$ contains at least one junior element of order $7$. Hence, ${\bf 1}$ has multiplicity at least $n-3$ as a subrepresentation of $M$, i.e., $E_{M(g)}(1)\cap E_{M(h)}(1)$ has codimension at most $3$ in $H^0(T_A)$.}

We now prove an auxiliary lemma for Proposition \ref{prop-dim4gen3}.

\lem[lem-gh3group]{Let $A$ be an abelian variety isogenous to ${E_j}^n$. Let $g,h\in\mathrm{Aut}(A)$ be two junior elements of order $3$ such that $\grp{g,h}$ contains no translation and no non-junior element fixing points, and $E_{M(g)}(1)\cap E_{M(h)}(1)$ has codimension 4 in $H^0(T_A)$. Then each non-trivial element of $\grp{g,h}$ has order $3$.}

\demo{Up to conjugating both $g$ and $h$ by the same translation, we can assume $g(0)=0$. Note that the linear automorphism $M(h)$ also fixes $0$, and let $W$ be the abelian variety fixed by $\grp{g,M(h)}$. It is a subvariety of codimension $4$ of $A$. Let $G:= \grp{g,h}$. Let $G_W$ be the subgroup of $G$ defined by Lemma \ref{lem-generalGW}. Again, we have $G_W=\grp{g,h}$.
Consider the faithful representation $\rho$ of $\grp{g,h}$ of rank four given by Lemma \ref{lem-generalGW}, with character $\chi$.

Let $k\in\grp{g,h}$. If $k$ has a fixed point in $A$, then $k$ is junior of order $3$. Otherwise, by Lemma \ref{lem-generalGW} (5), $1$ is an eigenvalue of $\rho(k)$. In that case, since the intersection $E_{\rho(g)}(j)\cap E_{\rho(h)}(j)$ has dimension $2$, one of the roots of unity $1,j,j^2$ must appear as an eigenvalue of multiplicity $2$ for $\rho(k)$. By Lemma \ref{lem-blrat}, \cite[Prop.2.4]{Wash}, and since $\rho(k)$ has determinant one, the characteristic polynomial of $\rho(k)$, defined over $\Q[j]$, is one of the following
$${\Phi_1}^4,\; {\Phi_1}^2{\Phi_2}^2,\;{\Phi_1}^2\Phi_3,\;{\Phi_1}^2\Phi_4,\;{\Phi_1}^2\Phi_6,\;\Phi_1(X)(X-j)^3,\; {\Phi_1}(X)(X-j^2)^3.$$
So the order of $k$ is $1$, $3$, or an even number. 

To conclude, it is enough to show that no element of $\grp{g,h}$ has order $2$. We prove it by contradiction: Suppose that $k\in \grp{g,h}$ is such that $\rho(k)$ is similar to $\diag(1,1,-1,-1)$. As the eigenspace $E_{\rho(g)}(j)$ has codimension 1, the matrix $\rho(gk)$ has $j$ and $-j$ as eigenvalues. In particular, it is not junior, and thus it fixes no point. But its characteristic polynomial is not one of the polynomials listed above either, a contradiction.}

\demode{Proposition \ref{prop-dim4gen3}.}{Arguing by contradiction, we assume that the subspace $E_{M(g)}(1)\cap E_{M(h)}(1)$ has codimension 4 in $H^0(T_A)$. By Lemma \ref{lem-gh3group}, the non-trivial four-by-four matrices in the group $\rho(\grp{g,h})$ are each similar to one of the following:
\begin{equation}\label{eq-type3}
\diag(1,j,j,j),\;\diag(1,j^2,j^2,j^2),\mbox{ or }\diag(1,1,j,j^2).
\end{equation}
Note in particular that $\diag(j,j,j^2,j^2)$ is not an option (since it does not have $1$ as an eigenvalue, it would come from an element $k\in\grp{g,h}$ which fixes a point, and yet it is not junior, a contradiction).

As $\grp{g,h}$ is a $3$-group, we can set $k\in Z(\grp{g,h})$ to be a central element of order $3$. Since $E_{M(g)}(1)\cap E_{M(h)}(1)$ has codimension 4 in $H^0(T_A)$, we know that $\grp{g}\cap\grp{h}=\{1\}$. Thus, up to renaming $g$ and $h$ into each other, we can assume $k\not\in\grp{g}$.
If $\rho(k)$ is similar to $\diag(1,j,j,j)$ or $\diag(1,j^2,j^2,j^2)$, then $\rho(gk)$ or $\rho(g^2k)$, respectively, has no 1 as an eigenvalue, which contradicts (\ref{eq-type3}).
Hence, $\rho(k)$ is similar to $\diag(1,1,j,j^2)$. Note that $E_{\rho(g)}(j)\cap E_{\rho(h)}(j)$ has dimension $2$, and that any element of $\rho(\grp{g,h})$ must be a homothethy on this 2-dimensional subspace. Given the matrix $\rho(k)$,   this means that $E_{\rho(g)}(j)\cap E_{\rho(h)}(j)=E_{\rho(k)}(1)$. Again, either $\rho(gk)$ or $\rho(g^2k)$ has no 1 as an eigenvalue, which contradicts (\ref{eq-type3}).
}



\newpage

\section*{Appendix}
\addcontentsline{toc}{section}{Appendix}

\subsection*{Groups of order dividing 240 with an automorphism of order 7}
\begin{lstlisting}[language = GAP, firstline=3]
order_list := [];
nb_groups_of_order_list := [];

for a in [0..4] do
   for b in [0..1] do
      for c in [0..1] do
         n := (2^a)*(3^b)*(5^c);
         Add(order_list, n);
         Add(nb_groups_of_order_list, NumberSmallGroups(n));
      od;
   od;
od;

have_aut7 := [];

for i in [1..Length(order_list)] do
   n := order_list[i];
   for v in [1..nb_groups_of_order_list[i]] do
      g := SmallGroup(n,v);
      s := SylowSubgroup(g,2);
      if StructureDescription(s) = "Q16" or StructureDescription(s) = "Q8" 
      or StructureDescription(s) = "C16" or StructureDescription(s) = "C8"
      or StructureDescription(s) = "C4" or StructureDescription(s) = "C2" 
      or StructureDescription(s) = "1" then
         h := AutomorphismGroup(g);
         if Order(h) mod 7 = 0 then
            Add(l,(n,v)); 
         fi;
      fi;
   od;
od;
\end{lstlisting}

\subsection*{Representations of $\Z_3\rtimes\Z_8$}
\begin{lstlisting}[language = GAP, firstline=3]
for v in [1..NumberSmallGroups(24)] do
   g := SmallGroup(24, v);
   if StructureDescription(g) = "C3 : C8" 
   then Add(groups_checked, v);
        tbl_conjcl := ConjugacyClasses(g);
        nb_conjcl := Size(tbl_conjcl);
#locating the unique element of order 2 
#among conjugacy classes of g
        index_2 := 0;
        for j in [1..nb_conjcl] do
            o := Order(Representative(tbl_conjcl[j]));
            if o = 2 then index_2 := j; fi;
        od;
#only keeping irreducible characters sending 
#the unique element of order 2 to -id
        T := Irr(g);
        Tbis := [];
        for k in [1..nb_conjcl] do
           if T[k][index_2] + T[k][1] = 0 
           then Add(Tbis,T[k]);
           fi;
        od;
        Print(Tbis);
        Print("\n\n");
   fi;
od;
\end{lstlisting}

\subsection*{Proposition \ref{prop-2sylowiscyc}: Five candidates for $F$}

\begin{lstlisting}[language = GAP, firstline=3]
order_list := [];
nb_groups_of_order_list := [];

for a in [3..4] do
   for b in [0..1] do
      for c in [0..1] do
         n := (2^a)*(3^b)*(5^c);
         Add(order_list, n);
         Add(nb_groups_of_order_list, NumberSmallGroups(n));
      od;
   od;
od;

right_sylows := [];
right_sylows_and_orders := [];

for i in [1..Length(order_list)] do
   n := order_list[i];
   for v in [1..nb_groups_of_order_list[i]] do
      g := SmallGroup(n,v);
      s := SylowSubgroup(g,2);
      if StructureDescription(s) = "Q16" or StructureDescription(s) = "Q8" then
         Add(right_sylows, [n,v]);
         Add(right_sylows_and_orders, [n,v]);
         tbl_conjcl := ConjugacyClassesByOrbits(g);
         nb_conjcl := Size(tbl_conjcl);
         remove_once_only := 0;
         is_15 := 0;
         is_20 := 0;
         is_24 := 0;
         for i in [1..nb_conjcl] do
            o := Order(Representative(tbl_conjcl[i]));
            if o = 15 then
               is_15 := 1;
            fi;
            if o = 20 then
               is_20 := 1;
            fi;
            if o = 24 then
               is_24 := 1;
            fi;         
            s := Size(tbl_conjcl[i]);
            if remove_once_only = 0 and 
            ((o = 2 and s > 1) or (o mod 60 = 0) or (o mod 40 = 0) 
            or (is_20 = 1 and o mod 15 = 0) or (is_24 = 1 and o mod 15 = 0) 
            or (is_15 = 1 and o mod 20 = 0) or (is_24 = 1 and o mod 20 = 0) 
            or (is_15 = 1 and o mod 24 = 0) or (is_20 = 1 and o mod 24 = 0)) then
               Remove(right_sylows_and_orders);
               remove_once_only := 1;
            fi;
         od;
         if remove_once_only = 0 and (1 - is_15)*(1 - is_20)*(1 - is_24) = 1 then
            Remove(right_sylows_and_orders);
            remove_once_only := 1;
         fi;
      fi;
   od;
od;
\end{lstlisting}

\subsection*{Proposition \ref{prop-2sylowiscyc}: Two candidates generated by elements of the right order}

\begin{lstlisting}[language = GAP, firstline=3]
testing := [[48,8],[48,27]];

for i in [1..2] do
g := SmallGroup(testing[i][1],testing[i][2]);
Print(StructureDescription(g));
tbl_conjcl := ConjugacyClasses(g);
nb_conjcl := Size(tbl_conjcl);
   nb_elts_order_24 := 0;
   for j in [1..nb_conjcl] do
      o := Order(Representative(tbl_conjcl[j]));
      s := Size(tbl_conjcl[j]);
      if o = 24 then
         nb_elts_order_24 := nb_elts_order_24 + s;
      fi;
   od;
Print("number of elements of order 24:");
Print(nb_elts_order_24);
Print(" ");
Print("\n");
od;

Print("\n\n");

testing := [[40,4],[40,11],[80,18]];

for i in [1..3] do
g := SmallGroup(testing[i][1],testing[i][2]);
Print(StructureDescription(g));
tbl_conjcl := ConjugacyClasses(g);
nb_conjcl := Size(tbl_conjcl);
   nb_elts_order_20 := 0;
   for j in [1..nb_conjcl] do
      o := Order(Representative(tbl_conjcl[j]));
      s := Size(tbl_conjcl[j]);
      if o = 20 then
         nb_elts_order_20 := nb_elts_order_20 + s;
      fi;
   od;
Print(" number of elements of order 20: ");
Print(nb_elts_order_20);
Print(" ");
Print("\n");
od;
\end{lstlisting}

\subsection*{Proposition \ref{prop-2sylowiscyc}: None admitting the right representation}

\begin{lstlisting}[language = GAP, firstline=3]
tables_char_irr := [[],[]];
indices_20 := [[],[]];
testing := [[40,11],[80,18]];

for i in [1..2] do
   g := SmallGroup(testing[i][1], testing[i][2]);
   tbl_conjcl := ConjugacyClasses(g);
   nb_conjcl := Size(tbl_conjcl);
#locating the unique element of order 2 
#among conjugacy classes of g
   index_2 := 0;
   for j in [1..nb_conjcl] do
      o := Order(Representative(tbl_conjcl[j]));
      if o = 2 then
         index_2 := j;
      fi;
   od;
#only keeping irreducible characters sending 
#the unique element of order 2 to -id
   T := Irr(g);
   Tbis := [];
   for k in [1..nb_conjcl] do
      if T[k][index_2] + T[k][1] = 0 then
         Add(Tbis,T[k]);
      fi;
   od;
   Add(tables_char_irr[i], Tbis);
   Print(StructureDescription(g));
   Print(" possible irreducible representations have characters: ");
   Print(tables_char_irr[i]);
   Print("\n\n");
od;
\end{lstlisting}

\subsection*{Pointwise stabilizers in codimension 6 as in Subsection \ref{subsec-pstabcod6}}

\begin{lstlisting}[language = GAP, firstline=3]
order_list := [];
nb_groups_of_order_list := [];

for a in [0,3] do
   for b in [0..3] do
      for c in [0..2] do
         if b <= 1 or c <= 1 then
            n := (2^a)*(3^b)*(7^c);
            Add(order_list, n);
            Add(nb_groups_of_order_list, NumberSmallGroups(n));
         fi;
      od;
   od;
od;

right_sylows := [];
right_sylows_and_orders := [];

for i in [1..Length(order_list)] do
   n := order_list[i];
   for v in [1..nb_groups_of_order_list[i]] do
      g := SmallGroup(n,v);
      s := SylowSubgroup(g,2);
      t := SylowSubgroup(g,3);
      u := SylowSubgroup(g,7);
      if (StructureDescription(s) = "1" or StructureDescription(s) = "Q8") 
      and (StructureDescription(t) = "1" or StructureDescription(t) = "C3" 
           or StructureDescription(t) = "C9" 
           or StructureDescription(t) = "C3 x C3" 
           or StructureDescription(t) = "C3 x C9")
      and (StructureDescription(u) = "1" or StructureDescription(u) = "C7" 
           or StructureDescription(u) = "C7 x C7")
      then Add(right_sylows, [n,v]);
           Add(right_sylows_and_orders, [n,v]);
#we now remove of the list right_sylows_and_orders candidates with elements 
#of inappropriate order 63, or with several elements of order 2
           tbl_conjcl := ConjugacyClassesByOrbits(g);
           nb_conjcl := Size(tbl_conjcl);
           remove_once_only := 0;
           for i in [1..nb_conjcl] do   
              o := Order(Representative(tbl_conjcl[i]));     
              s := Size(tbl_conjcl[i]);
              if remove_once_only = 0 and 
              ((o = 2 and s > 1) or (o mod 63 = 0)) 
              then Remove(right_sylows_and_orders);
                   remove_once_only := 1;
              fi;
           od;
      fi;
   od;
od;


describe := [];
for i in [1..Length(right_sylows_and_orders)] do
   g := SmallGroup(right_sylows_and_orders[i][1], right_sylows_and_orders[i][2]);
   Add(describe, StructureDescription(g));
   Print(StructureDescription(g));
   Print("\n\n");
od;
\end{lstlisting}

\subsection*{Groups of order dividing 168 as in Lemma \ref{lem-sl32appears}}

\begin{lstlisting}[language = GAP, firstline=3]
order_list := [];
nb_groups_of_order_list := [];

for a in [0..3] do
   for b in [0..1] do
      n := (2^a)*(3^b)*7;
      Add(order_list, n);
      Add(nb_groups_of_order_list, NumberSmallGroups(n));
   od;
od;

right_sylow := [];
right_sylow_description := [];
for i in [1..Length(order_list)] do
   n := order_list[i];
   for v in [1..nb_groups_of_order_list[i]] do
      g := SmallGroup(n,v);
      h := SylowSubgroup(g, 2);
      if StructureDescription(h) = "Q8" or StructureDescription(h) = "D8" 
      or StructureDescription(h) = "1"
      then Add(right_sylow, [n, v]); 
           Add(right_sylow_description, StructureDescription(g)); 
      fi;
   od;
od;

right_sylow_and_orders := [];
right_sylow_and_orders_description := [];
for element in right_sylow do
    n := element[1];
    v := element[2];
    g := SmallGroup(n,v);
   tbl_conjcl := ConjugacyClassesByOrbits(g);
   nb_conjcl := Size(tbl_conjcl);
   v_to_discard := 0;
   for i in [1..nb_conjcl] do
      o := Order(Representative(tbl_conjcl[i]));
      if (o = 14 or o = 21 or o = 12) and v_to_discard = 0 
      then v_to_discard := 1; 
      fi;
   od;
   if v_to_discard = 0 
   then Add(right_sylow_and_orders, [n, v]);
        Add(right_sylow_and_orders_description, StructureDescription(g));
   fi;
od;

Print(right_sylow_and_orders);
Print("\n");
Print(right_sylow_and_orders_description);
\end{lstlisting}

\newpage
\bibliographystyle{plain}
\bibliography{biblioTQ}

\end{document}